\documentclass[11pt]{article}

\usepackage{amsmath}
\usepackage{amssymb}
\usepackage{natbib}
\usepackage{array}
\usepackage{xcolor}
\usepackage{graphicx}
\usepackage{multirow}
\usepackage{hyperref}
\usepackage{longtable}
\usepackage{rotating}
\usepackage{tikz}
\usetikzlibrary{patterns}

\begin{document}

\title{A forward-looking matheuristic approach for the multi-period\\ 
  two-dimensional non-guillotine cutting stock problem\\ with usable leftovers\thanks{This work has been
    partially supported by FAPESP (grants 2013/07375-0, 2016/01860-1,
    and 2018/24293-0), FAPES (grant 116/2019), and CNPq (grants 146110/2013-7, 306083/2016-7 and
    302682/2019-8).}}

\author{
  Ernesto G. Birgin\thanks{Department of Computer Science, Institute of
    Mathematics and Statistics, University of S\~ao Paulo, Rua do
    Mat\~ao, 1010, Cidade Universit\'aria, 05508-090, S\~ao Paulo, SP,
    Brazil. e-mail: egbirgin@ime.usp.br. Corresponding author.}
  \and
  Oberlan C. Rom\~ao\thanks{Department of Computing and Electronics, University Center of Northern Esp\'{\i}rito Santo of the Federal University of Esp\'{\i}rito Santo, Rodovia BR 101 Norte, Km 60, Bairro Litor\^aneo, 29932-540, S\~ao Mateus, ES, Brazil. e-mail: oberlan.romao@ufes.br}
  \and
  D\'ebora P. Ronconi\thanks{Department of Production Engineering,
    Polytechnic School, University of S\~ao Paulo, Av. Luciano Gualberto, 1380,
    Cidade Universit\'aria, 05508-010, S\~ao Paulo SP, Brazil. e-mail:
    dronconi@usp.br} 
}

\date{January 28, 2022}

\maketitle

\begin{abstract}
In [E. G. Birgin, O. C. Rom\~ao, and D. P. Ronconi, The multi-period
  two-dimensional non-guillotine cutting stock problem with usable
  leftovers, \textit{International Transactions in Operational
    Research} 27(3), 1392--1418, 2020] the multi-period two-dimensional
non-guillotine cutting stock problem with usable leftovers was
introduced. At each decision instant, the problem consists in
determining a cutting pattern for a set of ordered items using a set
of objects that can be purchased or can be leftovers of previous
periods; the goal being the minimization of the overall cost of the
objects up to the considered time horizon. Among solutions with
minimum cost, a solution that maximizes the value of the leftovers at
the end of the considered horizon is sought. A forward-looking
matheuristic approach that applies to this problem is introduced in
the present work. At each decision instant, the objects and the
cutting pattern that will be used is determined, taking into account the
impact of this decision in future states of the system. More
specifically, for each potentially used object, an attempt is made to
estimate the utilization rate of its leftovers and thereby
determine whether the object should be used or not. The approach's
performance is compared to the performance of a myopic 
technique. Numerical experiments show the efficacy of the proposed
approach.\\

\noindent
\textbf{Key words:} Two-dimensional cutting stock with usable leftovers,
non-guillotine cutting and packing, multi-period scenario,
forward-looking or looking-ahead approach, matheuristic.
\end{abstract}

\section{Introduction}
\label{sec1}

In this paper, we consider the multi-period two-dimensional non-guillotined cutting stock problem with usable leftovers. In the problem, $P$ periods of time denoted by $[s-1,s]$ for $s=1,\dots,P$ are considered; period $[s-1,s]$ corresponding to $t_{s-1} \leq t \leq t_s$, where $t_0 < t_1 < \dots < t_P$ are given decision time instants. Small rectangular pieces of varying sizes (named items) can be ordered at any instant $t$ between $t_0$ and $t_{P-1}$. However, assuming the discrete time convention, if an item is ordered at an instant $t$ such that $t_{s-1} \leq t \leq t_s$ for some $s \in \{1,\dots,P-1\}$, then it is assumed the item was ordered at instant $t_s$. All items ordered at instant $t_s$ must be produced between $t_s$ and $t_{s+1}$ and delivered at instant $t_{s+1}$. Raw material is available in the form of large rectangular purchasable pieces (named purchasable objects) or as usable leftovers of previous periods, i.e. parts of objects purchased at previous periods that were not used to produce items. (Remains of the cutting process can be classified as usable leftovers or can be discarded as scrap. Usable leftovers will be formally defined in Section~\ref{sec2}, but roughly speaking they can not be very old and must satisfy size constraints.) At each instant $t_s$, ordered items are known and the problem consists in selecting objects to be purchased and existent leftovers to produce all ordered items. The cutting pattern of each object (leftover or purchased) must also be determined. The problem is said to be two-dimensional because it involves the width and the height of items and objects; while it is said to be non-guillotine because cuts are not restricted to be guillotine cuts. Objects as well as leftovers can produce new leftovers. The amount of leftovers in stock is maintained under control with a parameter $\xi \in \{0,1,\dots,P\}$ that determines that parts (leftovers, leftovers of leftovers, etc) of an object purchased at instant $t_s$ can only be used at instants $t_{s+1},\dots,t_{s+\xi}$. (If $\xi=0$, the problem has no leftovers at all; while, if $\xi=1$, leftovers can only be used in the period immediately following the period in which they were generated.) The goal is to minimize the overall cost of objects purchased to produce all orders from instant $t_0$ to instant $t_{P-1}$ and, among the minimum cost solutions, to choose one in which the value of the usable leftovers remaining at instant $t_P$ (end of the considered time horizon) is maximized.

In the current work, we propose a forward-looking matheuristic to solve medium- and large-sized instances of the problem described in the paragraph above. In a training phase, the method attempts to estimate the proportion of each generated usable leftover that will be effectively used to produce items ordered in forthcoming periods. With this information, at a given period, a more expensive object can be purchased if the estimated future use of its leftovers points to future savings. A subproblem is solved per period. The decision variables determine the objects the must be purchased, the leftovers from previous periods that will be used, and their cutting pattern. All ordered items must be produced; and the goal is to minimize an objective function that, by discounting the cost of leftovers that are assumed to be used in the near future to produce ordered items within the considered time horizon, minimize the effective cost of the raw material required to produced the period ordered items. The estimation of effective usage of leftovers being generated, that is required to estimate the actual cost of the raw material, constitutes the forward-looking ingredient of the method. At the end of each training cycle, the estimated utilization proportion of each leftover is compared with its actual utilization proportion, and the estimate is updated. The updating rule and the stopping criterion ensure that the number of training cycles is finite.

The proposed method is calibrated with the instances with four periods considered in~\cite{bromro}; and then evaluated on a new set of instances with four, eight, and twelve periods. The performance of the method is compared with a myopic approach on the new set of thirty instances with up to twelve periods. For the new (small) instances with four periods, an additional comparison with CPLEX is also presented. The myopic approach differs with the forward-looking approach only in the objective function being minimized at each period. While the forward-looking approach considers the possible future use of letfovers, the myopic approach greedily minimizes the cost of the objects necessary to produce the ordered items of the period. The problem includes a parameter that tells for how many periods, after being generated, a leftover is available for use. The larger the durability of the leftovers, the greater the opportunity for economy. Experiments show that the forward-looking approach outperforms the myopic approach by a large extent and that, the greater the number of periods or the larger the durability of usable leftovers, the greater the advantage.

The problem considered in the present work was proposed in~\cite{bromro}, where a mixed integer linear programming model was introduced and instances with up to four periods were solved using CPLEX. However, no solution method has yet been proposed to deal with larger instances of the problem. The single-period version of the problem was considered in~\cite{abmro2}, where a discussion related to alternative definitions of usable leftovers was presented. Several papers in the literature, many of them based on real-world applications, address the one-dimensional cutting stock problem with usable leftovers; see the pioneers’ works~\cite{rood1986,scheit1991} and the more recent works~\cite{cherri2013,cherri2014,poldi2010,tomat2017,marble2020,pipes2021,integrated2021}. On the other hand, only a few publications tackle the two-dimensional case considered in the present work.

In all publications dedicated to the one-dimensional problem mentioned in the previous paragraph, a multi-period scenario is considered and a single threshold determines whether a cutting pattern leftover is disposed of as trim-loss or is a usable leftover. In particular, \cite{tomat2017} focuses on determining the optimal amount of usable leftovers that should be kept in stock in order to make good use of the raw material and at the same time minimize the cost of stock handling. In~\cite{cherri2013}, a heuristic that prioritizes the use of leftovers in order to control their stock quantity is presented. A rolling horizon scheme for the same problem is proposed in~\cite{poldi2010}. The subproblem of each period is solved with a simplex method with column generation and different strategies are considered in order to obtain integer solutions through rounding. A survey that reviews published studies up to 2014 can be found in~\cite{cherri2014}. A recent work~\citep{integrated2021} integrates the problem with the lot-sizing problem. In the problem under consideration, it is possible to bring forward the production of items with known demand in a future period. A relax-and-fix approach is proposed that solves the subproblems with a simplex method with column generation. Other recent works present practical applications in the marble industry~\citep{marble2020} and in the use of leftover piping in construction~\citep{pipes2021}.

Exact and non-exact two- and three-stage two-dimensional cutting stock problems with leftovers are considered in~\cite{silva2010}. In the considered problem, a single item is cut from a raw material object at a time, through one or two guillotine cuts, generating zero, one, or two ``residual objects''. A MILP model that extends the one-cut model presented in~\cite{dyckhoff1981} for the one-dimensional cutting stock problem is introduced; and numerical experiments solving real-world instances of the furniture industry and instances from the literature are presented. MILP models are solved with CPLEX. On the one hand, the goal is minimizing the number of cuts. On the other hand, several extensions, such as minimizing the number of used raw material objects (that are all of the same type), minimizing the length of the cuts, minimizing waste, allowing rotations, and considering multiple type of objects are also considered. One of the extensions, that points to attributing a value to the leftovers, opens the possibility of embedding the considered problem in a multi-period framework, as its was later done by the same authors in~\cite{silva2014}. In~\cite{silva2014}, the problem is integrated with the lot-sizing problem with the aim of minimizing a total cost that includes material, waste and storage costs. In the problem under consideration, anticipating the production of items maximizes raw material utilization while incurring stock costs; and a balance between these conflicting objectives is sought by minimizing their pricing. Two MILP models that do not depend on cutting patterns generation and two heuristics based on the industrial practice are presented. In contrast to the problem considered in the present work, at each period, two-stage non-exact cutting patterns are generated. In a brief contribution~\citep{chen2015}, a single-period problem with three-stage cutting patterns is considered in which the leftovers consist of remnants of the first cutting stage, the objective being to minimize the difference between the object cost and the value of the usable leftovers generated. A real-world multi-period three-dimensional cutting problem related to the supply of steel blocks in the metalworking is considered in~\cite{viegas2016}. Since remnants from one period can be used to produce items ordered in future periods, the problem considers leftovers; the objective being to keep stock growth under control. For the problem at hand, constructive heuristic procedures are proposed.

The rest of this paper is organized as follows. Section~\ref{sec2}
provides a formal description of the multi-period two-dimensional non-guillotine cutting stock
problem with leftovers. Section~\ref{sec3} introduces the proposed
matheuristic with a looking-ahead feature. Section~\ref{sec4} presents
numerical experiments. Conclusions and lines for future research are
given in the last section.

\section{The multi-period two-dimensional non-guillotine cutting stock problem with leftovers}
\label{sec2}

\begin{figure}[ht!]
\begin{center}
\begin{tabular}{ccc}
\resizebox{0.35\textwidth}{!}{\begin{tikzpicture}[scale=0.25]
\path (0.0,0.0) coordinate (P0);
\path (26,0.0) coordinate (P1);
\path (0.0,21)coordinate (P2);
\path (26,21) coordinate (P3);
\draw [thick, fill=none](P0) rectangle (P3);
\coordinate (S1EE) at (0,15);
\coordinate (S1ED) at (16,21);
\draw [fill=black!10] (S1EE) rectangle (S1ED);
\coordinate (TXTS1) at (8,18);
\draw (TXTS1) node {Top leftover};
\coordinate (S2EE) at (16,0);
\coordinate (S2ED) at (26,21);
\draw [fill=black!10] (S2EE) rectangle (S2ED);
\coordinate (TXTS2) at (21,10);
\draw (TXTS2) node[rotate=90] {\begin{tabular}{c}Right-hand side\\ leftover\end{tabular}};
\coordinate (TXTS3) at (8,7.5);
\draw (TXTS3) node {Cutting area};
\draw (0,-0.1) node {$\phantom{.}$};
\end{tikzpicture}} & $\phantom{a}$ &
\resizebox{0.35\textwidth}{!}{\begin{tikzpicture}[scale=0.25]
\path (0.0,0.0) coordinate (P0);
\path (26,0.0) coordinate (P1);
\path (0.0,21)coordinate (P2);
\path (26,21) coordinate (P3);
\draw [thick, fill=none](P0) rectangle (P3);
\coordinate (S1EE) at (0,15);
\coordinate (S1ED) at (26,21);
\draw [fill=black!10] (S1EE) rectangle (S1ED);
\coordinate (TXTS1) at (13,18);
\draw (TXTS1) node {Top leftover};
\coordinate (S2EE) at (16,0);
\coordinate (S2ED) at (26,15);
\draw [fill=black!10] (S2EE) rectangle (S2ED);
\coordinate (TXTS2) at (21,7.5);
\draw (TXTS2) node[rotate=90] {\begin{tabular}{c}Right-hand side\\ leftover\end{tabular}};
\coordinate (TXTS3) at (8,7.5);
\draw (TXTS3) node {Cutting area};
\draw (0,-0.1) node {$\phantom{.}$};
\end{tikzpicture}} \\
(a) && (b) \\
\end{tabular}
\end{center}
\caption{Pictures (a) and~(b) illustrate the two possible ways in
  which two leftovers can be generated from an object by performing a
  vertical and a horizontal guillotine pre-cut. In case~(a), the
  vertical guillotine pre-cut is made first; while, in case~(b), the
  horizontal guillotine pre-cut is made first.}
\label{fig1}
\end{figure}
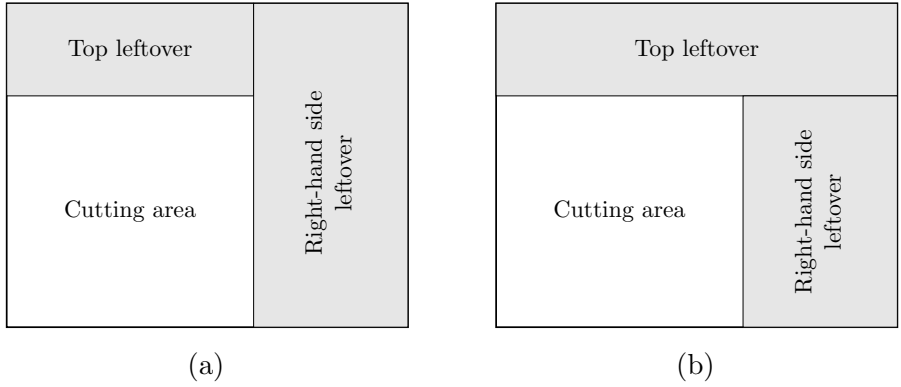

In this section, the multi-period two-di\-men\-sio\-nal non-guillotine
cutting stock problem with usable leftovers is described; and its mixed
integer linear programming formulation introduced in~\cite{bromro} is
presented. The (single-period) two-dimensional non-guillotine cutting
stock problem with leftovers was introduced in~\cite{abmro2} and
extended to the multi-period framework in~\cite{bromro}. One of the
main features of the problem is that, when an object is used to cut
items from it, two leftovers are obtained by performing a couple of
guillotine pre-cuts on the object that separate the leftovers from the
cutting area of the object (region from where the items will be
cut); see Figure~\ref{fig1}. Given a catalogue of items, we say a
leftover is usable if it can fit at least an item from the
catalogue. In this case, the leftover's value is given by its area
times the cost per unit of area of the object. Otherwise, the leftover
is disposable and has no value at all. It is worth noting that this
definition of leftovers implies that any part of the cutting area of
the object that is not used to produce an item is considered
waste. See~\cite{abm2016} and~\cite{abmro2} for other definitions of
leftovers in two-dimensional problems. \cite{abmro2} includes a detailed description of the single-period version of the problem, with several examples. Unlike the multi-period model presented in~\cite{bromro}, the model introduced in this section considers time instants~$s$ from~$p$ to~$P$. The possibility of choosing the initial and final instants of the model gives the necessary flexibility to formulate subproblems in algorithms of the rolling horizon type as the one that will be presented later.

\begin{figure}[ht!]
\begin{center}
\resizebox{0.8\textwidth}{!}{\definecolor{MyRed}{HTML}{F72C25} 
\definecolor{MyGreen}{RGB}{54,206,65}
\definecolor{MyBlue}{HTML}{0B88F8}

\begin{tikzpicture}[scale=0.35]
\draw[thick,dashed] (-3.5,10.0) -- (55,10.0);
\draw[thick,dashed] (-3.5,-1.5) -- (55,-1.5);
\draw[thick,dashed] (-3.5,-10.5) -- (55,-10.5);

\node[rotate=90] at (-2.9, 4.5) {Available objects};
\node[rotate=90] at (-2.9,-6.0) {Ordered items};

\draw[thick,dashed] (-2,-10.5) -- (-2,11.5);


\node[rotate=0] at (8, 10.8) {Instant $s=0$};


\node[rotate=0] at (4.0,8.9) {${\cal O}_{01}$};
\node[rotate=0] at (4.0,4.0) {10$\times$8};
\draw [thick](-1,0) rectangle (9,8);
\node[rotate=0] at (13,6.9) {${\cal O}_{02}$};
\node[rotate=0] at (13,3.0) {6$\times$6};
\draw [thick](10,0) rectangle (16,6);


\draw [fill=MyRed, thin](3.50,-8.00) rectangle (6.50,-4.00); 
\node[rotate=0] at (5.0,-9.5) {${\cal I}_{01}$};
\node[rotate=0] at (5.00,-6.00) {3$\times$4};

\draw [fill=MyRed, thin](8.50,-8.00) rectangle (11.50,-7.00); 
\node[rotate=0] at (10.0,-9.5) {${\cal I}_{02}$};
\node[rotate=0] at (10.00,-7.50) {3$\times$1};

\draw[thick,dashed] (17,-10.5) -- (17,11.5);


\node[rotate=0] at (27, 10.8) {Instant $s=1$};
%
\node[rotate=0] at (23.0,8.9) {${\cal O}_{11}$};
\node[rotate=0] at (23.0,4.0) {10$\times$8};
\draw [thick](18,0) rectangle (28,8);
\node[rotate=0] at (32,6.9) {${\cal O}_{12}$};
\node[rotate=0] at (32,3.0) {6$\times$6};
\draw [thick](29,0) rectangle (35,6);
%
%
\draw [fill=MyGreen, thin](20.5,-8.00) rectangle (24.5,-6.00); 
\node[rotate=0] at (22.5,-9.5) {${\cal I}_{11}$};
\node[rotate=0] at (22.5,-7.0) {4$\times$2};
\draw [fill=MyGreen, thin](25.5,-8.0) rectangle (28.5,-4.0); 
\node[rotate=0] at (27,-9.5) {${\cal I}_{12}$};
\node[rotate=0] at (27,-6.0) {3$\times$4};
\draw [fill=MyGreen, thin](29.5,-8.0) rectangle (32.5,-7.0); 
\node[rotate=0] at (31,-9.5) {${\cal I}_{13}$};
\node[rotate=0] at (31,-7.5) {3$\times$1};
\draw [fill=MyGreen, thin](29.5,-6.0) rectangle (32.5,-5.0); 
\node[rotate=0] at (31,-3.5) {${\cal I}_{14}$};
\node[rotate=0] at (31,-5.5) {3$\times$1};
\draw[thick,dashed] (36,-10.5) -- (36,11.5);
%
%
\node[rotate=0] at (46, 10.8) {Instant $s=2$};
%
%
\node[rotate=0] at (42.0,8.9) {${\cal O}_{21}$};
\node[rotate=0] at (42.0,4.0) {10$\times$8};
\draw [thick](37,0) rectangle (47,8);
\node[rotate=0] at (51,6.9) {${\cal O}_{22}$};
\node[rotate=0] at (51,3.0) {6$\times$6};
\draw [thick](48,0) rectangle (54,6);
%
%
\draw [fill=MyBlue, thin](41,-8.0) rectangle (45,-6.0); 
\node[rotate=0] at (43,-9.5) {${\cal I}_{22}$};
\node[rotate=0] at (43.00,-7.0) {4$\times$2};
\draw [fill=MyBlue, thin](41,-5.0) rectangle (45,-3.0); 
\node[rotate=0] at (39,-4) {${\cal I}_{21}$};
\node[rotate=0] at (43.00,-4.0) {4$\times$2};
\draw [fill=MyBlue, thin](47,-8.0) rectangle (50,-4.0); 
\node[rotate=0] at (48.5,-9.5) {${\cal I}_{23}$};
\node[rotate=0] at (48.50,-6.0) {3$\times$4};
\draw[thick,dashed] (55,-10.5) -- (55,11.5);

\end{tikzpicture}}
\end{center}
\caption{Illustration of a small instance with $p=0$, $P=3$, and $\xi=P-p=3$, meaning that
  usable leftovers generated at any period remain usable up to
  instant~$P$. The picture shows the available purchasable objects and
  the ordered items at each instant~$s \in \{0,1,2\}$. The numbers of
  available purchasable objects and ordered items at each instant are
  given by $m_0=m_1=m_2=2$ and $n_0=2$, $n_1=4$ and $n_2=3$,
  respectively. The cost per unit of area of all the objects is one
  (i.e.\ $c_{01} = c_{02} = c_{11} = c_{12} = c_{21} = c_{22} = 1$)
  and the catalogue with $d=1$ item is composed by an item with $\bar
  w_1 = 3$ and $\bar h_1 = 1$.}
\label{fig2}
\end{figure}
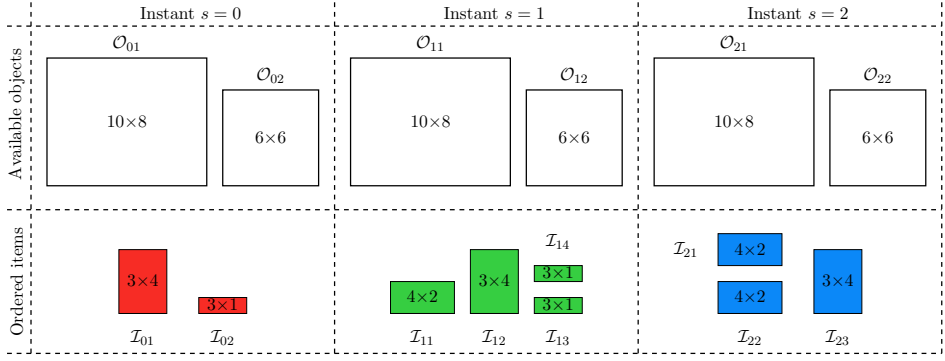

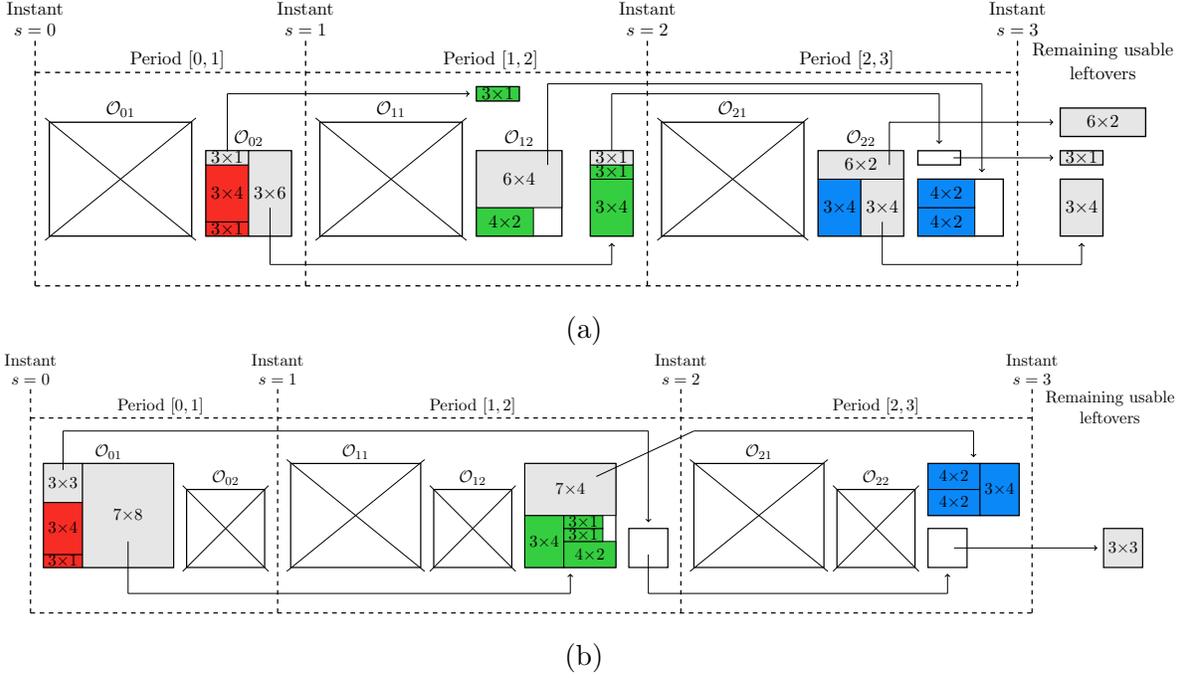
\begin{figure}[ht!]
\begin{center}
\begin{tabular}{c}
\resizebox{\textwidth}{!}{\definecolor{MyRed}{HTML}{F72C25} 
\definecolor{MyGreen}{RGB}{54,206,65}
\definecolor{MyBlue}{HTML}{0B88F8}

\begin{tikzpicture}[scale=0.3]

\draw[thick,dashed] (-2.0,11.5) -- (67,11.5);
\draw[thick,dashed] (-2.0,-3.5) -- (67,-3.5);

\draw[thick,dashed] (-2.0,-3.5) -- (-2.0,14);
\node[rotate=0] at (-2, 16) {Instant};
\node[rotate=0] at (-2, 14.5) {$s=0$};


\node[rotate=0] at (8, 12.4) {Period $[0,1]$};


\node[rotate=0] at (4.0,8.9) {${\cal O}_{01}$};
\draw [thick](-1,0) rectangle (9,8);
\node[rotate=0] at (13,6.9) {${\cal O}_{02}$};
\draw [thick](10,0) rectangle (16,6);

\draw [fill=MyRed, thin](10.00,0.00) rectangle (13.00,1.00); 
\node[rotate=0] at (11.5,0.50) {3$\times$1};
\draw [fill=MyRed, thin](10.00,1.00) rectangle (13.00,5.00); 
\node[rotate=0] at (11.5,3.0) {3$\times$4};

\draw [fill=black!10] (10.00,5.00) rectangle (13.00,6.00); 
\node[rotate=0] at (11.5,5.50) {3$\times$1};
\draw [fill=black!10] (13.00,0.00) rectangle (16.00,6.00); 
\node[rotate=0] at (14.5,3.0) {3$\times$6};

\draw[thin] (-1.3,8.3) -- (9.3,-0.3);
\draw[thin] (9.3,8.3) -- (-1.3,-0.3);

\node[rotate=0] at (11.50,5.5) (St02) {};
\node[rotate=0] at (11.50,5.5) (Sr02) {};

\draw[thick,dashed] (17.0,-3.5) -- (17.0,14);
\node[rotate=0] at (17.0, 16) {Instant};
\node[rotate=0] at (17.0, 14.5) {$s=1$};


\node[rotate=0] at (30, 12.4) {Period $[1,2]$};


\node[rotate=0] at (23.0,8.9) {${\cal O}_{11}$};
\draw [thick](18,0) rectangle (28,8);
\node[rotate=0] at (32,6.9) {${\cal O}_{12}$};
\draw [thick](29,0) rectangle (35,6);


\draw [thick](37,0.00) rectangle (40,6.00); 
\draw [thick](29,9.50) rectangle (32,10.50); 


\draw [fill=MyGreen, thin](37,0.00) rectangle (40,4.00); 
\node[rotate=0] at (38.50,2.00) {3$\times$4};
\draw [fill=MyGreen, thin](37.0,4.00) rectangle (40,5.0); 
\node[rotate=0] at (38.5,4.5) {3$\times$1};
\draw [fill=MyGreen, thin](29,9.50) rectangle (32,10.50); 
\node[rotate=0] at (30.5,10) {3$\times$1};
\draw [fill=MyGreen, thin](29.0,0.00) rectangle (33,2.0); 
\node[rotate=0] at (31.0,1.0) {4$\times$2};


\node[rotate=0] at (15.5,6.00) (St21) {};

\draw [fill=black!10] (37,5.00) rectangle (40,6.00); 
\node[rotate=0] at (38.5,5.5) {3$\times$1};
\draw [fill=black!10] (29,2.00) rectangle (35,6.00); 
\node[rotate=0] at (32,4.0) {6$\times$4};

\draw[thin] (17.7,8.3) -- (28.3,-0.3);
\draw[thin] (28.3,8.3) -- (17.7,-0.3);

\draw[-, thin] (St02) to (11.5,10);
\draw[->, thin] (11.5,10) to (28.5,10);
\draw[-, thin] (14.5,2.0) to (14.5,-2.0);
\draw[-, thin] (14.5,-2.0) to (38.5,-2.0);
\draw[->, thin] (38.5,-2.0) to (38.5,-0.5);

\draw[thick,dashed] (41,-3.5) -- (41,14);
\node[rotate=0] at (41, 16) {Instant};
\node[rotate=0] at (41, 14.5) {$s=2$};
%

\node[rotate=0] at (55, 12.4) {Period $[2,3]$};


\node[rotate=0] at (47.0,8.9) {${\cal O}_{21}$};
\draw [thick](42,0) rectangle (52,8);
\node[rotate=0] at (56,6.9) {${\cal O}_{22}$};
\draw [thick](53,0) rectangle (59,6);

\draw[thin] (41.7,8.3) -- (52.3,-0.3);
\draw[thin] (52.3,8.3) -- (41.7,-0.3);


\draw [thick](60,0.00) rectangle (66,4.00); 
\draw [thick](60,5.00) rectangle (63,6.00);

\draw [fill=MyBlue, thin](60,0.00) rectangle (64,2.00); 
\node[rotate=0] at (62.00,1.0) {4$\times$2};
\draw [fill=MyBlue, thin](60,2.00) rectangle (64,4.00); 
\node[rotate=0] at (62.00,3.0) {4$\times$2};
\draw [fill=MyBlue, thin](53,0.00) rectangle (56,4.00); 
\node[rotate=0] at (54.5,2.0) {3$\times$4};

\draw [fill=black!10] (56,0.00) rectangle (59,4.00); 
\node[rotate=0] at (57.5,2.0) {3$\times$4};
\draw [fill=black!10] (53,4.00) rectangle (59,6.00); 
\node[rotate=0] at (56,5.0) {6$\times$2};
%
\draw[-, thin] (38.5,6.0) to (38.5,10.0);
\draw[-, thin] (38.5,10.0) to (61.5,10.0);
\draw[->, thin] (61.5,10.0) to (61.5,6.5);

\draw[-, thin] (34,5.0) to (34,10.7);
\draw[-, thin] (34,10.7) to (64.5,10.7);
\draw[->, thin] (64.5,10.7) to (64.5,4.5);

\draw[thick,dashed] (67,-3.5) -- (67,14);
\node[rotate=0] at (67, 16) {Instant};
\node[rotate=0] at (67, 14.5) {$s=3$};

\node[rotate=0] at (73, 13) {Remaining usable};
\node[rotate=0] at (73, 11.5) {leftovers};
%
%
\draw [thick,fill=black!10](70,7) rectangle (76,9);
\node[rotate=0] at (73,8) {6$\times$2};

\draw [thick,fill=black!10](70,5) rectangle (73,6);
\node[rotate=0] at (71.5,5.5) {3$\times$1};

\draw [thick,fill=black!10](70,0) rectangle (73,4);
\node[rotate=0] at (71.5,2.0) {3$\times$4};

\draw[->, thin] (62.5,5.50) to (69.5,5.50);

\draw[-, thin] (57.5,1.0) to (57.5,-2.0);
\draw[-, thin] (57.5,-2.0) to (71.5,-2.0);
\draw[->, thin] (71.5,-2.0) to (71.5,-0.5);

\draw[-, thin] (58,5.0) to (58,8.0);
\draw[-, thin] (58,8.0) to (58,8.0);
\draw[->, thin] (58,8.0) to (69.5,8.0);

\end{tikzpicture}}\\[2mm]
(a) \\
\resizebox{\textwidth}{!}{\definecolor{MyRed}{HTML}{F72C25} 
\definecolor{MyGreen}{RGB}{54,206,65}
\definecolor{MyBlue}{HTML}{0B88F8}

\begin{tikzpicture}[scale=0.3]

\draw[thick,dashed] (-2.0,11.5) -- (75,11.5);
\draw[thick,dashed] (-2.0,-3.5) -- (75,-3.5);

\draw[thick,dashed] (-2.0,-3.5) -- (-2.0,14);
\node[rotate=0] at (-2, 16) {Instant};
\node[rotate=0] at (-2, 14.5) {$s=0$};


\node[rotate=0] at (8, 12.4) {Period $[0,1]$};


\node[rotate=0] at (4.0,8.9) {${\cal O}_{01}$};
\draw [thick](-1,0) rectangle (9,8);
\node[rotate=0] at (13,6.9) {${\cal O}_{02}$};
\draw [thick](10,0) rectangle (16,6);

\draw [fill=MyRed, thin](-1.00,0.00) rectangle (2.00,1.00); 
\node[rotate=0] at (0.5,0.50) {3$\times$1};
\draw [fill=MyRed, thin](-1.00,1.00) rectangle (2.00,5.00); 
\node[rotate=0] at (0.5,3.0) {3$\times$4};

\draw [fill=black!10] (-1.00,5.00) rectangle (2.00,8.00); 
\node[rotate=0] at (0.5,6.50) {3$\times$3};
\draw [fill=black!10] (2.00,0.00) rectangle (9.00,8.00); 
\node[rotate=0] at (5.5,4.0) {7$\times$8};

\draw[thin] (9.7,6.3) -- (16.3,-0.3);
\draw[thin] (16.3,6.3) -- (9.7,-0.3);

\draw[thick,dashed] (17.0,-3.5) -- (17.0,14);
\node[rotate=0] at (17.0, 16) {Instant};
\node[rotate=0] at (17.0, 14.5) {$s=1$};


\node[rotate=0] at (32, 12.4) {Period $[1,2]$};


\node[rotate=0] at (23.0,8.9) {${\cal O}_{11}$};
\draw [thick](18,0) rectangle (28,8);
\node[rotate=0] at (32,6.9) {${\cal O}_{12}$};
\draw [thick](29,0) rectangle (35,6);


\draw [thick](44,0.00) rectangle (47,3.00); 
\draw [thick](36,0.00) rectangle (43,8.00); 


\draw [fill=MyGreen, thin](36,0.00) rectangle (39,4.00); 
\node[rotate=0] at (37.50,2.00) {3$\times$4};
\draw [fill=MyGreen, thin](39.0,0.00) rectangle (43,2.0); 
\node[rotate=0] at (41.0,1.0) {4$\times$2};
\draw [fill=MyGreen, thin](39.0,2.00) rectangle (42,3.0); 
\node[rotate=0] at (40.5,2.5) {3$\times$1};
\draw [fill=MyGreen, thin](39.0,3.00) rectangle (42,4.0); 
\node[rotate=0] at (40.5,3.5) {3$\times$1};


\draw [fill=black!10] (36,4.00) rectangle (43,8.00); 
\node[rotate=0] at (39.5,6.0) {7$\times$4};

\draw[thin] (17.7,8.3) -- (28.3,-0.3);
\draw[thin] (28.3,8.3) -- (17.7,-0.3);

\draw[thin] (28.7,6.3) -- (35.3,-0.3);
\draw[thin] (35.3,6.3) -- (28.7,-0.3);

\draw[-, thin] (0.5,7.50) to (0.5,10.5);
\draw[-, thin] (0.5,10.5) to (45.5,10.5);
\draw[->, thin] (45.5,10.5) to (45.5,3.5);
\draw[-, thin] (5.5,2.0) to (5.5,-2.0);
\draw[-, thin] (5.5,-2.0) to (39.5,-2.0);
\draw[->, thin] (39.5,-2.0) to (39.5,-0.5);

\draw[thick,dashed] (48,-3.5) -- (48,14);
\node[rotate=0] at (48, 16) {Instant};
\node[rotate=0] at (48, 14.5) {$s=2$};


\node[rotate=0] at (63, 12.4) {Period $[2,3]$};


\node[rotate=0] at (54.0,8.9) {${\cal O}_{21}$};
\draw [thick](49,0) rectangle (59,8);
\node[rotate=0] at (63,6.9) {${\cal O}_{22}$};
\draw [thick](60,0) rectangle (66,6);

\draw[thin] (48.7,8.3) -- (59.3,-0.3);
\draw[thin] (59.3,8.3) -- (48.7,-0.3);

\draw[thin] (59.7,6.3) -- (66.3,-0.3);
\draw[thin] (66.3,6.3) -- (59.7,-0.3);


\draw [thick](67,0.00) rectangle (70,3.00); 
\draw [thick](67,4.00) rectangle (74,8.00);


\draw [fill=MyBlue, thin](67,4.00) rectangle (71,6.00); 
\node[rotate=0] at (69.00,5.0) {4$\times$2};
\draw [fill=MyBlue, thin](67,6.00) rectangle (71,8.00); 
\node[rotate=0] at (69.00,7.0) {4$\times$2};
\draw [fill=MyBlue, thin](71,4.00) rectangle (74,8.00); 
\node[rotate=0] at (72.5,6.0) {3$\times$4};
%
%

\draw[-, thin] (41.5,7.0) to (49,10.5);
\draw[-, thin] (49,10.5) to (70.5,10.5);
\draw[->, thin] (70.5,10.5) to (70.5,8.5);
\draw[-, thin] (45.5,1.00) to (45.5,-2.00);
\draw[-, thin] (45.5,-2.00) to (68.5,-2.00);
\draw[->, thin] (68.5,-2.00) to (68.5,-0.5);

\draw[thick,dashed] (75,-3.5) -- (75,14);
\node[rotate=0] at (75, 16) {Instant};
\node[rotate=0] at (75, 14.5) {$s=3$};

\node[rotate=0] at (81, 13) {Remaining usable};
\node[rotate=0] at (81, 11.5) {leftovers};
%
%

\draw [thick,fill=black!10](80.5,0) rectangle (83.5,3);
\node[rotate=0] at (82,1.5) {3$\times$3};

\draw[->, thin] (69,1.5) to (80.0,1.5);

\end{tikzpicture}}\\[2mm]
(b)
\end{tabular}
\end{center}
\caption{Illustration of two solutions that, at each period, may cut
  ordered items from purchasable objects or from usable leftovers from
  previous periods. (a) Greedy solution obtained by a myopic method
  that, \textit{at each decision instant}, minimizes the cost of the
  purchasable objects required to cut the ordered items of that
  instant, assuming that usable leftovers from previous periods are
  free. (b) Solution with minimum total cost of the required
  purchasable objects and, in addition, maximum value of the usable
  leftovers at instant $P=3$. The cost of the purchased objects in the
  solution in (a) is 108; while the same cost is 80 in (b).}
\label{fig3}
\end{figure}

Let $p$ and $P$ satisfying $p < P$ be the first and the last instant
to be considered, respectively. For each instant $s=p,\dots,P-1$,
there are given~$m_s$ purchasable objects~${\cal O}_{sj}$ with
width~$W_{sj}$, height~$H_{sj}$, and cost~$c_{sj}$ per unit of area
($j=1,\dots,m_s$) and a set of $n_s$ ordered items ${\cal I}_{si}$
with width~$w_{si}$ and height~$h_{si}$ ($i=1,\dots,n_s$). A catalogue
composed by~$d$ items $\bar {\cal I}_i$ with width~$\bar w_i$ and
height~$\bar h_i$ ($i=1,\dots,d$) is also given. A parameter $\xi \in
[0,P-p]$ says that leftovers generated within a period~$[s,s+1)$
  remain valid up to period $[s+\xi,s+\xi+1)$. By definition, each
    object generates two leftovers. This means that the number of
    objects at instant~$s$ is given by
\begin{equation} \label{mpbar}
\bar m_s = m_s + 2 \, \hat m_{s-1} \mbox{ for } s = p,\dots,P,
\end{equation}
where 
\begin{equation} \label{mphat}
\hat m_s = \sum_{\ell=0}^{\min\{s-p,\xi-1\}} 2^{\ell} m_{s-\ell},
\mbox{ for } s=p,\dots,P-1,
\end{equation}
stands for the number of objects that, at period $[s,s+1)$, generate
leftovers, $\hat m_{p-1}=0$ (i.e.\ no leftovers coming from previous
periods at the first considered instant $s=p$), and $m_P=0$ (i.e.\ no
purchasable objects at the last considered instant $s=P$).  Note that,
since, by definition, there are no purchasable objects at instant~$P$,
$\bar m_P$ represents the \textit{number of leftovers} available at
instant~$P$. The problem consists in minimizing the overall cost of
the purchasable objects required to produce the items ordered at
instants $p,\dots,P-1$ making use of leftovers; and, among all
solutions with minimum cost, maximizing the value of the usable
leftovers at instant~$P$. See Figures~\ref{fig2} and~\ref{fig3}. Figure~\ref{fig2} describes a toy instance of the problem; while Figure~\ref{fig3} exhibits two different feasible solutions.

Purchasable objects~${\cal O}_{sj}$ ($s=p,\dots,P-1$, $j=1,\dots,m_s$)
have a given cost~$c_{sj}$ per unit of area. The value of an usable
leftover is given by its area times its cost per unit of area; and the
cost per unit of area of a leftover corresponds to the cost per unit
of area of the purchasable object from which the leftover comes
from. In order to make this relation, we associate to each
(purchasable or leftover) object~${\cal O}_{sj}$ ($s=p,\dots,P$,
$j=1,\dots,\bar m_s$) an expiration date~$e_{sj}$ in such a way that,
if~${\cal O}_{sj}$ is a purchasable object, we define $e_{sj}=\xi$;
while if~${\cal O}_{sj}$ is a leftover then we define~$e_{sj}$ as the
expiration date of the object from which it comes from reduced by
one. Clearly, $e_{sj} \geq 0$, since objects with null expiration date
do not generate leftovers. Let $j_1^s \leq j_2^s \leq \dots \leq
j_{\hat m_s}^s$ be the indices of the $\hat m_s$ objects that generate
leftovers in the period $[s,s+1)$; and let us define that, at instant
  $s+1$, objects~${\cal O}_{s+1,m_{s+1}+2k-1}$ and~${\cal
    O}_{s+1,m_{s+1}+2k}$ correspond to the ``top leftover'' and to the
  ``right-hand-side leftover'' of object~${\cal O}_{s,j_k^s}$,
  respectively. Thus, $c_{s+1,m_{s+1}+2k-1} = c_{s+1,m_{s+1}+2k} =
  c_{s,j_k^s}$ and $e_{s+1,m_{s+1}+2k-1} = e_{s+1,m_{s+1}+2k} =
  e_{s,j_k^s} - 1$.  The relevant costs are the costs $c_{P,j}$
  ($j=m_P+1,\dots,\bar m_P$) that correspond to the value (per unit of
  area) of the leftovers available at instant~$P$, i.e. at the end of
  the considered time horizon, that are the leftovers whose value must
  be maximized. For a given instant $s$ ($s=p,\dots,P-1$) and the
  expiration dates $e_{sj}$ of the $\bar m_s$ objects available at the
  instant, the $\hat m_s \leq \bar m_s$ indices $j^s_1, j^s_2,\dots$
  of the objects that potentially generate leftovers can be computed
  as follows.  Start with $k = 0$ and, for $j$ from $1$ to $\bar m_s$,
  if $e_{sj}>0$ then increase $k$ by one and set $j_k^s=j$. Finish by
  setting $\hat m_s=k$.

The description of the problem's variables follows. Variables
$v_{sij}\in \{0,1\}$ ($s=p,\dots,P-1$, $j=1,\dots,\bar m_s$,
$i=1,\dots,n_s$) assign items to objects ($v_{sij}=1$ if item ${\cal
  I}_{si}$ is assigned to object ${\cal O}_{sj}$; and $v_{sij}=0$
otherwise). Variables $u_{sj} \in \{0,1\}$
($s=p,\dots,P-1$,$j=1,\dots,m_s$) identify whether at least an item is
assigned to object~${\cal O}_{sj}$ or not ($u_{sj}=1$ and $u_{sj}=0$,
respectively). Variables $\eta_{sj} \in \{0,1\}$
($s=p,\dots,P-1$,$j=1,\dots,\bar m_s$) determine if the vertical
pre-cut that separates the cutting area from the leftover in
object~${\cal O}_{sj}$ is made before the horizontal pre-cut
($\eta_{sj}=1$) or if the horizontal pre-cut precedes the vertical
pre-cut ($\eta_{sj}=0$). Variables $t_{sj}$ and $r_{sj} \in
\mathbb{R}$ ($s=p,\dots,P-1$, $j=1,\dots,\bar m_s$) determine the
height of the top leftover and the width of the right-hand-side
leftover of object~${\cal O}_{sj}$, respectively. Variables $\bar
W_{sj}$ and $\bar H_{sj} \in \mathbb{R}$ ($s=p,\dots,P$,
$j=1,\dots,\bar m_s$) represent the width and the height of
object~${\cal O}_{sj}$. (This is relevant to the objects that are
leftovers of objects purchased at previous periods, since the
dimensions of purchasable objects are constant, i.e.\ $\bar W_{sj} =
W_{sj}$ and $\bar H_{sj} = H_{sj}$ for every~$s$ whenever $1 \leq j
\leq m_s$.)  Variables $\pi_{sii'}$ and $\tau_{sii'} \in \{0,1\}$
($s=p,\dots,P-1$, $i=1,\dots,n_s$, $i'=i+1,\dots,n_s$) are auxiliary
variables used to avoid the overlapping between items. Variables
$\gamma_j \in \mathbb{R}$ ($j=1,\dots,\bar m_P$) are related to the
value of the area of the leftovers at instant~$P$, i.e.\ at the end of
the considered time horizon. Variables $\theta_{j\ell} \in \{0,1\}$
and $\omega_{j\ell} \in \mathbb{R}$ ($j=1,\dots,\bar m_P$,
$\ell=1,\dots,L$) are auxiliary variables used to linearize the
computation of these areas (product of the leftovers variable
dimensions), where $L= \lfloor \log_2(\hat W) \rfloor + 1$, $\hat W =
\max \{ W_{sj} \; | \; s=p,\dots,P-1, j=1,\dots,m_s \}$, and, for
further reference, $\hat H = \max \{ H_{sj} \; | \; s=p,\dots,P-1,
j=1,\dots,m_s \}$. The auxiliary variables $\zeta_{ji}\in \{0,1\}$
($j=1,\dots,\bar m_P$, $i=1,\dots,d$) are used to nullify the value of
the area of a leftover at instant~$P$ if it can not fit any item from
the catalogue.

The problem consists in minimizing
\begin{equation} \label{fo2}
\left( \sum_{s=p}^{P-1} \sum_{j=1}^{m_s} c_{sj} W_{sj} H_{sj} \right)
\left( \sum_{s=p}^{P-1} \sum_{j=1}^{m_s} c_{sj} W_{sj} H_{sj} u_{sj} \right) -
\sum_{j=m_P + 1}^{\bar m_P} c_{Pj} \gamma_j
\end{equation}
subject to
\begin{align} 
\sum_{j=1}^{\bar m_s} v_{sij} = 1, & \; s=p,\dots,P-1, \; i=1,\dots,n_s, \label{constr3} \\
u_{sj} \geq v_{sij}, & \; s=p,\dots,P-1, \; j=1,\dots,\bar m_s, \; i=1,\dots,n_s, \label{constr4} \\
u_{sj} \leq \sum_{i=1}^{n_s} v_{sij}, & \; s=p, \dots, P-1, \; j=1,\dots,\bar m_s, \label{constr4b}
\end{align}
\begin{equation} \label{constr1}
0 \leq t_{sj} \leq \bar H_{sj} \mbox{ and } 0 \leq r_{sj} \leq \bar W_{sj}, \; j=1,\dots,\bar m_s,
\end{equation}
\begin{equation} \label{constr5}
  \frac{1}{2} w_{si} \leq x_{si} \leq \bar W_{sj} - r_{sj} + ( 1 - v_{sij} ) \hat W - \frac{1}{2} w_{si},
  \; s=p,\dots,P-1, i=1,\dots,n_s, j=1,\dots,\bar m_s,
\end{equation}
\begin{equation} \label{constr5b}
  \frac{1}{2} h_{si} \leq y_{si} \leq \bar H_{sj} - t_{sj} + ( 1 - v_{sij} ) \hat H - \frac{1}{2} h_{si},
  \; s=p,\dots,P-1, i=1,\dots,n_s, j=1,\dots,\bar m_s,
\end{equation}
\begin{equation} \label{constr2}
\begin{array}{rcccl}
0 &\leq& \bar H_{s+1,\ell_1} &\leq& \hat H u_{sj},\\
t_{sj} - (1-u_{sj}) \hat H &\leq& \bar H_{s+1,\ell_1} &\leq& t_{sj} + (1-u_{sj}) \hat H,\\[2mm]
0 &\leq& \bar W_{s+1,\ell_1} &\leq& \hat W u_{sj},\\
\bar W_{sj} - r_{sj} - ( 1 - \eta_{sj} ) \hat W - (1-u_{sj}) \hat W &\leq& \bar W_{s+1,\ell_1} &\leq& \bar W_{sj} - r_{sj} + ( 1 - \eta_{sj} ) \hat W + (1-u_{sj}) \hat W,\\
\bar W_{sj} - \eta_{sj} \hat W - (1-u_{sj}) \hat W &\leq& \bar W_{s+1,\ell_1} &\leq& \bar W_{sj} + \eta_{sj} \hat W + (1-u_{sj}) \hat W,\\[2mm]
0 &\leq& \bar W_{s+1,\ell_2} &\leq& \hat W u_{sj},\\
r_{sj} - (1-u_{sj}) \hat W &\leq&\bar W_{s+1,\ell_2} &\leq& r_{sj} + (1-u_{sj}) \hat W,\\[2mm]
0 &\leq& \bar H_{s+1,\ell_2} &\leq& \hat H u_{sj},\\
\bar H_{sj} - ( 1 - \eta_{sj} ) \hat H - (1-u_{sj}) \hat H &\leq& \bar H_{s+1,\ell_2} &\leq& \bar H_{sj} + ( 1 - \eta_{sj} ) \hat H + (1-u_{sj}) \hat H,\\
\bar H_{sj} - t_{sj} - \eta_{sj} \hat H - (1-u_{sj}) \hat H &\leq& \bar H_{s+1,\ell_2} &\leq& \bar H_{sj} - t_{sj} + \eta_{sj} \hat H + (1-u_{sj}) \hat H,\\
\end{array}
\end{equation}
for $s=p,\dots,P-1$ and $j=j_k^s \leq m_s$ for $k=1,\dots\hat m_s$, 
with $\ell_1=m_{s+1} + 2k-1$ and $\ell_2=m_{s+1} + 2k$,
\begin{equation} \label{constr2b}
\begin{array}{rcccl}
\bar H_{sj} - \hat H u_{sj} &\leq& \bar H_{s+1,\ell_1} &\leq& \bar H_{sj} + \hat H u_{sj},\\
t_{sj} - (1-u_{sj}) \hat H &\leq& \bar H_{s+1,\ell_1} &\leq& t_{sj} + (1-u_{sj}) \hat H,\\[2mm]
\bar W_{sj} - \hat W u_{sj} &\leq& \bar W_{s+1,\ell_1} &\leq& \bar W_{sj} + \hat W u_{sj},\\
\bar W_{sj} - r_{sj} - ( 1 - \eta_{sj} ) \hat W - (1-u_{sj}) \hat W &\leq& \bar W_{s+1,\ell_1} &\leq& \bar W_{sj} - r_{sj} + ( 1 - \eta_{sj} ) \hat W + (1-u_{sj}) \hat W,\\
\bar W_{sj} - \eta_{sj} \hat W - (1-u_{sj}) \hat W &\leq& \bar W_{s+1,\ell_1} &\leq& \bar W_{sj} + \eta_{sj} \hat W + (1-u_{sj}) \hat W,\\[2mm]
0 &\leq& \bar W_{s+1,\ell_2} &\leq& \hat W u_{sj},\\
r_{sj} - (1-u_{sj}) \hat W &\leq&\bar W_{s+1,\ell_2} &\leq& r_{sj} + (1-u_{sj}) \hat W,\\[2mm]
0 &\leq& \bar H_{s+1,\ell_2} &\leq& \hat H u_{sj},\\
\bar H_{sj} - ( 1 - \eta_{sj} ) \hat H - (1-u_{sj}) \hat H &\leq& \bar H_{s+1,\ell_2} &\leq& \bar H_{sj} + ( 1 - \eta_{sj} ) \hat H + (1-u_{sj}) \hat H,\\
\bar H_{sj} - t_{sj} - \eta_{sj} \hat H - (1-u_{sj}) \hat H &\leq& \bar H_{s+1,\ell_2} &\leq& \bar H_{sj} - t_{sj} + \eta_{sj} \hat H + (1-u_{sj}) \hat H,\\
\end{array}
\end{equation}
for $s=p,\dots,P-1$ and $j=j_k^s > m_s$ for $k=1,\dots\hat m_s$, 
with $\ell_1=m_{s+1} + 2k-1$ and $\ell_2=m_{s+1} + 2k$,
\begin{equation} \label{constr6}
\begin{array}{rclcl}
  x_{si} - x_{si'} &\geq& \frac{1}{2} (w_{si} + w_{si'}) & - & \hat W \left[ (1-v_{sij}) + (1-v_{si'j}) + \pi_{sii'} + \tau_{sii'} \right],\\[2mm]
- x_{si} + x_{si'} &\geq& \frac{1}{2} (w_{si} + w_{si'}) & - & \hat W \left[ (1-v_{sij}) + (1-v_{si'j}) + \pi_{sii'} + (1-\tau_{sii'}) \right],\\[2mm]
  y_{si} - y_{si'} &\geq& \frac{1}{2} (h_{si} + h_{si'}) & - & \hat H \left[ (1-v_{sij}) + (1-v_{si'j}) + (1-\pi_{sii'}) + \tau_{sii'} \right],\\[2mm]
- y_{si} + y_{si'} &\geq& \frac{1}{2} (h_{si} + h_{si'}) & - & \hat H \left[ (1-v_{sij}) + (1-v_{si'j}) + (1-\pi_{sii'}) + (1-\tau_{sii'}) \right],
\end{array}
\end{equation}
for $s=p,\dots,P-1$, $j=1,\dots,\bar m_s$, $i=1,\dots,n_s$, $i'=i+1,\dots,n_s$,
\begin{equation} \label{ufa}
0 \leq \omega_{j\ell} \leq \bar H_{Pj} \mbox{  and  }
\bar H_{Pj} - (1-\theta_{j\ell}) \hat H \leq \omega_{j\ell} \leq \theta_{j\ell} \hat H
\mbox{ for } j=m_P+1,\dots,\bar m_P, \ell=1,\dots,L,
\end{equation}
\begin{equation} \label{ufa2}
\bar w_i \leq \bar W_{Pj} + \hat W ( 1 - \zeta_{ji} ) \mbox{  and  }
\bar h_i \leq \bar H_{Pj} + \hat H ( 1 - \zeta_{ji} ) \mbox{  for  }
j=m_P+1,\dots,\bar m_P, \; i=1,\dots,d,
\end{equation}
\begin{equation} \label{ufa3}
0 \leq \gamma_j \leq \sum_{\ell=1}^L 2^{\ell-1} \omega_{j\ell} \mbox{  and  }
\gamma_j \leq \left( \sum_{i=1}^d \zeta_{ji} \right) \hat W \hat H
\mbox{  for  } j=m_P+1,\dots,\bar m_P,
\end{equation}
and
\begin{equation} \label{Wbar}
\bar W_{Pj} = \sum_{\ell=1}^{L} 2^{\ell-1} \theta_{j\ell} \mbox{  for  } j=m_P+1,\dots,\bar m_P.
\end{equation}

The objective function~(\ref{fo2}) is given by the cost of the used
\textit{purchasable} objects multiplied by an strict upper bound on
the value of the leftovers at instant~$P$ minus the value of the
leftovers at that instant. Assuming integrality of the constants that
define the instance (see \cite[\S3.7]{bromro}), this composition has
the desired effect of minimizing the cost of the purchased objects
and, among solutions with the same cost, maximizing the value of the
leftovers at instant~$P$. Constraints~(\ref{constr3}) say that each
item must be assigned to exactly one
object. Constraints~(\ref{constr4}) and~(\ref{constr4b}) say that an
object~${\cal O}_{sj}$ is used (i.e.\ $u_{sj}=1$) if and only if at
least an item is allocated to the object. At a first glance, since the
cost of the used objects is being minimized,
constrains~(\ref{constr4b}) may appear to be superfluous. However,
forcing $u_{sj}=0$ when no item is assigned to object~${\cal O}_{sj}$
prevents purchasing and cutting an object to which no item is being
assigned in period~$s$. Constraints~(\ref{constr1}) define the
height~$t_{sj}$ of the top leftover and the width~$r_{sj}$ of the
right-hand-side leftover of object~${\cal
  O}_{sj}$. Constraints~(\ref{constr5},\ref{constr5b}) assume, without
loss of generality, that objects have its bottom-left corner in the
origin of the Cartesian two-dimensional
space. Constraints~(\ref{constr5},\ref{constr5b}) say that if an item~${\cal I}_{si}$
is assigned to an object~${\cal O}_{sj}$, that has dimensions $\bar
W_{sj}$ and $\bar H_{sj}$, then the center~$(x_{si},y_{si})$ of the
item must be placed within the cutting area of the object that goes
from $(0,0)$ to $(\bar W_{sj} - r_{sj}, \bar H_{sj} -
t_{sj})$. Moreover, the constraints say the center of each item must
be far from the borders of the cutting area, so the whole item can be
placed within the object's cutting area. In
constraints~(\ref{constr2}), restrictions on the dimensions of the
leftovers of purchasable objects with positive expiration date are
given; while in~(\ref{constr2b}) the same is done with the dimensions
of leftovers of objects that are leftovers of previous periods. The
difference is that, in the first case, leftovers of a purchasable
object must have null dimensions if the purchasable object is not used
(purchased); while, in the second case, if an object that is a
leftover is not used and its expiration date is strictly positive,
then it must pass to the next instant as its own top or
right-hand-side leftover. Constraints~(\ref{constr6}) model the
non-overlapping of items assigned to the same
object. Constraints~(\ref{ufa},\ref{ufa2},\ref{ufa3},\ref{Wbar}) model
the value~$\gamma_j$ of the $j$-th leftover of the last instant~$P$,
i.e.\ object~${\cal O}_{Pj}$. Recall that, in case a leftover can fit
at least an item from the catalogue, its value is given by its area
(product of its variable dimensions) times the value per unit of area
of the purchasable object that generated the leftover. Otherwise, the
value of the leftover is null. (See \cite[\S3.7.1]{bromro} for
details.) In (\ref{ufa},\ref{ufa2},\ref{ufa3},\ref{Wbar}), the
index~$j$ starts from $m_P+1$. This is the same as saying that it
starts at~$1$, since $m_P=0$ by definition. However, we opted by
writing this way because it simplifies the re-definition of the
meaning of variables $\gamma$ in the next section.  Note also that
variables $\omega$, $\theta$, $\zeta$, and $\gamma$, differently from
all other variables in the model, do not have an index~$s$ that
relates them to an instant of the multi-period scenario. This is
because they all refer to the last instant~$P$. Note that the
\textit{area} of the leftovers of the last instant of the considered
horizon plays a fundamental role in the objective
function~(\ref{fo2}); while for all other instants (including
instant~$P$) only the (variable) dimensions of the leftovers are
required, but not their area.

\section{Forward-looking proposed heuristic}
\label{sec3}

The mixed integer linear programming (MILP)
problem~(\ref{fo2}--\ref{Wbar}) will be named~${\cal M}(p,P)$ from now
on. This notation allow us to refer to the single-period
problem~${\cal M}(\kappa,\kappa+1)$ for some $\kappa \in \{ p, \dots,
P-1 \}$. In problem ${\cal M}(\kappa,\kappa+1)$, it is assumed that
(a) all decisions of instants $s=p,\dots,\kappa-1$ have already been
taken; (b) quantities and dimensions of the ordered items and available
objects (that may be purchasable or leftovers from previous periods)
of instant~$\kappa$ are known; and (c) the last instant of the
considered horizon is pushed back and artificially considered as if it
were $P=\kappa+1$. Thus, the single-period problem~${\cal
  M}(\kappa,\kappa+1)$ coincides with the single-period problem introduced
in~\cite{abmro2}. This means that problem ${\cal M}(\kappa,\kappa+1)$
consists in determining a cutting pattern to produce all items
ordered at instant $\kappa$ minimizing the cost of the purchased
objects and, among solutions with minimum cost, choosing one that
maximizes the value of the leftovers at instant $\kappa+1$.  The
particularity of ${\cal M}(\kappa,\kappa+1)$ with respect to the
single-period problem introduced in~\cite{abmro2} is that in ${\cal
  M}(\kappa,\kappa+1)$ there are some objects that can be used for
free.  This is because the summation in~(\ref{fo2}) goes from~$1$ up
to~$m_\kappa$; meaning that the costs of objects numbered
from~$m_\kappa + 1$ up to~$\bar m_\kappa$, that are the leftovers of
previous periods, are not included in the objective function. Special
attention must also be given to the role of variables~$\gamma_j$ in
${\cal M}(\kappa,\kappa+1)$. On the one hand, in ${\cal M}(p,P)$,
their indices goes from~$1$ (because $m_P=0$ by definition) to $\bar
m_P$ and they represent the areas of the leftovers at instant~$P$. On
the other hand, in ${\cal M}(\kappa,\kappa+1)$, since $P$ is redefined
as if it were $\kappa+1$, the indices of variables $\gamma$ go from
$m_{\kappa+1}+1$ to $\bar m_{\kappa+1}$; and variables $\gamma$
represent the areas of the leftovers at instant $\kappa+1$.

If we assume that the available computational capacity is enough to
solve (with an exact commercial solver) instances with no more than a
single period, a heuristic approach to tackle the original multi-period problem must be considered. At each
instant~$\kappa$, a decision has to be made. The decision consists in
selecting a set of objects (between the $m_\kappa$ purchasable objects
${\cal O}_{\kappa j}$ for $j=1,\dots,m_\kappa$ or leftovers ${\cal
  O}_{\kappa j}$ for $j=m_\kappa+1,\dots,\bar m_\kappa$ from previous
periods) and a cutting pattern to produce, along period
$[\kappa,\kappa+1)$, the~$n_\kappa$ items ordered at
  instant~$\kappa$. The simplest (matheuristic) approach would be to
  solve the single-period problem~${\cal M}(\kappa,\kappa+1)$, for
  $\kappa=p,\dots,P-1$. Substituting~$P$ by $\kappa+1$ in~(\ref{fo2}),
  we have that the objective function of problem ${\cal
    M}(\kappa,\kappa+1)$ is given by
\begin{equation} \label{fo2b}
\left( \sum_{s=p}^{\kappa} \sum_{j=1}^{m_s} c_{sj} W_{sj} H_{sj} \right)
\left( \sum_{s=p}^{\kappa} \sum_{j=1}^{m_s} c_{sj} W_{sj} H_{sj} u_{sj} \right) -
\sum_{j=m_{\kappa+1}+1}^{\bar m_{\kappa+1}} c_{\kappa+1,j} \gamma_j.
\end{equation}
Since in problem~${\cal M}(\kappa,\kappa+1)$ it is assumed that all decisions of 
instants $s=p,\dots,\kappa-1$ have already been taken, we have that
$u_{sj}$ for $s=p,\dots,\kappa-1$ and $j=1,\dots,\bar m_\kappa$ are constant.
Thus, minimizing (\ref{fo2b}) is equivalent to minimizing
\begin{equation} \label{fo2c}
C_\kappa \sum_{j=1}^{m_\kappa} c_{\kappa j} W_{\kappa j} H_{\kappa j} u_{\kappa j} -
\sum_{j=m_{\kappa+1}+1}^{\bar m_{\kappa+1}} c_{\kappa+1,j} \gamma_j,
\end{equation}
where, as in~(\ref{fo2}),
\[
C_\kappa = \sum_{s=p}^{\kappa} \sum_{j=1}^{m_s} c_{sj} W_{sj} H_{sj}
\]
is a constant. Note that $C_\kappa$ corresponds to the total cost of
all purchasable objects existent from the first instant~$p$ up to
instant~$\kappa$. Therefore, it is a strict upper bound on the value
of the leftovers that could have been generated up to instant
$\kappa+1$. Thus, multiplying the first summation in~(\ref{fo2c})
by~$C_\kappa$ has the desired effect of making one unit of this
summation to be more relevant that the whole second summation in~(\ref{fo2c}). It is in this
way that the cost of the used purchasable objects is minimized
and, among solutions with minimum cost, a solution that maximizes the
value of the leftovers at the end of the considered horizon, in this
case instant~$\kappa+1$, is sought. Note that this interpretation requires the first summation in~(\ref{fo2c}) to assume integer values only; see~\cite{abmro2} for details.

The main drawback of a myopic/greedy strategy like the one described
above is that the overall cost is not being minimized at all. This
strategy was used to find the solution depicted in
Figure~\ref{fig3}(a) to the instance described in
Figure~\ref{fig2}. Its flaw is to ignore the effect in the future of
the decisions made at each instant~$\kappa$. Figure~\ref{fig3}(b)
shows that, by buying a more expensive object at instant~$\kappa=0$, a
better solution can be found. In addition, note that, at each
instant~$\kappa$, the number of available objects~$m_\kappa$ is
finite. If we redefine~$m_1=0$ for the instance in Figure~\ref{fig2}
(i.e.\ no purchasable objects available at instant $\kappa=1$), then
the choice of purchasing the small object~${\cal O}_{02}$ at
instant~$\kappa=0$ produces an infeasible solution. This is because
the $3 \times 6$ leftover of~${\cal O}_{02}$ is not enough to produce
the items ordered at~$\kappa=1$ and, since we redefined $m_1=0$, no
other object is available at~$\kappa=1$. So, the myopic approach is
unable to find a feasible solution to the modified instance.

Assume that we are at an instant~$\kappa$ and that at that instant
there are two different objects (one cheaper and smaller and another
more expensive but larger) that can be used to produce the $n_\kappa$
ordered items. Buying the cheapest object would be the myopic
choice. However, assume that buying and using the more expensive
object produces two leftovers that, by being used in forthcoming
periods, produce an overall saving. Quantifying this saving and using
it to decide which object to buy at instant~$\kappa$ is the
looking-ahead strategy we are looking for. An optimistic view would
consist in subtracting from the cost of each object the value of its leftovers. 
We say this view is
optimistic because it assumes that 100\% of the object's leftovers
will be used to produce items (and, thus, savings) in forthcoming
periods. In a more realistic view, each leftover has a different
utilization rate that depends on its dimensions and on the ordered
items in the forthcoming periods.

At any instant~$\kappa+1$, objects ${\cal O}_{\kappa+1,j}$ with index
$j$ between $m_{\kappa+1}+1$ and $m_{\kappa+1}+2 m_\kappa$ correspond
to the $2 m_\kappa$ leftovers of the $m_\kappa$ purchasable objects
that were available at instant~$\kappa$. Therefore, at
instant~$\kappa$, $\gamma_{2j-1}$ and $\gamma_{2j}$ correspond to the
area of the two leftovers of the purchasable object~${\cal O}_{\kappa
  j}$ for $j=1,\dots,m_\kappa$ (nullified when the object is not
purchased or when the leftover does not fit any item from the
catalog). Thus, if object~${\cal O}_{\kappa j}$ is used, then its optimistic 
amortized cost, that assumes that 100\% of its leftovers will be
used, is given by
\begin{equation} \label{costcutarea}
c_{\kappa j} W_{\kappa j} H_{\kappa j} u_{ \kappa j} - 
c_{\kappa j} \gamma_{2j-1} - c_{\kappa j} \gamma_{2j}.
\end{equation}
The value of~(\ref{costcutarea}) is null if object~${\cal O}_{\kappa
  j}$ is not used because in this case $u_{\kappa j} = \gamma_{2j-1} =
\gamma_{2j} = 0$. If utilization rates $\delta_{\kappa,2j-1},
\delta_{\kappa,2j} \in [0,1]$ for $j=1,\dots,m_\kappa$ were known,
then we would be able to compute, at instant~$\kappa$, the more realistic 
amortized cost
\begin{equation} \label{costcutarea2}
c_{\kappa j} W_{\kappa j} H_{\kappa j} u_{ \kappa j} - 
c_{\kappa j} \left( \delta_{\kappa,2j-1} \gamma_{2j-1} + \delta_{\kappa,2j} \gamma_{2j} \right)
\end{equation}
of using object~${\cal O}_{\kappa j}$ to produce the ordered
items. Since we need the summation of costs to assume integer values, 
we would approximate~(\ref{costcutarea2}) by
\begin{equation} \label{costcutarea3}
c_{\kappa j} W_{\kappa j} H_{\kappa j} u_{ \kappa j} - \lfloor 
c_{\kappa j} \left( \delta_{\kappa,2j-1} \gamma_{2j-1} + \delta_{\kappa,2j} \gamma_{2j} \right) \rfloor.
\end{equation}
However, since $\gamma_{2j-1}$ and $\gamma_{2j}$ ($j=1,\dots,m_\kappa$) are variables of the problem, (\ref{costcutarea3}) can not be included in the objective function. (It is not a linear function of continuous and integer variables.) Thus, we need new
integer variables $\lambda_j$ ($j=1,\dots,m_\kappa$) and constraints
\begin{equation} \label{newconstr}
\lambda_j \leq c_{\kappa j} \left(
\delta_{\kappa,2j-1} \gamma_{2j-1} + \delta_{\kappa,2j} \gamma_{2j} \right)
\mbox{ for } j=1,\dots,m_\kappa;
\end{equation}
so we can write the approximation~(\ref{costcutarea3}) of~(\ref{costcutarea2}) as 
\begin{equation} \label{costcutarea4}
c_{\kappa j} W_{\kappa j} H_{\kappa j} u_{ \kappa j} - \lambda_j.
\end{equation}
We call~(\ref{costcutarea4}) the amortized cost of object~${\cal O}_{\kappa j}$.
Thus, including estimations of the leftovers utilization rates, the objective function~(\ref{fo2c}) of
problem~${\cal M}(\kappa,\kappa+1)$ can be substituted by
\begin{equation} \label{fo2d}
C_\kappa \sum_{j=1}^{m_\kappa} 
\left( c_{\kappa j} W_{\kappa j} H_{\kappa j} u_{ \kappa j} - \lambda_j \right)
- \sum_{j=m_\kappa+1}^{\bar m_{\kappa+1}} c_{\kappa+1,j} \gamma_j.
\end{equation}
We call ${\cal M}(\delta; \kappa,\kappa+1)$, the single-period problem
${\cal M}(\kappa,\kappa+1)$ in which the objective function is
replaced with~(\ref{fo2d}) and constraints~(\ref{newconstr}) are included. 
Note that~(\ref{newconstr}) and, in consequence (\ref{fo2d}), depends on the
unknown constants $\delta_{\kappa,2j-1}$ and $\delta_{\kappa,2j}$ for
$j=1,\dots,m_\kappa$.

\begin{figure}[ht!]
\begin{center}
\resizebox{0.9\textwidth}{!}{\definecolor{MyRed}{HTML}{F72C25} 
\definecolor{MyGreen}{RGB}{54,206,65}
\definecolor{MyBlue}{HTML}{0B88F8}

\begin{tikzpicture}[scale=0.35]

\draw[thick,dashed] (-2.5,12.0) -- (62,12.0);
\draw[thick,dashed] (-2.5,-1.0) -- (62,-1.0);
\draw[thick,dashed] (-2.5,-10.5) -- (62,-10.5);

\node[rotate=90] at (-1.9,  5.5) {Available objects};
\node[rotate=90] at (-1.9, -5.8) {Ordered items};

\draw[thick,dashed] (-1,-10.5) -- (-1,14.0);


\node[rotate=0] at (16, 13) {Instant $s=0$};


\node[rotate=0] at (5.0,17/2+1) {${\cal O}_{01}$};
\node[rotate=0] at (5.0,4.25) {21$\times$17};
\draw [thick](0,0) rectangle (20/2, 17/2);
\node[rotate=0] at (15.75,19/2+1) {${\cal O}_{02}$};
\node[rotate=0] at (15.75,4.75) {19$\times$19};
\draw [thick](11,0) rectangle (11+19/2, 19/2);
\node[rotate=0] at (27.5,13/2+1) {${\cal O}_{03}$};
\node[rotate=0] at (27.5,3.25) {24$\times$13};
\draw [thick](21.5,0) rectangle (21.5+24/2, 13/2);


\draw [fill=MyRed, thin](11,-8.00) rectangle (11+19/2,-8+11/2);
\node[rotate=0] at (15.75,-9.25) {${\cal I}_{01}$};
\node[rotate=0] at (15.75,-5.25) {19$\times$11};

\draw[thick,dashed] (34.5,-10.5) -- (34.5,14);


\node[rotate=0] at (41.5, 13) {Instant $s=1$};


\node[rotate=0] at (41.5, 9) {${\cal O}_{11}$};
\node[rotate=0] at (41.5, 4) {10$\times$16};
\draw [thick](39,0) rectangle (39+10/2,16/2);


\draw [fill=MyGreen, thin](35.5,-8.00) rectangle (35.5+7/2,-8+6.0/2);
\node[rotate=0] at (37.25,-6.5) {7$\times$6};
\node[rotate=0] at (37.25,-9.25) {${\cal I}_{11}$};
\draw [fill=MyGreen, thin](40,-8.00) rectangle (40+7/2,-8+5.0/2);
\node[rotate=0] at (41.75,-6.75) {7$\times$5};
\node[rotate=0] at (41.75,-9.25) {${\cal I}_{12}$};
\draw [fill=MyGreen, thin](44.5,-8.00) rectangle (44.5+7/2,-8+4.0/2);
\node[rotate=0] at (46.25,-7.0) {7$\times$4};
\node[rotate=0] at (46.25,-9.25) {${\cal I}_{13}$};

\draw[thick,dashed] (49,-10.5) -- (49,14);


\node[rotate=0] at (55.75, 13) {Instant $s=2$};


\node[rotate=0] at (55.75,7.0) {${\cal O}_{21}$};
\node[rotate=0] at (55.75,3.0) {10$\times$12};

\draw [thick](53.25,0) rectangle (53.25+10/2,12/2);


\draw [fill=MyBlue, thin](52,-8.0) rectangle (52+6/2,-8+5.0/2);
\node[rotate=0] at (53.5,-6.75) {6$\times$5};
\node[rotate=0] at (53.5,-9.25) {${\cal I}_{21}$};
\draw [fill=MyBlue, thin](57,-8.0) rectangle (57+6/2,-8+5.0/2);
\node[rotate=0] at (58.5,-6.75) {6$\times$5};
\node[rotate=0] at (58.5,-9.25) {${\cal I}_{22}$};

\draw[thick,dashed] (62,-10.5) -- (62,14);

\end{tikzpicture}}
\end{center}
\caption{Illustration of a small instance with $p=0$, $P=3$. The
  figure displays the available purchasable objects and the ordered
  items at each instant $s \in \{p,\dots,P-1\}$.}
\label{fig4}
\end{figure}
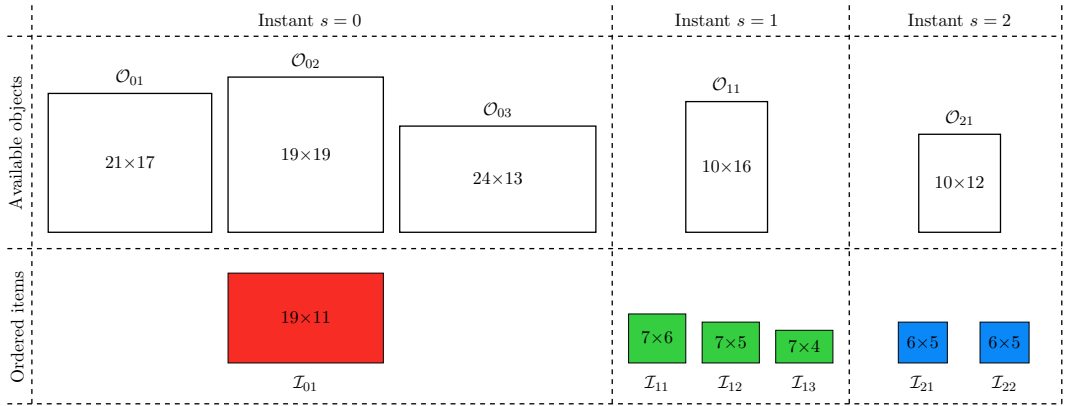

Let us illustrate the idea of amortized costs with an
example. Figure~\ref{fig4} displays the available purchasable objects
and the ordered items of a small instance with $p=0$, $P=3$, and
$\xi=P-p=3$, meaning that usable leftovers generated at any period
remain usable up to instant~$P$. The picture shows the available
purchasable objects and the ordered items at each instant~$s \in
\{0,1,2\}$. The numbers of available purchasable objects and ordered
items at each instant are given by $m_0=3$, $m_1=m_2=1$ and $n_0=1$,
$n_1=3$ and $n_2=2$, respectively. The cost per unit of area of all
the objects is one (i.e.\ $c_{01} = c_{02} = c_{03} = c_{11} = c_{21}
= 1$) and the catalogue with $d=2$ item is composed by two items with
$\bar w_1 = 7$, $\bar h_1 = 4$, $\bar w_2 = 6$, and $\bar h_2 = 5$.

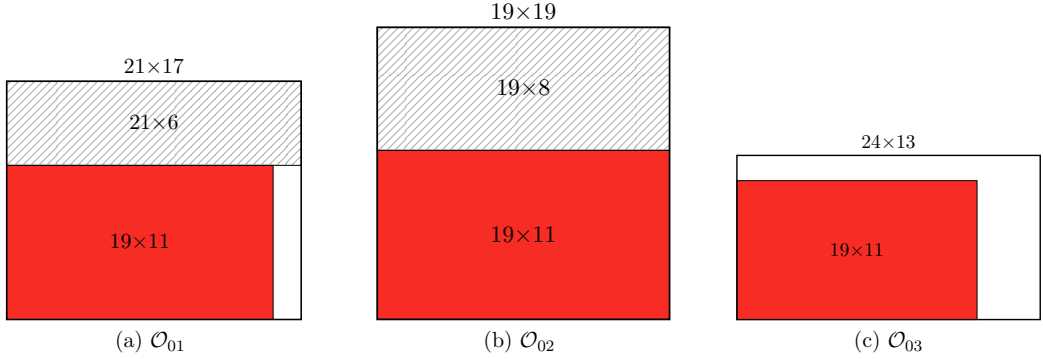
\begin{figure}[ht!]
\begin{center}
\resizebox{0.9\textwidth}{!}{
\begin{tabular}{ccccc}
\resizebox{0.35\textwidth}{!}{\definecolor{MyRed}{HTML}{F72C25} 
\definecolor{Gray}{gray}{0.70}

\begin{tikzpicture}[scale=0.25]
\node[rotate=0] at (10.50,18.0) {21$\times$17};
\draw [fill=MyRed, thin](0.00,0.00) rectangle (19.00,11.00);
\node[rotate=0] at (9.50,5.50) {19$\times$11};

\draw [pattern=north east lines, pattern color=Gray, thin](0.00,11.00) rectangle (21.00,17.00);
\node[rotate=0] at (10.50,14.00) {21$\times$6};

\draw [thick](0,0) rectangle (21,17);
\end{tikzpicture}} & $\phantom{a}$ &
\resizebox{0.35\textwidth}{!}{\definecolor{MyRed}{HTML}{F72C25} 
\definecolor{Gray}{gray}{0.70}

\begin{tikzpicture}[scale=0.25]
\node[rotate=0] at (9.50,20.0) {19$\times$19};
\draw [fill=MyRed, thin](0.00,0.00) rectangle (19.00,11.00);
\node[rotate=0] at (9.50,5.50) {19$\times$11};

\draw [pattern=north east lines, pattern color=Gray, thin](0.00,11.00) rectangle (19.00,19.00);
\node[rotate=0] at (9.50,15.00) {19$\times$8};

\draw [thick](0,0) rectangle (19,19);
\end{tikzpicture}} & $\phantom{a}$ &
\resizebox{0.35\textwidth}{!}{\definecolor{MyRed}{HTML}{F72C25} 

\begin{tikzpicture}[scale=0.25]
\node[rotate=0] at (12.00,14.0) {24$\times$13};
\draw [fill=MyRed, thin](0.00,0.00) rectangle (19.00,11.00);
\node[rotate=0] at (9.50,5.50) {19$\times$11};

\draw [thick](0,0) rectangle (24,13);
\end{tikzpicture}} \\
(a) ${\cal O}_{01}$ && (b) ${\cal O}_{02}$ && (c) ${\cal O}_{03}$\\
\end{tabular}}
\end{center}
\caption{Dashed regions represent the usable leftovers in the
  assignment of item~${\cal I}_{01}$ to the three purchasable objects
  available at instant~$s=0$.}
\label{fig5}
\end{figure}

At instant $s=0$, item ${\cal I}_{01}$ can be assigned to any of the
three available purchasable objects ${\cal O}_{01}$, ${\cal O}_{02}$,
or ${\cal O}_{03}$.  Dashed regions in Figure~\ref{fig5}(a--c)
represent the usable leftovers in each possible assignment. In case
(b) there is only a top usable leftover simply because
$W_{02}=w_{01}$. In case (a) there is also a top usable leftover
only. This is because the right-hand-side leftover has width $W_{02} -
w_{01} < \min\{ \bar w_1, \bar w_2\}$. Thus, it can not fit any item
of the catalogue and, therefore, it is \textit{not} usable. In case
(c), the situation described in case (a) occurs for both, the top and
the right-hand-side leftovers; thus none of them are usable. Since all
the three objects have a unitary cost per unit of area
(i.e.\ $c_{01}=c_{02}=c_{03}=1$), purchasing objects ${\cal O}_{01}$,
${\cal O}_{02}$, and ${\cal O}_{03}$ costs $W_{01} \times H_{01} = 21
\times 17 = 357$, $W_{02} \times H_{02} = 19 \times 19 = 361$, and
$W_{03} \times H_{03} = 24 \times 13 = 312$, respectively. The greedy
choice mandates to buy object ${\cal O}_{03}$, that is the cheapest
one. However, assuming that usable leftovers will be 100\% used to
produce items in forthcoming periods and reducing the value of the
leftovers from the cost of their respective objects, we obtain, for the configurations depicted in Figure~\ref{fig5}, the amortized costs $357 - 21 \times 6 = 231$ and $361 - 19 \times 8 =
209$ for objects ${\cal O}_{01}$ and ${\cal O}_{02}$,
respectively. The amortized cost of object ${\cal O}_{03}$ whose usage
generates no usable leftovers coincides with its actual cost. Thus,
the optimistic forward-looking approach would recommend to purchase
object~${\cal O}_{02}$.

If the myopic approach is applied to the instance of
Figure~\ref{fig4}, then the solution found is to purchase object
${\cal O}_{03}$ at instant $s=0$ and objects ${\cal O}_{11}$ and
${\cal O}_{21}$ at instants $s=1$ and $s=2$, respectively.  This
solution has an overall cost of $592$ and has no usable leftovers at
instant $s=3$. If the optimistic forward-looking approach, that
assumes that 100\% of the usable leftovers will be used in forthcoming
periods, is used, then the solution found is the one illustrated in
Figure~\ref{fig6}(a).  (To simplify the presentation, unused objects
are not being displayed in the figure.)  In this solution, the object
with the smallest amortized cost is chosen at instant~$s=0$,
i.e.\ object~${\cal O}_{02}$. At instant~$s=1$, object~${\cal O}_{11}$
is purchased and ordered items are produced from the purchased object
and from the leftover of the previous period.  At instant~$s=2$ no
object is purchased and the ordered items are produced from a leftover
of the leftover of the object bought at instant~$s=0$. The overall
cost of the solution is 521 and a leftover with value 70 remains
available at instant~$P=3$. (This solution is clearly better than the solution obtained by the myopic approach.) However, it can be noted that the
assumption that 100\% of the leftover of object~${\cal O}_{02}$ would
be used in the next periods turned out to be false. In fact, the
leftover of area~$152$ was used to produce items whose areas
totalize~$102$, i.e.\ an utilization rate of $102/152 \approx
0.67$. If we consider this utilization rate for object ${\cal
  O}_{02}$, then its amortized cost for the configuration depicted in Figure~\ref{fig5}(b) becomes $361 - 102 = 259$.  The
amortized cost of object ${\cal O}_{01}$ (for the configuration in Figure~\ref{fig5}(a)) remains the same, i.e.\ 231,
since there is no new information to update the presumed utilization
rate of 100\% of its usable leftover.  The amortized cost of object
${\cal O}_{03}$ (for the configuration in Figure~\ref{fig5}(a)) continues being 312 as well. Thus, if the problem is
solved once again, object ${\cal O}_{01}$ is chosen at instant~$s=0$
to produce the ordered items of instant~$s=0$. Then, its leftover is
used to produce all ordered items of instant~$s=1$; and object~${\cal
  O}_{21}$ is purchased to produce the items ordered at
instant~$s=2$. This solution, depicted at Figure~\ref{fig6}(b), has an
overall cost of~$477$ and it has no usable leftovers at
instant~$s=3$. In this solution, the actual utilization rate of the
leftover of object ${\cal O}_{02}$ is $314/357 \approx 0.88$; which
increases its amortized cost for the configuration depicted in 
Figure~\ref{fig5}(b) from 231 to $357 - \lfloor (314/357)
\times 126 \rfloor = 247$. Anyway, it continues to be the cheapest
purchasable object at instant~$s=0$. Thus, a new cycle would produce
the same solution.

\begin{figure}[ht!]
\begin{center}
\begin{tabular}{c}
\resizebox{\textwidth}{!}{\definecolor{MyRed}{HTML}{F72C25} 
\definecolor{MyGreen}{RGB}{54,206,65}
\definecolor{MyBlue}{HTML}{0B88F8}
\definecolor{Gray}{gray}{0.70}

\begin{tikzpicture}[scale=0.3]

\draw[thick,dashed] (-1,22.5) -- (66,22.5);
\draw[thick,dashed] (-1,-3.5) -- (66,-3.5);

\draw[thick,dashed] (-1,-3.5) -- (-1,25.0);
\node[rotate=0] at (-1, 27.5) {Instant};
\node[rotate=0] at (-1, 26.0) {$s=0$};


\node[rotate=0] at (9.5, 23.8) {Period [0,1]};
\node[rotate=0] at (9.50,20.0) {19$\times$19};

\draw [fill=MyRed, thin](0.00,0.00) rectangle (19.00,11.00);
\node[rotate=0] at (9.50,5.50) {19$\times$11};

\draw [pattern=north east lines, pattern color=Gray, thin](0.00,11.00) rectangle (19.00,19.00);
\node[rotate=0] at (9.50,15.00) {19$\times$8};

\draw [thick](0,0) rectangle (19,19);

\draw[-, thick] (13.0, 15.0) to (30.5, 15.0);
\draw[->, thick] (30.5,15.0) to (30.5, 10.0);

\draw[thick,dashed] (20,-3.5) -- (20,25.0);			
\node[rotate=0] at (20, 27.5) {Instant};
\node[rotate=0] at (20, 26.0) {$s=1$};


\node[rotate=0] at (36, 23.8) {Period [1,2]};

\draw [fill=MyGreen, thin](21+0.00,0.00) rectangle (21+7.00,6.00);
\node[rotate=0] at (21+3.50,3.0) {7$\times$6};

\draw [fill=MyGreen, thin](21+20.00,0.00) rectangle (21+27.00,5.00);
\node[rotate=0] at (21+23.50,2.50) {7$\times$5};
\draw [fill=MyGreen, thin](21+20.00,5.00) rectangle (21+27.00,9.00);
\node[rotate=0] at (21+23.50,7.00) {7$\times$4};

\draw [pattern=north east lines, pattern color=Gray, thin](21+7.00,0.00) rectangle (21+19.00,8.00);
\node[rotate=0] at (21+13.00,4.00) {12$\times$8};	
\draw [pattern=north east lines, pattern color=Gray, thin](21+20.00,9.00) rectangle (21+30.00,16.00);
\node[rotate=0] at (21+25.00,12.50) {10$\times$7};		

\node[rotate=0] at (30.50,9.0) {19$\times$8};
\draw [thick](21,0) rectangle (40,8);

\node[rotate=0] at (21+25.00,17.0) {10$\times$16};
\draw [thick](21+20,0) rectangle (21+30,16);

\draw[->, thick] (49.0, 12.5) to (53.5, 15.0);

\draw[-, thick] (34.4,2.00) to (34.4, -2.0);
\draw[-, thick] (34.4, -2.0) to (59, -2.0);
\draw[->, thick] (59, -2.0) to (59, -0.5);

\draw[thick,dashed] (52,-3.5) -- (52,25.0);
\node[rotate=0] at (52, 27.5) {Instant};
\node[rotate=0] at (52, 26.0) {$s=2$};


\node[rotate=0] at (59, 23.8) {Period [2,3]};


\draw [fill=MyBlue, thin](53+0.00,0.00) rectangle (53+6.00,5.00);
\node[rotate=0] at (53+3.00,2.50) {6$\times$5};
\draw [fill=MyBlue, thin](53+6.00,0.00) rectangle (53+12.00,5.00);
\node[rotate=0] at (53+9.00,2.50) {6$\times$5};

\node[rotate=0] at (53+6.00,9.0) {12$\times$8};
\draw [thick](53+0,0) rectangle (53+12,8);

\draw [pattern=north east lines, pattern color=Gray, thin](54,12) rectangle (53+11,19);
\draw [thick](54,12) rectangle (53+11,19);
\node[rotate=0] at (54+5.00,15.5) {10$\times$7};

\draw[->, thick] (62,15.5) to (66.8, 15.5);

\draw[thick,dashed] (66,-3.5) -- (66,25.0);
\node[rotate=0] at (66, 27.5) {Instant};
\node[rotate=0] at (66, 26.0) {$s=3$};

\node[rotate=0] at (72, 23.8) {Remaining usable};
\node[rotate=0] at (72, 22.3) {leftovers};

\draw [pattern=north east lines, pattern color=Gray, thin](67+0.00,12.00) rectangle (67+10.00,19.00);
\draw [thick](67+0.00,12.00) rectangle (67+10.00,19.00);
\node[rotate=0] at (67+5.0,15.50) {10$\times$7};

\end{tikzpicture}}\\[2mm]
(a) \\[2mm]
\resizebox{0.9\textwidth}{!}{\definecolor{MyRed}{HTML}{F72C25} 
\definecolor{MyGreen}{RGB}{54,206,65}
\definecolor{MyBlue}{HTML}{0B88F8}
\definecolor{Gray}{gray}{0.70}

\begin{tikzpicture}[scale=0.3]
  
\draw[thick,dashed] (-1.0,22.5) -- (57,22.5);
\draw[thick,dashed] (-1.0,-3.5) -- (57,-3.5);

\draw[thick,dashed] (-1,-3.5) -- (-1,25.0);
\node[rotate=0] at (-1, 27.5) {Instant};
\node[rotate=0] at (-1, 26.0) {$s=0$};

\node[rotate=0] at (10.5, 23.5) {Period [0,1]};

\draw [fill=MyRed, thin](0.00,0.00) rectangle (19.00,11.00);
\node[rotate=0] at (9.50,5.50) {19$\times$11};

\draw [pattern=north east lines, pattern color=Gray, thin](0.00,11.00) rectangle (21.00,17.00);
\node[rotate=0] at (10.50,14.00) {21$\times$6};

\node[rotate=0] at (10.50,17.80) {21$\times$17};
\draw [thick](0,0) rectangle (21,17);

\draw[-, thick] (13.0, 14.0) to (33.5, 14.0);
\draw[->, thick] (33.5,14.0) to (33.5, 8.0);

\draw[thick,dashed] (22,-3.5) -- (22,25.0);			
\node[rotate=0] at (22, 27.5) {Instant};
\node[rotate=0] at (22, 26.0) {$s=1$};


\node[rotate=0] at (33, 23.5) {Period [1,2]};

\draw [fill=MyGreen, thin](23+0.00,0.00) rectangle (23+7.00,6.00);
\node[rotate=0] at (23+3.50,3.0) {7$\times$6};
\draw [fill=MyGreen, thin](23+7.00,0.00) rectangle (23+14.00,5.00);
\node[rotate=0] at (23+10.50,2.50) {7$\times$5};
\draw [fill=MyGreen, thin](23+14.00,0.00) rectangle (23+21.00,4.00);
\node[rotate=0] at (23+17.50,2.00) {7$\times$4};

\node[rotate=0] at (33.50,6.80) {21$\times$6};
\draw [thick](23,0) rectangle (44,6);

\draw[thick,dashed] (45,-3.5) -- (45,25.0);
\node[rotate=0] at (45, 27.5) {Instant};
\node[rotate=0] at (45, 26.0) {$s=2$};

\node[rotate=0] at (51, 23.5) {Period [2,3]};

\draw [fill=MyBlue, thin](46+0.00,0.00) rectangle (46+6.00,5.00);
\node[rotate=0] at (46+3.00,2.50) {6$\times$5};
\draw [fill=MyBlue, thin](46+0.00,5.00) rectangle (46+6.00,10.00);
\node[rotate=0] at (46+3.00,7.50) {6$\times$5};

\node[rotate=0] at (46+5.00,12.80) {10$\times$12};
\draw [thick](46+0,0) rectangle (46+10,12);

\draw[thick,dashed] (57,-3.5) -- (57,25.0);
\node[rotate=0] at (57, 27.5) {Instant};
\node[rotate=0] at (57, 26.0) {$s=3$};

\node[rotate=0] at (63, 23.5) {Remaining usable};
\node[rotate=0] at (63, 22.0) {leftovers};

\end{tikzpicture}}\\[2mm]
(b)
\end{tabular}
\end{center}
\caption{Different feasible solutions to the instance of
  Figure~\ref{fig4}. (a) Solution obtained with the optimistic
  forward-looking approach in which it is assumed that 100\% of each
  usable leftover is used to produce items in forthcoming periods. (b)
  Solution obtained with an adaptive forward-looking approach that
  cycles updating the utilization rate of the leftovers.}
\label{fig6}
\end{figure}

The proposed forward-looking matheuristic approach consists in a
sequence of training cycles.  In each cycle, the $P-p$ single-period
problems ${\cal M}(\delta,\kappa,\kappa+1)$ for $\kappa=p,\dots,P-1$
are solved with fixed values of $\delta_{\kappa,2j-1}$ and
$\delta_{\kappa,2j}$ for $\kappa=p,\dots,P-1$ and
$j=1,\dots,m_\kappa$.  In the $0$th cycle, $\delta^0_{\kappa,2j-1} =
\delta^0_{\kappa,2j} = \delta^{\mathrm{ini}}$ for all $\kappa$ and
$j$, where $\delta^{\mathrm{ini}}\in [0,1]$ is a given constant. At
the end of the $\eta$th cycle, it is possible to compute the actual
fractions $f^\eta_{\kappa,2j-1}$ and $f^\eta_{\kappa,2j}$ of each of
the two leftover ${\cal O}_{\kappa+1,m_{\kappa+1}+2j-1}$ and ${\cal
  O}_{\kappa+1,m_{\kappa+1}+2j}$ of a purchasable object ${\cal
  O}_{\kappa j}$ that were effectively used to produce items in
forthcoming periods for all $\kappa$ and $j$. Note that here we are
talking about items directly produced from the leftovers ${\cal
  O}_{\kappa+1,m_{\kappa+1}+2j-1}$ and ${\cal
  O}_{\kappa+1,m_{\kappa+1}+2j}$ and also about items produced from
leftovers of these leftovers up to $\xi$ periods after purchasing the
purchasable object ${\cal O}_{\kappa j}$. Thus, each
$\delta^\eta_{\kappa,2j-1}$ and $\delta^\eta_{\kappa,2j}$ can be
updated using $f^\eta_{\kappa,2j-1}$ and $f^\eta_{\kappa,2j}$. In
particular, we define
\begin{equation} \label{deltaupdate}
\delta^{\eta+1}_{\kappa,2j-1} = (1-\sigma^\eta)
\delta^\eta_{\kappa,2j-1} + \sigma^\eta f^\eta_{\kappa,2j-1} \; \mbox{
  and } \; \delta^{\eta+1}_{\kappa,2j} = (1-\sigma^\eta)
\delta^\eta_{\kappa,2j} + \sigma^\eta f^\eta_{\kappa,2j},
\end{equation}
where $\sigma \in (0,1)$ is a given constant and $\sigma^\eta$ means
$\sigma$ to the power of~$\eta$. This means that, at the end of each
cycle, new estimations $\delta^{\eta+1}_{\kappa,2j-1}$ and
$\delta^{\eta+1}_{\kappa,2j}$ of the utilization rates of the two
leftovers of object ${\cal O}_{\kappa j}$ for all $\kappa$ and $j$ are
computed as convex combination (parameterized by $\sigma^\eta$) of
their previous values $\delta^{\eta}_{\kappa,2j-1}$ and
$\delta^{\eta}_{\kappa,2j}$ and their actual values
$f^{\eta}_{\kappa,2j-1}$ and $f^{\eta}_{\kappa,2j}$ in the solution
found in the current cycle. Since consecutive cycles with the same values of $\delta$'s
produce the same solution, it makes sense to use
\begin{equation} \label{stopcrit}
\max_{\{\kappa=p,\dots,P-1, j=1,\dots,m_\kappa\}} 
\left\{ \left| \delta^{\eta+1}_{\kappa,2j-1} - \delta^\eta_{\kappa,2j-1} \right|, 
        \left| \delta^{\eta+1}_{\kappa,2j} - \delta^\eta_{\kappa,2j} \right| \right\} \leq \epsilon,
\end{equation}
where $\epsilon>0$ is a given constant, as a stopping criterion.

The forward-looking approach considers the utilization rates of
the top and the right-hand-side leftovers of purchasable objects. We
say these are first-order leftovers. In opposition, when a leftover is
a leftover of a leftover, we say it is a high-order leftover.  When an
item is produced from a first-order leftover, its area plays a role in 
the utilization rate of the first-order leftover
itself. On the other hand, when an item is produced from a high-order
leftover, its area plays a role in the utilization rate of
the first-order leftover that is the ancestor of the used high-order
leftover. Therefore, computing the utilization rate of the
first-order leftovers requires to keep track of their successor
leftovers or, equivalently, to keep track of the ancestors of the
high-order leftovers. Assume we are in the $\eta$th cycle of the
forward-looking approach and that the current instant is
instant~$\kappa$. Before solving the single-period problem ${\cal
  M}(\delta,\kappa,\kappa+1)$ we proceed as follows. (The supra-index
$\eta$ will be omitted for simplicity.) Let $j^\kappa_1 \leq
j^\kappa_2 \leq \dots \leq j^\kappa_{\hat m_\kappa}$ be the indices of
the $\hat m_\kappa$ objects that generate leftovers, that correspond
to the indices~$j$ of objects~${\cal O}_{\kappa j}$ ($j=1,\dots,\bar
m_\kappa$) such that $e_{\kappa j} > 0$. On the one hand, every $j_k
\leq m_\kappa$ is a purchasable object. This means that its two
leftovers are first-order leftovers. So, in this case, we initialize
the used area of the two leftovers as
\[
a_{\kappa+1,m_{\kappa+1}+2k-1}=a_{\kappa+1,m_{\kappa+1}+2k}=0
\]
and the ancestor (or origin) of the two leftovers as themselves, i.e.\
\[
o_{\kappa+1,m_{\kappa+1}+2k-1}= m_{\kappa+1}+2k-1 \; \mbox{ and } \; 
o_{\kappa+1,m_{\kappa+1}+2k}= m_{\kappa+1}+2k.
\]
On the other hand, every $j_k > m_\kappa$ is a leftover that is
generating high-order leftovers. So, in this case, we simply
set the ancestor (or origin) of the two leftovers as
\[
o_{\kappa+1,m_{\kappa+1}+2k-1} = o_{\kappa+1,m_{\kappa+1}+2k} = (\kappa,j_k).
\]
(Note that the ``ancestor'' is a pair that saves the instant and the
index of the first-order leftover that generated the high-order
leftover.) After these initializations, we are ready to solve the
single-period problem ${\cal M}(\delta,\kappa,\kappa+1)$. After
solving it, we can also set the area of the two first-order
leftovers as
\[
A_{\kappa+1,m_{\kappa+1}+2k-1}=\gamma_{m_{\kappa+1}+2k-1} \; \mbox{ and } \;
A_{\kappa+1,m_{\kappa+1}+2k}=\gamma_{m_{\kappa+1}+2k},
\]
for every $j_k \leq m_\kappa$. Then, for each item ${\cal I}_{\kappa
  i}$ ($i=1,\dots,n_\kappa$), we proceed as follows. Variables
$v_{\kappa i j} \in \{0,1\}$ indicate to which object the item was
assigned. By~(\ref{constr1}), only one of the $v_{\kappa i j}$ is equal to one and all
the other are null. Let $j$ be the index such that $v_{\kappa i
  j}=1$. If $j > m_\kappa$, then item ${\cal I}_{\kappa i}$ was
produced from a leftover. So, we add its area, given by $w_{\kappa i}
\times h_{\kappa i}$ to the used area of the ancestor $o_{\kappa j}$
of the leftover ${\cal O}_{\kappa j}$ (that may be itself or not),
i.e.\
\[
a_{o_{\kappa j}} \leftarrow a_{o_{\kappa j}} + w_{\kappa i} \times h_{\kappa i}.
\]
Note that $o_{\kappa j}$ is a pair of the form $o_{\kappa j} =
([o_{\kappa j}]_1,[o_{\kappa j}]_2)$. So, notation $a_{o_{\kappa j}}$
means $a_{[o_{\kappa j}]_1,[o_{\kappa j}]_2}$. At the end of the
current $\eta$th cycle, we are ready to compute the actual 
utilization rates of the first-order leftovers given by
\[
f^{\eta}_{\kappa+1,j} = \frac{a_{\kappa+1,j}}{A_{\kappa+1,j}} 
\mbox{ for } \kappa=p,\dots,P-1 \mbox{ and } j = m_{\kappa+1}, \dots, 2 m_\kappa.
\]
Then, the $\delta$'s are updated as
in~(\ref{deltaupdate}). If~(\ref{stopcrit}) holds, the method
stops. Otherwise, we update $\eta \leftarrow \eta + 1$ and start a new
cycle. The method also stops if in ten consecutive cycles the best solution found so far is not updated.

\section{Numerical experiments}
\label{sec4}

In this section, we aim to evaluate the performance of the proposed forward-looking approach.
The single-period models ${\cal M}(\kappa,\kappa+1)$ and ${\cal
  M}(\delta,\kappa,\kappa+1)$ were implemented in C/C++ using the ILOG
Concert Technology. The myopic and the proposed forward-looking matheuristic approaches
were also implemented in C/C++. Models and code are available at \url{https://github.com/oberlan/bromro2}. Code was compiled with g++ from gcc
version 7.5.0 (GNU compiler collection) with the -O3 option
enable. Numerical experiments were conducted using a machine with
Intel(R) Xeon(R) Silver 4114 CPU @ 2.20GHz with 160GB of RAM memory,
and Ubuntu Server 18.04 operating system. 
Single-period instances within the myopic and the forward-looking approaches were solved using
IBM ILOG CPLEX 12.10.0. A solution is reported as optimal by CPLEX when 
\[
\text{absolute gap} = \text{best feasible solution} - \text{best lower bound} \leq \varepsilon_{\mathrm{abs}}
\]
or
\begin{equation} \label{cplexrelgap}
\text{relative gap} = \frac{|\text{best feasible solution} - \text{best lower bound}|}{10^{-10} + |\text{best feasible solution}|} \leq \varepsilon_{\mathrm{rel}},
\end{equation}
where, by default, $\varepsilon_{\mathrm{abs}} = 10^{-6}$ and $\varepsilon_{\mathrm{rel}} = 10^{-4}$, and ``best feasible solution'' means the smallest value of the objective function related to a feasible solution generated by the method. The objective functions~\eqref{fo2} and~\eqref{fo2d} of models ${\cal M}(\kappa,\kappa+1)$ and ${\cal M}(\delta,\kappa,\kappa+1)$, respectively, for $\kappa=p,\dots,P-1$, assume large integer values at feasible points. Thus, a stopping criterion based on a relative error less than or equal to $\varepsilon_{\mathrm{rel}} = 10^{-4}$ has the undesired effect of stopping the method prematurely. On the other hand, due to the integrality of the objective function values, an absolute error strictly smaller than~$1$ is enough to prove the optimality of the incumbent solution. Therefore, in the numerical experiments, we considered $\varepsilon_{\mathrm{abs}} = 1 - 10^{-6}$ and $\varepsilon_{\mathrm{rel}} = 0$. In addition, \textsc{NodeFileInd} and \textsc{WorkMem} parameters were set to~$3$ and~$32{,}000$, respectively; so the Branch \& Bound tree is partially transferred to disk if memory is exhausted. All other parameters of the solver were used with their default values.

\subsection{Parameters tuning}

In a first set of experiments, we aim to analyze the behavior of the forward-looking approach for variations of its two parameters $\delta_{\mathrm{ini}}$ and $\sigma$. Recall that $\delta_{\mathrm{ini}} \in [0,1]$ corresponds to the initial value of the leftovers utilization fraction; while $\sigma \in (0,1)$ plays a role in the utilization fraction update rule in~(\ref{deltaupdate}). In the numerical experiments of this section,
we considered the twenty five instances with four periods introduced in~\cite{bromro}, varying their leftovers ``expiration date'' parameter $\xi \in \{1,2,3,4\}$. The experiments in~\cite{bromro} show that, when applied to these one hundred instances, CPLEX found an optimal solution in 91 cases. Therefore, we applied the forward-looking approach with all combinations of $\delta_{\mathrm{ini}}$ and $\sigma \in \{ 0.5, 0.55, \dots, 1.0\}$ to these 91 instances and computed the gap to the known optimal solution computed by CPLEX. 

\begin{figure}[ht!]
\begin{center}
\includegraphics[width=\textwidth]{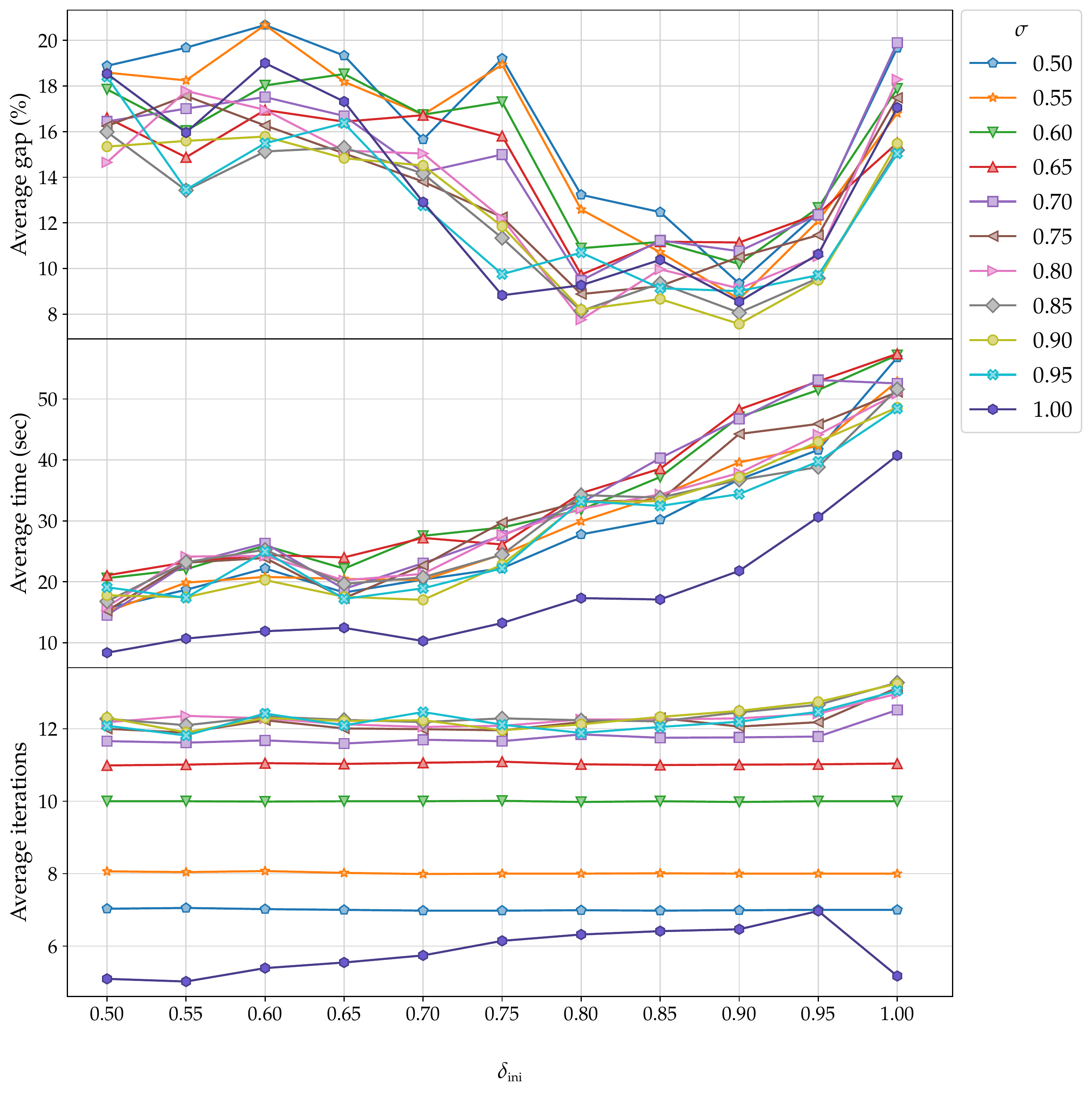}
\caption{Average gap (to optimal solution computed with CPLEX), CPU time (in seconds), and number of cycles of the forward-looking approach for variations of its parameters~$\delta_{\mathrm{ini}}$ and~$\sigma$.}
\label{fig:parameterstuning}
\end{center}
\end{figure}

Figure~\ref{fig:parameterstuning} (top) shows the average gap (over the 91 instances) for each combination of~$\delta_{\mathrm{ini}}$ and~$\sigma$. The figure shows that best results are obtained for the combination $(\delta_{\mathrm{ini}},\sigma)=(0.9,0.9)$. The graphic also shows that, as desired, small variations in the parameters produce a small variation in the average results of the method. It should be noted that the number of cycles (or iterations) $\eta$ that are performed until the satisfaction of the stopping rule~(\ref{stopcrit}) depends on~$\delta_{\mathrm{ini}}$ and~$\sigma$. Figure~\ref{fig:parameterstuning} (middle and bottom) displays the average number of cycles~$\eta$ and the average elapsed CPU time in seconds, as a function of $\delta_{\mathrm{ini}}$ and~$\sigma$. On the one hand, the CPU time has a low dependence on $\sigma$ and, roughly speaking, is an increasing function of $\delta_{\mathrm{ini}}$. On the other hand, the number of cycles has a low dependence on $\delta_{\mathrm{ini}}$ and increases as $\sigma$ increases. Note that, when $\sigma=1$, the rule~(\ref{deltaupdate}) reduces to, at each cycle, discarding information of previous cycles and defining the utilization fraction as the actual utilization fraction of the cycle. In this case, the stopping rule~(\ref{stopcrit}) is satisfied if and only if the utilization rates of all objects are the same for two consecutive cycles. Figure~\ref{fig:parameterstuning} shows that, actually, this phenomenon occur; but it produces a premature stopping with lower quality solutions. However, regardless of the metrics related to computational cost, based on the quality of the solutions obtained, we selected $(\delta_{\mathrm{ini}},\sigma)=(0.9,0.9)$ for the rest of the experiments.

\subsection{Forward-looking versus myopic approach}

In a second set of experiments, we compare the introduced forward-looking approach with $(\delta_{\mathrm{ini}},\sigma)=(0.9,0.9)$ against the myopic approach, that only differs with the forward-looking approach in the objective function that is minimized in each subproblem. In this comparison, a new set of thirty instances with four, eight, and twelve periods is considered. Instances were generated with the random generator introduced in~\cite{bromro}. In order to allow reproducibility, a table describing each instance is given in the Appendix. Table~\ref{tab:BV_CV_CO} shows the number of binary variables, continuous variables, and constraints of each instance when $\xi \in \{1,2,3,4\}$ and, for the instances with eight or twelve periods, $\xi=P$. Note that instances with twelve periods and $\xi=P$ have around $400{,}000$ binary variables, $300{,}000$ continuous variables, and $4{,}000{,}000$ constraints. 

\begin{table}[ht!]
\caption{Number of binary variables (BV), continuous variables (CV), and constraints (CO) of the thirty considered instances.}
\label{tab:BV_CV_CO}
\begin{center}
\resizebox{\textwidth}{!}{
\begin{tabular}{|c|c|rrr|rrr|rrr|rrr|rrr|}
\hline
\multicolumn{2}{|c|}{\multirow{2}{*}{Inst.}} & 
\multicolumn{3}{c|}{$\xi = 1$} & 
\multicolumn{3}{c|}{$\xi = 2$} & 
\multicolumn{3}{c|}{$\xi = 3$} & 
\multicolumn{3}{c|}{$\xi = 4$} & 
\multicolumn{3}{c|}{$\xi = P$} \\
\cline{3-17}
\multicolumn{2}{|c|}{} & \multicolumn{1}{c}{BV} & \multicolumn{1}{c}{CV} & \multicolumn{1}{c|}{CO} & 
\multicolumn{1}{c}{BV} & \multicolumn{1}{c}{CV} & \multicolumn{1}{c|}{CO} & \multicolumn{1}{c}{BV} & 
\multicolumn{1}{c}{CV} & \multicolumn{1}{c|}{CO} & \multicolumn{1}{c}{BV} & \multicolumn{1}{c}{CV} & 
\multicolumn{1}{c|}{CO} & \multicolumn{1}{c}{BV} & \multicolumn{1}{c}{CV} & \multicolumn{1}{c|}{CO} \\ 
\hline
\multirow{10}{*}{\rotatebox{90}{4 periods}}
& 1 & 369 & 150 & 2,664 & 609 & 294 & 5,688 & 897 & 518 & 8,168 & 1,185 & 838 & 9,352 & \multicolumn{3}{c|}{\multirow{10}{*}{\begin{minipage}{4cm}Since instances from~1 to~10 have $P=4$ periods, the case $\xi=P$ coincides with the case $\xi=4$.\end{minipage}}} \\
& 2 & 270 & 150 & 1,683 & 498 & 310 & 3,787 & 786 & 566 & 5,555 & 1,218 & 1,046 & 7,331 &  &  &  \\
& 3 & 298 & 176 & 1,854 & 450 & 304 & 3,122 & 626 & 496 & 4,074 & 754 & 656 & 4,634 &  &  &  \\
& 4 & 397 & 152 & 2,649 & 529 & 240 & 3,805 & 721 & 384 & 5,205 & 1,041 & 704 & 6,453 &  &  &  \\
& 5 & 487 & 150 & 3,752 & 695 & 254 & 6,932 & 951 & 430 & 9,396 & 1,335 & 910 & 11,076 &  &  &  \\
& 6 & 290 & 202 & 1,809 & 546 & 402 & 3,845 & 898 & 754 & 5,757 & 1,042 & 914 & 6,349 &  &  &  \\
& 7 & 572 & 214 & 4,443 & 844 & 358 & 8,667 & 1,164 & 630 & 11,683 & 1,308 & 790 & 12,275 &  &  &  \\
& 8 & 503 & 154 & 3,328 & 675 & 282 & 5,456 & 979 & 426 & 11,560 & 1,235 & 746 & 12,680 &  &  &  \\
& 9 & 318 & 196 & 2,044 & 538 & 380 & 3,672 & 706 & 556 & 4,520 & 1,138 & 1,036 & 6,296 &  &  &  \\
& 10 & 345 & 162 & 2,072 & 525 & 290 & 3,584 & 749 & 434 & 5,784 & 1,069 & 754 & 7,032 &  &  &  \\
\hline
\multirow{10}{*}{\rotatebox{90}{8 periods}}
& 11 & 1,028 & 444 & 9,014 & 1,848 & 868 & 19,982 & 3,368 & 1,668 & 40,422 & 5,672 & 2,820 & 70,806 & 28,904 & 21,764 & 265,142 \\
& 12 & 1,116 & 394 & 9,701 & 1,872 & 754 & 20,881 & 3,040 & 1,378 & 35,801 & 4,848 & 2,338 & 58,953 & 30,096 & 19,874 & 324,841 \\
& 13 & 593 & 362 & 3,824 & 1,105 & 722 & 8,004 & 1,889 & 1,298 & 14,092 & 3,281 & 2,418 & 22,780 & 20,625 & 16,818 & 113,308 \\
& 14 & 921 & 374 & 7,804 & 1,609 & 734 & 17,444 & 2,721 & 1,358 & 32,308 & 4,673 & 2,414 & 60,884 & 23,297 & 18,286 & 238,260 \\
& 15 & 986 & 390 & 8,311 & 1,702 & 742 & 17,911 & 2,982 & 1,430 & 33,255 & 5,334 & 2,710 & 62,487 & 25,910 & 17,558 & 228,343 \\
& 16 & 974 & 408 & 7,886 & 1,782 & 840 & 19,586 & 2,982 & 1,528 & 36,114 & 5,174 & 2,616 & 69,122 & 31,094 & 26,168 & 257,986 \\
& 17 & 1,251 & 394 & 10,836 & 2,071 & 714 & 26,772 & 3,455 & 1,386 & 50,388 & 5,631 & 2,282 & 91,972 & 27,359 & 16,362 & 432,452 \\
& 18 & 839 & 380 & 6,413 & 1,467 & 756 & 13,393 & 2,483 & 1,460 & 23,449 & 3,859 & 2,420 & 36,057 & 18,547 & 15,924 & 130,777 \\
& 19 & 1,020 & 400 & 8,012 & 1,660 & 720 & 16,656 & 2,780 & 1,296 & 31,432 & 4,620 & 2,320 & 53,288 & 22,956 & 17,488 & 202,888 \\
& 20 & 1,141 & 414 & 10,206 & 1,825 & 774 & 19,074 & 2,977 & 1,350 & 34,826 & 5,089 & 2,374 & 66,490 & 30,401 & 19,334 & 377,914 \\
\hline
\multirow{10}{*}{\rotatebox{90}{12 periods}}
& 21 & 1,184 & 514 & 8,941 & 2,056 & 978 & 19,957 & 3,728 & 1,842 & 42,077 & 6,784 & 3,442 & 82,925 & 343,904 & 246,834 & 3,855,917 \\
&  22 & 1,559 & 576 & 13,531 & 2,595 & 1,080 & 29,079 & 4,483 & 1,944 & 58,567 & 7,827 & 3,544 & 108,343 & 307,763 & 248,728 & 2,474,167 \\
&23 & 1,158 & 530 & 8,965 & 2,066 & 1,050 & 19,405 & 3,794 & 1,994 & 40,397 & 6,626 & 3,594 & 73,149 & 326,178 & 295,370 & 2,276,765 \\
& 24 & 1,258 & 562 & 9,857 & 2,198 & 1,058 & 21,645 & 3,838 & 1,986 & 40,837 & 7,086 & 3,714 & 81,909 & 370,446 & 314,050 & 2,672,821 \\
& 25 & 1,443 & 584 & 12,671 & 2,403 & 1,096 & 25,275 & 4,283 & 2,104 & 50,827 & 7,387 & 3,928 & 88,299 & 359,931 & 319,320 & 2,927,211 \\
& 26 & 1,230 & 524 & 9,706 & 2,218 & 1,028 & 22,226 & 3,970 & 1,892 & 44,954 & 7,202 & 3,588 & 83,226 & 395,682 & 263,684 & 3,072,506 \\
& 27 & 1,452 & 558 & 11,777 & 2,480 & 1,054 & 26,525 & 4,472 & 2,030 & 56,773 & 7,928 & 3,790 & 108,821 & 482,392 & 405,326 & 4,270,261 \\
& 28 & 1,587 & 546 & 13,404 & 2,567 & 1,010 & 28,464 & 4,471 & 1,874 & 59,328 & 8,135 & 3,570 & 119,344 & 417,927 & 269,042 & 5,343,952 \\
& 29 & 1,488 & 656 & 12,636 & 2,596 & 1,224 & 27,628 & 4,588 & 2,264 & 54,436 & 8,300 & 4,152 & 106,004 & 480,652 & 339,576 & 6,202,740 \\
& 30 & 1,299 & 630 & 10,782 & 2,363 & 1,198 & 24,670 & 4,259 & 2,238 & 49,086 & 7,315 & 4,126 & 82,830 & 435,731 & 336,414 & 4,289,870 \\ 
\hline
\end{tabular}}
\end{center}
\end{table}

Tables~\ref{tab:result_xi1}--\ref{tab:result_xiP} show the results. The tables show, for the myopic and the forward-looking approaches, the best objective function value found (i.e.\ the value of~\eqref{fo2}), the corresponding cost of the purchased objects, the corresponding value of the leftovers at the final instant of the time horizon, and the CPU time in seconds. In addition, for the forward-looking approach, tables show the gap given by
\begin{equation} \label{flookmyogap}
100 \left( \frac{F_{\mathrm{flook}}-F_{\mathrm{myopic}}}{F_{\mathrm{myopic}}} \right)\%,
\end{equation}
where $F_{\mathrm{flook}}$ is the best objective function value found by the forward-looking approach and $F_{\mathrm{myopic}}$ is the best objective function value found by the myopic approach. It is important to notice that, by definition, the objective function~\eqref{fo2} is dominated by the objects' cost (which is multiplied by an upper bound on the value of the leftovers at the last time instant); while the value of the leftovers at the last time instant plays a ``tie-breaking role''. Thus, a tiny gap may represent a situation where both methods have found a solution with the same cost of the objects but with a relevant difference in the value of the leftovers at instant~$P$. Also note that Tables~\ref{tab:result_xi1}--\ref{tab:result_xiP} do not include averages in the columns corresponding to the leftovers values. This is because, in the considered problem, the main goal is to find a solution that minimizes the overall cost of the objects and, among solutions with minimum costs of the objects, a solution that maximizes the value of the leftovers at instant~$P$. Thus, it makes no sense to compare the value of the leftovers at instant~$P$ of solutions with different objects cost. It would be very easy to construct a solution with high objects cost and plenty of leftovers at the end of the considered time horizon. Given two solutions, the one with lower objects cost is better than the other; and in case the objects cost is identical, the one with the higher value of the leftovers at instant~$P$ is preferable. Solutions must be compared with this objective in mind; so the gaps must be examined carefully.

\begin{table}[ht!]
\caption{Myopic approach versus forward-looking approach considering the scenario with smallest possible use of leftovers, i.e.\ $\xi = 1$.}
\label{tab:result_xi1}
\begin{center}
\resizebox{\textwidth}{!}{
\begin{tabular}{|c|r|rrrr|rrrrr|}
\hline
\multicolumn{2}{|c|}{\multirow{3}{*}{Inst.}} & 
\multicolumn{4}{c|}{Myopic approach} & 
\multicolumn{5}{c|}{Forward-looking approach} \\ 
\cline{3-11}
\multicolumn{2}{|c|}{} & \multicolumn{1}{c}{Best objective} & \multicolumn{1}{c}{Objects} & \multicolumn{1}{c}{Leftovers} & \multicolumn{1}{c|}{CPU} &  \multicolumn{1}{c}{Best objective} & \multicolumn{1}{c}{Objects} & \multicolumn{1}{c}{Leftovers} & \multicolumn{1}{c}{CPU} & \multicolumn{1}{c|}{\multirow{2}{*}{gap (\%)}} \\
\multicolumn{2}{|c|}{} & \multicolumn{1}{c}{function value} & \multicolumn{1}{c}{cost} & \multicolumn{1}{c}{value} & \multicolumn{1}{c|}{time} & \multicolumn{1}{c}{function value} & \multicolumn{1}{c}{cost} & \multicolumn{1}{c}{value} & \multicolumn{1}{c}{time} & \multicolumn{1}{c|}{} \\ 
\hline
\hline
\multirow{10}{*}{\rotatebox{90}{4 periods}}
&1 & \textbf{314,108,050} & 9,155 & 0 & 60.1 & 400,703,843 & 11,679 & 2,647 & 732.3 & 27.5688 \\
&2 & \textbf{187,422,365} & 6,715 & 0 & 30.9 & \textbf{187,422,365} & 6,715 & 0 & 122.5 & 0.0000 \\
&3 & \textbf{340,487,089} & 8,951 & 0 & 3.8 & \textbf{340,487,089} & 8,951 & 0 & 237.9 & 0.0000 \\
&4 & \textbf{309,586,584} & 9,677 & 0 & 1.3 & \textbf{309,586,584} & 9,677 & 0 & 677.9 & 0.0000 \\
&5 & 444,536,794 & 15,954 & 5,462 & 60.3 & \textbf{182,258,424} & 6,541 & 0 & 1,443.5 & -59.0004 \\
&6 & 236,240,392 & 6,246 & 2,066 & 0.2 & \textbf{148,039,222} & 3,914 & 0 & 124.9 & -37.3353 \\
&7 & \textbf{607,520,858} & 13,433 & 0 & 15.8 & \textbf{607,520,858} & 13,433 & 0 & 916.3 & 0.0000 \\
&8 & 241,124,382 & 12,191 & 1,407 & 96.8 & \textbf{191,042,687} & 9,659 & 2,674 & 2,260.8 & -20.7701 \\
&9 & \textbf{226,123,995} & 4,757 & 0 & 0.8 & \textbf{226,123,995} & 4,757 & 0 & 221.0 & 0.0000 \\
&10 & \textbf{354,815,285} & 10,884 & 3,115 & 8.9 & \textbf{354,815,285} & 10,884 & 3,115 & 470.3 & 0.0000 \\
\hline
\multicolumn{2}{|c|}{Avg.} & 326,196,579 & 9,796 &  & 27.9 & 294,800,035 & 8,621 & & 720.7 & -8.9537\\
\hline
\multirow{10}{*}{\rotatebox{90}{8 periods}}
&11 & 1,550,317,180 & 16,165 & 3,310 & 180.3 & \textbf{1,482,704,276} & 15,460 & 2,484 & 4,664.5 & -4.3612 \\
&12 & \textbf{1,625,463,920} & 17,980 & 0 & 103.7 & 1,764,776,484 & 19,521 & 0 & 2,880.4 & 8.5706 \\
&13 & \textbf{1,102,076,378} & 11,453 & 0 & 0.7 & \textbf{1,102,076,378} & 11,453 & 0 & 627.0 & 0.0000 \\
&14 & 1,423,459,632 & 16,701 & 0 & 18.7 & \textbf{1,360,217,488} & 15,959 & 0 & 2,970.7 & -4.4428 \\
&15 & \textbf{1,156,701,480} & 15,396 & 0 & 159.0 & 1,169,398,450 & 15,565 & 0 & 3,086.2 & 1.0977 \\
&16 & \textbf{1,037,649,354} & 12,633 & 0 & 163.6 & 1,299,831,032 & 15,825 & 2,818 & 4,514.8 & 25.2669 \\
&17 & \textbf{1,236,188,630} & 17,285 & 0 & 124.9 & \textbf{1,236,188,630} & 17,285 & 0 & 3,578.7 & 0.0000 \\
&18 & \textbf{1,271,449,952} & 15,649 & 0 & 61.6 & \textbf{1,271,449,952} & 15,649 & 0 & 1,689.7 & 0.0000 \\
&19 & \textbf{1,489,848,521} & 17,883 & 2,092 & 125.9 & 1,589,990,435 & 19,085 & 0 & 3,001.7 & 6.7216 \\
&20 & \textbf{1,464,089,337} & 17,855 & 2,808 & 63.7 & 1,555,845,819 & 18,974 & 3,207 & 2,160.4 & 6.2671 \\
\hline
\multicolumn{2}{|c|}{Avg.} & 1,335,724,438 & 15,683 &  & 104.3 & 1,364,070,347 & 16,200 & & 3,001.5 & 3.6503\\
\hline
\multirow{10}{*}{\rotatebox{90}{12 periods}}
&21 & \textbf{2,905,035,501} & 22,879 & 2,645 & 61.3 & 3,012,458,150 & 23,725 & 0 & 2,638.6 & 3.6978 \\
&22 & \textbf{2,526,326,584} & 22,230 & 1,766 & 181.6 & 2,592,808,909 & 22,815 & 1,766 & 3,926.3 & 2.6316 \\
&23 & \textbf{2,586,793,620} & 22,189 & 0 & 74.0 & 2,910,185,329 & 24,963 & 1,211 & 2,340.8 & 12.5016 \\
&24 & \textbf{2,745,092,742} & 23,139 & 2,523 & 73.7 & 2,753,399,715 & 23,209 & 0 & 2,387.0 & 0.3026 \\
&25 & 3,911,466,834 & 28,039 & 1,705 & 135.8 & \textbf{3,770,293,527} & 27,027 & 0 & 3,244.8 & -3.6092 \\
&26 & 3,966,384,615 & 27,042 & 735 & 124.7 & \textbf{3,927,662,020} & 26,778 & 1,130 & 3,847.7 & -0.9763 \\
&27 & \textbf{3,462,474,633} & 26,709 & 0 & 240.9 & 3,711,377,673 & 28,629 & 0 & 2,842.8 & 7.1886 \\
&28 & 3,106,309,844 & 28,536 & 4,972 & 159.8 & \textbf{2,956,637,816} & 27,161 & 0 & 3,652.1 & -4.8183 \\
&29 & \textbf{2,682,280,094} & 19,795 & 1,791 & 135.6 & 2,802,335,761 & 20,681 & 1,782 & 2,796.8 & 4.4759 \\
&30 & 3,821,604,621 & 24,685 & 3,654 & 182.3 & \textbf{3,437,821,791} & 22,206 & 99 & 1,904.0 & -10.0425 \\
\hline
\multicolumn{2}{|c|}{Avg.} & 3,171,376,909 & 24,524 &  & 137.0 & 3,187,498,069 & 24,719 & & 2,958.1 & 1.1352\\
\hline
\hline
\multicolumn{2}{|c|}{Avg.} & 1,611,099,309 & 16,740 &  & 88.4 & 1,621,848,666 & 16,606 &  & 2,198.7 & -1.3022 \\
\hline
\end{tabular}}
\end{center}
\end{table}

\begin{table}[ht!]
\caption{Myopic approach versus forward-looking approach considering the scenario with low use of leftovers, i.e.\ $\xi = 2$.}
\label{tab:result_xi2}
\begin{center}
\resizebox{\textwidth}{!}{
\begin{tabular}{|c|r|rrrr|rrrrr|}
\hline
\multicolumn{2}{|c|}{\multirow{3}{*}{Inst.}} & 
\multicolumn{4}{c|}{Myopic approach} & 
\multicolumn{5}{c|}{Forward-looking approach} \\ 
\cline{3-11}
\multicolumn{2}{|c|}{} & \multicolumn{1}{c}{Best objective} & \multicolumn{1}{c}{Objects} & \multicolumn{1}{c}{Leftovers} & \multicolumn{1}{c|}{CPU} &  \multicolumn{1}{c}{Best objective} & \multicolumn{1}{c}{Objects} & \multicolumn{1}{c}{Leftovers} & \multicolumn{1}{c}{CPU} & \multicolumn{1}{c|}{\multirow{2}{*}{gap (\%)}} \\
\multicolumn{2}{|c|}{} & \multicolumn{1}{c}{function value} & \multicolumn{1}{c}{cost} & \multicolumn{1}{c}{value} & \multicolumn{1}{c|}{time} & \multicolumn{1}{c}{function value} & \multicolumn{1}{c}{cost} & \multicolumn{1}{c}{value} & \multicolumn{1}{c}{time} & \multicolumn{1}{c|}{} \\ 
\hline
\hline
\multirow{10}{*}{\rotatebox{90}{4 periods}}
&1 & 300,655,883 & 8,763 & 2,647 & 0.9 & \textbf{277,053,250} & 8,075 & 0 & 994.4 & -7.8504 \\
&2 & \textbf{183,066,191} & 6,559 & 2,058 & 0.8 & 187,421,443 & 6,715 & 922 & 441.1 & 2.3791 \\
&3 & 340,482,337 & 8,951 & 4,752 & 63.1 & \textbf{339,152,364} & 8,916 & 3,360 & 535.5 & -0.3906 \\
&4 & 309,582,278 & 9,677 & 4,306 & 76.5 & \textbf{277,209,196} & 8,665 & 1,484 & 781.3 & -10.4570 \\
&5 & 274,293,216 & 9,844 & 0 & 120.1 & \textbf{182,257,329} & 6,541 & 1,095 & 2,407.4 & -33.5538 \\
&6 & 181,132,639 & 4,789 & 1,708 & 2.4 & \textbf{179,167,551} & 4,737 & 0 & 350.0 & -1.0849 \\
&7 & \textbf{527,061,892} & 11,654 & 1,912 & 133.9 & \textbf{527,061,892} & 11,654 & 1,912 & 1,656.1 & 0.0000 \\
&8 & \textbf{166,697,412} & 8,428 & 0 & 36.9 & \textbf{166,697,412} & 8,428 & 0 & 936.7 & 0.0000 \\
&9 & 226,123,365 & 4,757 & 630 & 7.3 & \textbf{226,122,767} & 4,757 & 1,228 & 429.2 & -0.0003 \\
&10 & \textbf{284,400,266} & 8,724 & 2,134 & 61.5 & \textbf{284,400,266} & 8,724 & 2,134 & 618.2 & 0.0000 \\
\hline
\multicolumn{2}{|c|}{Avg.} & 279,349,548 & 8,215 &  & 50.3 & 264,654,347 & 7,721 & & 915.0 & -5.0958\\
\hline
\multirow{10}{*}{\rotatebox{90}{8 periods}}
&11 & 1,425,351,269 & 14,862 & 3,703 & 246.3 & \textbf{1,200,933,896} & 12,522 & 1,036 & 2,632.4 & -15.7447 \\
&12 & 1,492,298,384 & 16,507 & 444 & 301.5 & \textbf{1,492,297,743} & 16,507 & 1,085 & 3,898.8 & 0.0000 \\
&13 & 1,041,838,902 & 10,827 & 0 & 5.8 & \textbf{741,805,859} & 7,709 & 375 & 687.4 & -28.7984 \\
&14 & \textbf{1,151,398,287} & 13,509 & 801 & 53.7 & 1,151,398,632 & 13,509 & 456 & 3,121.8 & 0.0000 \\
&15 & 1,190,883,867 & 15,851 & 1,763 & 137.6 & \textbf{1,104,410,088} & 14,700 & 912 & 3,913.4 & -7.2613 \\
&16 & \textbf{1,037,649,024} & 12,633 & 330 & 161.1 & 1,064,753,959 & 12,963 & 935 & 4,125.3 & 2.6121 \\
&17 & 1,137,778,191 & 15,909 & 1,671 & 190.0 & \textbf{1,083,926,464} & 15,156 & 344 & 6,128.1 & -4.7331 \\
&18 & 1,203,279,118 & 14,810 & 3,762 & 108.9 & \textbf{1,025,673,954} & 12,624 & 798 & 3,779.6 & -14.7601 \\
&19 & \textbf{1,111,449,959} & 13,341 & 2,092 & 193.6 & 1,257,579,545 & 15,095 & 0 & 3,526.3 & 13.1477 \\
&20 & 1,282,624,633 & 15,642 & 3,725 & 126.9 & \textbf{1,242,694,261} & 15,155 & 584 & 3,393.1 & -3.1132 \\
\hline
\multicolumn{2}{|c|}{Avg.} & 1,207,455,163 & 14,389 &  & 152.5 & 1,136,547,440 & 13,594 & & 3,520.6 & -5.8651\\
\hline
\multirow{10}{*}{\rotatebox{90}{12 periods}}
&21 & 2,573,632,748 & 20,269 & 3,258 & 137.4 & \textbf{2,457,199,628} & 19,352 & 1,220 & 5,564.3 & -4.5241 \\
&22 & \textbf{2,286,762,128} & 20,122 & 2,562 & 195.5 & 2,380,519,253 & 20,947 & 2,562 & 3,351.6 & 4.1000 \\
&23 & 2,324,372,040 & 19,938 & 0 & 232.9 & \textbf{2,259,553,182} & 19,382 & 378 & 5,180.9 & -2.7887 \\
&24 & 2,704,400,499 & 22,796 & 2,961 & 173.4 & \textbf{2,517,790,605} & 21,223 & 0 & 2,540.5 & -6.9002 \\
&25 & 3,310,219,229 & 23,729 & 0 & 188.8 & \textbf{2,870,233,075} & 20,575 & 0 & 3,798.6 & -13.2918 \\
&26 & 3,384,818,287 & 23,077 & 688 & 129.9 & \textbf{3,229,489,415} & 22,018 & 735 & 4,707.2 & -4.5890 \\
&27 & \textbf{2,952,610,016} & 22,776 & 2,296 & 308.0 & 3,214,994,976 & 24,800 & 2,624 & 3,507.9 & 8.8865 \\
&28 & 2,991,904,474 & 27,485 & 2,686 & 231.7 & \textbf{2,717,695,698} & 24,966 & 3,198 & 5,075.0 & -9.1650 \\
&29 & 2,369,810,419 & 17,489 & 1,548 & 245.5 & \textbf{2,201,786,731} & 16,249 & 1,516 & 3,951.7 & -7.0902 \\
&30 & 3,189,962,135 & 20,605 & 940 & 174.1 & \textbf{2,837,294,505} & 18,327 & 0 & 3,256.4 & -11.0555 \\
\hline
\multicolumn{2}{|c|}{Avg.} & 2,808,849,198 & 21,829 &  & 201.7 & 2,668,655,707 & 20,784 & & 4,093.4 & -4.6418\\
\hline
\hline
\multicolumn{2}{|c|}{Avg.} & 1,431,884,636 & 14,811 &  & 134.9 & 1,356,619,165 & 14,033 &  & 2,843.0 & -5.2009 \\ 
\hline
\end{tabular}}
\end{center}
\end{table}

\begin{table}[ht!]
\caption{Myopic approach versus forward-looking approach considering the scenario with medium use of leftovers, i.e.\ $\xi = 3$.}
\label{tab:result_xi3}
\begin{center}
\resizebox{\textwidth}{!}{
\begin{tabular}{|c|r|rrrr|rrrrr|}
\hline
\multicolumn{2}{|c|}{\multirow{3}{*}{Inst.}} & 
\multicolumn{4}{c|}{Myopic approach} & 
\multicolumn{5}{c|}{Forward-looking approach} \\ 
\cline{3-11}
\multicolumn{2}{|c|}{} & \multicolumn{1}{c}{Best objective} & \multicolumn{1}{c}{Objects} & \multicolumn{1}{c}{Leftovers} & \multicolumn{1}{c|}{CPU} &  \multicolumn{1}{c}{Best objective} & \multicolumn{1}{c}{Objects} & \multicolumn{1}{c}{Leftovers} & \multicolumn{1}{c}{CPU} & \multicolumn{1}{c|}{\multirow{2}{*}{gap (\%)}} \\
\multicolumn{2}{|c|}{} & \multicolumn{1}{c}{function value} & \multicolumn{1}{c}{cost} & \multicolumn{1}{c}{value} & \multicolumn{1}{c|}{time} & \multicolumn{1}{c}{function value} & \multicolumn{1}{c}{cost} & \multicolumn{1}{c}{value} & \multicolumn{1}{c}{time} & \multicolumn{1}{c|}{} \\ 
\hline
\hline
\multirow{10}{*}{\rotatebox{90}{4 periods}}
&1 & \textbf{177,005,290} & 5,159 & 0 & 60.5 & 277,052,202 & 8,075 & 1,048 & 743.2 & 56.5220 \\
&2 & 183,066,047 & 6,559 & 2,202 & 2.6 & \textbf{165,540,141} & 5,931 & 0 & 534.1 & -9.5735 \\
&3 & 340,482,702 & 8,951 & 4,387 & 63.4 & \textbf{205,638,144} & 5,406 & 690 & 400.9 & -39.6039 \\
&4 & 309,582,332 & 9,677 & 4,252 & 121.9 & \textbf{187,633,080} & 5,865 & 0 & 950.3 & -39.3915 \\
&5 & 274,289,096 & 9,844 & 4,120 & 122.4 & \textbf{182,257,427} & 6,541 & 997 & 2,654.2 & -33.5528 \\
&6 & 181,132,281 & 4,789 & 2,066 & 2.3 & \textbf{92,931,111} & 2,457 & 0 & 215.9 & -48.6943 \\
&7 & \textbf{352,310,540} & 7,790 & 0 & 135.6 & \textbf{352,310,540} & 7,790 & 0 & 1,113.0 & 0.0000 \\
&8 & \textbf{166,694,832} & 8,428 & 2,580 & 96.5 & 166,694,950 & 8,428 & 2,462 & 1,568.4 & 0.0001 \\
&9 & \textbf{226,122,641} & 4,757 & 1,354 & 8.8 & \textbf{226,122,641} & 4,757 & 1,354 & 387.7 & 0.0000 \\
&10 & \textbf{178,974,000} & 5,490 & 0 & 65.3 & \textbf{178,974,000} & 5,490 & 0 & 684.5 & 0.0000 \\
\hline
\multicolumn{2}{|c|}{Avg.} & 238,965,976 & 7,144 &  & 67.9 & 203,515,424 & 6,074 & & 925.2 & -11.4294\\
\hline
\multirow{10}{*}{\rotatebox{90}{8 periods}}
&11 & 1,231,334,604 & 12,839 & 2,530 & 150.6 & \textbf{1,118,166,238} & 11,659 & 1,816 & 2,543.6 & -9.1907 \\
&12 & 1,661,892,542 & 18,383 & 4,190 & 301.7 & \textbf{1,459,391,772} & 16,143 & 0 & 4,490.6 & -12.1849 \\
&13 & 920,593,767 & 9,567 & 375 & 51.8 & \textbf{776,062,690} & 8,065 & 0 & 1,226.3 & -15.6998 \\
&14 & \textbf{1,019,203,389} & 11,958 & 867 & 50.2 & 1,019,203,408 & 11,958 & 848 & 3,878.2 & 0.0000 \\
&15 & 1,190,882,635 & 15,851 & 2,995 & 198.4 & \textbf{1,048,738,758} & 13,959 & 912 & 3,914.1 & -11.9360 \\
&16 & 1,210,381,894 & 14,736 & 3,674 & 143.9 & \textbf{966,517,321} & 11,767 & 525 & 4,190.4 & -20.1477 \\
&17 & 1,292,683,743 & 18,075 & 4,107 & 242.6 & \textbf{1,083,926,384} & 15,156 & 424 & 4,302.7 & -16.1491 \\
&18 & \textbf{911,276,358} & 11,216 & 1,210 & 173.6 & 1,025,673,954 & 12,624 & 798 & 4,277.2 & 12.5536 \\
&19 & \textbf{1,111,449,683} & 13,341 & 2,368 & 206.6 & 1,343,385,248 & 16,125 & 4,627 & 3,499.5 & 20.8678 \\
&20 & 1,218,090,995 & 14,855 & 4,150 & 242.6 & \textbf{1,045,977,464} & 12,756 & 1,780 & 3,820.2 & -14.1298 \\
\hline
\multicolumn{2}{|c|}{Avg.} & 1,176,778,961 & 14,082 &  & 176.2 & 1,088,704,324 & 13,021 & & 3,614.3 & -6.6017\\
\hline
\multirow{10}{*}{\rotatebox{90}{12 periods}}
&21 & 2,263,564,302 & 17,827 & 1,196 & 174.9 & \textbf{2,177,222,273} & 17,147 & 905 & 4,195.6 & -3.8144 \\
&22 & 2,254,372,691 & 19,837 & 3,174 & 225.2 & \textbf{2,151,866,309} & 18,935 & 1,766 & 7,007.4 & -4.5470 \\
&23 & \textbf{2,093,542,769} & 17,958 & 871 & 182.2 & 2,198,114,815 & 18,855 & 1,085 & 4,105.2 & 4.9950 \\
&24 & 2,704,399,467 & 22,796 & 3,993 & 192.5 & \textbf{2,198,543,471} & 18,532 & 349 & 2,637.3 & -18.7049 \\
&25 & 3,374,945,006 & 24,193 & 2,687 & 209.1 & \textbf{2,750,262,215} & 19,715 & 0 & 4,346.5 & -18.5094 \\
&26 & 2,790,050,551 & 19,022 & 1,299 & 218.3 & \textbf{2,658,923,500} & 18,128 & 900 & 3,584.4 & -4.6998 \\
&27 & \textbf{2,719,263,555} & 20,976 & 2,157 & 312.6 & 2,804,696,495 & 21,635 & 0 & 4,507.4 & 3.1418 \\
&28 & 2,947,923,389 & 27,081 & 5,947 & 329.4 & \textbf{2,331,585,459} & 21,419 & 1,205 & 4,223.7 & -20.9075 \\
&29 & 2,280,785,228 & 16,832 & 1,268 & 247.8 & \textbf{2,163,304,424} & 15,965 & 971 & 6,646.0 & -5.1509 \\
&30 & 2,677,059,585 & 17,292 & 1,395 & 244.3 & \textbf{2,546,550,563} & 16,449 & 1,372 & 3,181.9 & -4.8751 \\
\hline
\multicolumn{2}{|c|}{Avg.} & 2,610,590,654 & 20,381 &  & 233.6 & 2,398,106,952 & 18,678 & & 4,443.5 & -7.3072\\
\hline
\hline
\multicolumn{2}{|c|}{Avg.} & 1,342,111,864 & 13,869 &  & 159.3 & 1,230,108,900 & 12,591 &  & 2,994.3 & -8.4461 \\
\hline
\end{tabular}}
\end{center}
\end{table}

\begin{table}[ht!]
\caption{Myopic approach versus forward-looking approach considering the scenario with high use of leftovers, i.e.\ $\xi = 4$.}
\label{tab:result_xi4}
\begin{center}
\resizebox{\textwidth}{!}{
\begin{tabular}{|c|r|rrrr|rrrrr|}
\hline
\multicolumn{2}{|c|}{\multirow{3}{*}{Inst.}} & 
\multicolumn{4}{c|}{Myopic approach} & 
\multicolumn{5}{c|}{Forward-looking approach} \\ 
\cline{3-11}
\multicolumn{2}{|c|}{} & \multicolumn{1}{c}{Best objective} & \multicolumn{1}{c}{Objects} & \multicolumn{1}{c}{Leftovers} & \multicolumn{1}{c|}{CPU} &  \multicolumn{1}{c}{Best objective} & \multicolumn{1}{c}{Objects} & \multicolumn{1}{c}{Leftovers} & \multicolumn{1}{c}{CPU} & \multicolumn{1}{c|}{\multirow{2}{*}{gap (\%)}} \\
\multicolumn{2}{|c|}{} & \multicolumn{1}{c}{function value} & \multicolumn{1}{c}{cost} & \multicolumn{1}{c}{value} & \multicolumn{1}{c|}{time} & \multicolumn{1}{c}{function value} & \multicolumn{1}{c}{cost} & \multicolumn{1}{c}{value} & \multicolumn{1}{c}{time} & \multicolumn{1}{c|}{} \\ 
\hline
\hline
\multirow{10}{*}{\rotatebox{90}{4 periods}}
&1 & \textbf{177,003,277} & 5,159 & 2,013 & 68.2 & 277,048,397 & 8,075 & 4,853 & 1,045.8 & 56.5216 \\
&2 & 183,066,038 & 6,559 & 2,211 & 2.6 & \textbf{165,538,679} & 5,931 & 1,462 & 508.7 & -9.5743 \\
&3 & 340,482,702 & 8,951 & 4,387 & 63.4 & \textbf{205,637,388} & 5,406 & 1,446 & 366.0 & -39.6042 \\
&4 & 309,582,269 & 9,677 & 4,315 & 122.0 & \textbf{309,582,225} & 9,677 & 4,359 & 872.1 & 0.0000 \\
&5 & 274,288,961 & 9,844 & 4,255 & 123.0 & \textbf{182,257,457} & 6,541 & 967 & 1,686.5 & -33.5528 \\
&6 & 181,131,635 & 4,789 & 2,712 & 2.5 & \textbf{92,930,797} & 2,457 & 314 & 353.7 & -48.6943 \\
&7 & \textbf{352,308,306} & 7,790 & 2,234 & 193.9 & 352,308,700 & 7,790 & 1,840 & 1,553.6 & 0.0001 \\
&8 & \textbf{166,694,901} & 8,428 & 2,511 & 96.9 & 166,694,948 & 8,428 & 2,464 & 1,641.1 & 0.0000 \\
&9 & \textbf{226,122,426} & 4,757 & 1,569 & 8.9 & \textbf{226,122,426} & 4,757 & 1,569 & 470.5 & 0.0000 \\
&10 & 178,973,172 & 5,490 & 828 & 65.3 & \textbf{178,972,975} & 5,490 & 1,025 & 669.4 & -0.0001 \\
\hline
\multicolumn{2}{|c|}{Avg.} & 238,965,369 & 7,144 &  & 74.7 & 215,709,399 & 6,455 & & 916.7 & -7.4904\\
\hline
\multirow{10}{*}{\rotatebox{90}{8 periods}}
&11 & \textbf{997,133,908} & 10,397 & 774 & 80.5 & 1,007,107,737 & 10,501 & 1,169 & 3,315.6 & 1.0002 \\
&12 & 1,555,759,725 & 17,209 & 2,711 & 307.0 & \textbf{1,283,103,972} & 14,193 & 0 & 6,345.4 & -17.5256 \\
&13 & 1,006,137,497 & 10,456 & 1,559 & 60.9 & \textbf{741,805,436} & 7,709 & 798 & 1,203.2 & -26.2720 \\
&14 & \textbf{1,019,202,506} & 11,958 & 1,750 & 210.5 & 1,019,203,304 & 11,958 & 952 & 4,411.3 & 0.0001 \\
&15 & 1,190,882,363 & 15,851 & 3,267 & 185.6 & \textbf{1,010,121,118} & 13,445 & 1,732 & 5,956.5 & -15.1788 \\
&16 & 1,210,381,894 & 14,736 & 3,674 & 201.3 & \textbf{1,037,648,275} & 12,633 & 1,079 & 2,826.9 & -14.2710 \\
&17 & 1,137,777,475 & 15,909 & 2,387 & 288.6 & \textbf{1,031,360,898} & 14,421 & 180 & 5,989.4 & -9.3530 \\
&18 & 1,203,278,753 & 14,810 & 4,127 & 188.2 & \textbf{1,025,673,270} & 12,624 & 1,482 & 6,018.4 & -14.7601 \\
&19 & 1,111,449,235 & 13,341 & 2,816 & 208.1 & \textbf{1,026,389,387} & 12,320 & 2,133 & 5,096.1 & -7.6531 \\
&20 & 1,282,623,697 & 15,642 & 4,661 & 308.5 & \textbf{1,049,996,019} & 12,805 & 1,176 & 4,192.9 & -18.1369 \\
\hline
\multicolumn{2}{|c|}{Avg.} & 1,171,462,705 & 14,031 &  & 203.9 & 1,023,240,942 & 12,261 & & 4,535.6 & -12.2150\\
\hline
\multirow{10}{*}{\rotatebox{90}{12 periods}}
&21 & \textbf{2,197,791,998} & 17,309 & 968 & 226.6 & 2,243,630,156 & 17,670 & 424 & 3,523.3 & 2.0856 \\
&22 & 2,254,372,691 & 19,837 & 3,174 & 199.3 & \textbf{1,940,033,795} & 17,071 & 0 & 4,644.8 & -13.9435 \\
&23 & \textbf{2,061,483,504} & 17,683 & 636 & 223.2 & 2,073,956,989 & 17,790 & 1,211 & 5,502.0 & 0.6051 \\
&24 & 2,301,874,270 & 19,403 & 635 & 189.4 & \textbf{2,173,748,840} & 18,323 & 265 & 2,756.2 & -5.5661 \\
&25 & 2,981,413,301 & 21,372 & 2,071 & 141.9 & \textbf{2,779,137,217} & 19,922 & 1,705 & 4,666.2 & -6.7846 \\
&26 & 2,929,977,991 & 19,976 & 1,809 & 253.6 & \textbf{2,658,922,840} & 18,128 & 1,560 & 4,117.1 & -9.2511 \\
&27 & \textbf{2,727,819,075} & 21,042 & 2,679 & 364.1 & 2,727,819,151 & 21,042 & 2,603 & 7,119.4 & 0.0000 \\
&28 & 2,792,803,602 & 25,656 & 5,934 & 337.3 & \textbf{2,421,391,653} & 22,244 & 1,211 & 4,623.5 & -13.2989 \\
&29 & 2,491,626,343 & 18,388 & 2,821 & 178.7 & \textbf{2,147,856,651} & 15,851 & 1,402 & 5,261.6 & -13.7970 \\
&30 & 2,677,058,982 & 17,292 & 1,998 & 245.5 & \textbf{2,262,001,370} & 14,611 & 595 & 4,625.8 & -15.5042 \\
\hline
\multicolumn{2}{|c|}{Avg.} & 2,541,622,176 & 19,796 &  & 236.0 & 2,342,849,866 & 18,265 & & 4,684.0 & -7.5455\\
\hline
\hline
\multicolumn{2}{|c|}{Avg.} & 1,317,350,083 & 13,657 &  & 171.5 & 1,193,933,402 & 12,327 &  & 3,378.8 & -9.0836 \\
\hline
\end{tabular}}
\end{center}
\end{table}

\begin{table}[ht!]
\caption{Myopic approach versus forward-looking approach considering the scenario with unrestricted use of leftovers, i.e.\ $\xi = P$.}
\label{tab:result_xiP}
\begin{center}
\resizebox{\textwidth}{!}{
\begin{tabular}{|c|r|rrrr|rrrrr|}
\hline
\multicolumn{2}{|c|}{\multirow{3}{*}{Inst.}} & 
\multicolumn{4}{c|}{Myopic approach} & 
\multicolumn{5}{c|}{Forward-looking approach} \\ 
\cline{3-11}
\multicolumn{2}{|c|}{} & \multicolumn{1}{c}{Best objective} & \multicolumn{1}{c}{Objects} & \multicolumn{1}{c}{Leftovers} & \multicolumn{1}{c|}{CPU} &  \multicolumn{1}{c}{Best objective} & \multicolumn{1}{c}{Objects} & \multicolumn{1}{c}{Leftovers} & \multicolumn{1}{c}{CPU} & \multicolumn{1}{c|}{\multirow{2}{*}{gap (\%)}} \\
\multicolumn{2}{|c|}{} & \multicolumn{1}{c}{function value} & \multicolumn{1}{c}{cost} & \multicolumn{1}{c}{value} & \multicolumn{1}{c|}{time} & \multicolumn{1}{c}{function value} & \multicolumn{1}{c}{cost} & \multicolumn{1}{c}{value} & \multicolumn{1}{c}{time} & \multicolumn{1}{c|}{} \\ 
\hline
\hline
\multirow{10}{*}{\rotatebox{90}{8 periods}}
&11 & 1,215,891,809 & 12,678 & 4,459 & 189.5 & \textbf{909,955,304} & 9,488 & 824 & 4,170.5 & -25.1615 \\
&12 & 1,555,758,322 & 17,209 & 4,114 & 306.5 & \textbf{1,254,444,657} & 13,876 & 1,247 & 5,958.0 & -19.3676 \\
&13 & 773,366,591 & 8,037 & 1,771 & 68.5 & \textbf{594,579,167} & 6,179 & 1,287 & 1,887.5 & -23.1181 \\
&14 & 1,019,201,343 & 11,958 & 2,913 & 206.2 & \textbf{900,474,723} & 10,565 & 1,357 & 4,569.7 & -11.6490 \\
&15 & 1,190,882,133 & 15,851 & 3,497 & 156.2 & \textbf{1,003,810,128} & 13,361 & 1,802 & 5,501.9 & -15.7087 \\
&16 & 1,210,381,894 & 14,736 & 3,674 & 201.4 & \textbf{980,726,922} & 11,940 & 798 & 4,128.7 & -18.9738 \\
&17 & 1,137,777,262 & 15,909 & 2,600 & 288.4 & \textbf{1,025,352,729} & 14,337 & 837 & 6,173.7 & -9.8811 \\
&18 & 1,203,277,781 & 14,810 & 5,099 & 188.3 & \textbf{925,900,439} & 11,396 & 1,769 & 3,356.6 & -23.0518 \\
&19 & 1,111,448,881 & 13,341 & 3,170 & 268.7 & \textbf{883,095,903} & 10,600 & 697 & 8,360.2 & -20.5455 \\
&20 & 1,190,621,519 & 14,520 & 3,961 & 305.4 & \textbf{873,944,862} & 10,658 & 480 & 6,159.4 & -26.5976 \\
\hline
\multicolumn{2}{|c|}{Avg.} & 1,160,860,754 & 13,905 &  & 217.9 & 935,228,483 & 11,240 & & 5,026.6 & -19.4055\\
\hline
\multirow{10}{*}{\rotatebox{90}{12 periods}}
&21 & 1,983,206,578 & 15,619 & 328 & 173.5 & \textbf{1,873,119,997} & 14,752 & 451 & 6,038.4 & -5.5509 \\
&22 & 1,813,886,558 & 15,961 & 1,287 & 262.2 & \textbf{1,727,516,865} & 15,201 & 780 & 9,290.1 & -4.7616 \\
&23 & 1,741,938,045 & 14,942 & 315 & 250.5 & \textbf{1,691,575,639} & 14,510 & 161 & 8,961.1 & -2.8912 \\
&24 & 2,301,871,943 & 19,403 & 2,962 & 187.4 & \textbf{1,969,220,958} & 16,599 & 1,407 & 3,904.1 & -14.4513 \\
&25 & 2,883,203,059 & 20,668 & 3,609 & 202.4 & \textbf{2,434,989,295} & 17,455 & 660 & 6,673.6 & -15.5457 \\
&26 & 2,790,048,502 & 19,022 & 3,348 & 193.4 & \textbf{2,290,036,133} & 15,613 & 642 & 8,996.0 & -17.9213 \\
&27 & 2,727,820,154 & 21,042 & 1,600 & 247.7 & \textbf{2,391,282,662} & 18,446 & 1,440 & 5,000.3 & -12.3372 \\
&28 & 2,303,933,956 & 21,165 & 3,284 & 309.0 & \textbf{2,039,308,091} & 18,734 & 213 & 11,452.9 & -11.4858 \\
&29 & 1,989,452,967 & 14,682 & 2,079 & 160.0 & \textbf{1,970,076,007} & 14,539 & 2,110 & 5,093.2 & -0.9740 \\
&30 & 2,677,058,826 & 17,292 & 2,154 & 244.8 & \textbf{2,153,321,736} & 13,909 & 99 & 5,106.9 & -19.5639 \\
\hline
\multicolumn{2}{|c|}{Avg.} & 2,321,242,059 & 17,980 &  & 223.1 & 2,054,044,738 & 15,976 & & 7,051.6 & -10.5483\\
\hline
\hline
\multicolumn{2}{|c|}{Avg.} & 1,741,051,406 & 15,942 &  & 220.5 & 1,494,636,611 & 13,608 &  & 6,039.1 & -14.9769 \\
\hline
\end{tabular}}
\end{center}
\end{table}

From what was recalled in the previous paragraph, by the definition of the problem, to win means to find a solution with strictly lower cost of the objects or with equal cost of the objects and strictly higher value of the leftovers at instant~$P$. To tie means to find a solution with the same cost of the objects and the same value of the leftovers at instant~$P$. If the method does not win or does not tie, then it loses. In Tables~\ref{tab:result_xi1}--\ref{tab:result_xiP}, values in bold correspond to the cases in which the method wins or ties. Table~\ref{tab:summary} summarizes the results. Each cell of the table is of the form ``W/T/L G(\%)'', i.e.\ for each combination of number of periods $P \in \{4,8,12\}$ and parameter $\xi \in \{1,2,3,4,P\}$ (comprising 10 instances), it displays the number of instances in which the forward-looking strategy wins, ties, and looses (with respect to the myopic approach), and the average gap given by~\eqref{flookmyogap}. Figures in the table shows that, the larger the chance of taking advantage of leftovers (i.e.\ the larger $\xi$), the larger the number of victories and the larger the gap. Clearly, the way to estimate the future impact of current decisions is heuristic in nature. This fact, associated with an instance in which there is little chance of using leftovers from previous periods (small~$\xi$) occasionally leads the myopic method to obtain better results. This is an expected behavior that does not diminish the value of the proposed method. In the case $\xi=P$, which is the extreme case of the type of instances for which the method was developed, the forward looking approach find better solutions in all instances, with an average gap of, approximately, 15\%.

\begin{table}[ht!]
\caption{Summary of the comparison between the myopic and the forward-looking approaches in the set of thirty instances with 4, 8, and 12 periods and $\xi \in \{1,2,3,4,P\}$.}
\label{tab:summary}
\begin{center}
\resizebox{\textwidth}{!}{
\begin{tabular}{|c|cr|cr|cr|cr|cr|}
\hline
\multirow{2}{*}{Periods} & \multicolumn{2}{c|}{$\xi=1$} & \multicolumn{2}{c|}{$\xi=2$} & \multicolumn{2}{c|}{$\xi=3$} & \multicolumn{2}{c|}{$\xi=4$} & \multicolumn{2}{c|}{$\xi=P$} \\
\cline{2-11}
& W/T/L & G(\%) & W/T/L & G(\%) & W/T/L & G(\%) & W/T/L & G(\%) & W/T/L & G(\%) \\
\hline
\hline
 4 & 3/6/1 &-8.95 & 6/3/1 &-5.01 & 5/3/2 &-11.43 & 6/1/3 & -7.49 & -- & -- \\
 8 & 2/3/5 & 3.65 & 7/0/3 &-5.87 & 7/0/3 & -6.60 & 8/0/2 &-12.22 & 10/0/0 &-19.41 \\
12 & 4/0/6 & 1.14 & 8/0/2 &-4.64 & 8/0/2 & -7.31 & 7/0/3 & -7.55 & 10/0/0 &-10.55 \\
\hline
\hline
Avg. & 9/9/12 & -1.30 & 21/3/6 & -5.20 & 20/3/7 & -8.45 & 21/1/8 & -9.08 & 20/0/0 & -14.98 \\
\hline
\end{tabular}}
\end{center}
\end{table}

\subsection{Assessing the quality of small instances' solutions}

In the previous section, numerical experiments made clear that the forward-looking approach outperforms the myopic approach; and the greater the possibility of economy using leftovers (i.e.\ the larger the parameter $\xi$), the greater the advantage of the method. Since both methods differ in the looking-ahead objective function being minimized at each period, it is clear that this characteristic is well succeeded in that which it is intended to accomplish. On the other hand, we know nothing about how far from the optimal solution are the solutions that the method finds. In this section we perform an experiment comparing the solutions found by the forward-looking approach with the solutions found with CPLEX.

We consider in this experiment the ten instances with four periods and $\xi \in \{1,2,3,4\}$. These problems, i.e.\ the corresponding multi-period models ${\cal M}(p,P)$, were solved with CPLEX, considering a time limit of two hours. The left-hand side of Table~\ref{tab:ADPvsCPLEX} shows the results. The table shows the ceiling of the best lower bound, the best objective function value found, the relative gap~\eqref{cplexrelgap}, and the CPU time in seconds. In addition, Since the value of the objective function~\eqref{fo2} mixes the cost of the objects and the value of the leftovers at instant~$P$ and, thus, it is not very informative by itself, the table shows the cost of the objects and the value of the leftovers associated with each solution found. The right-hand side of the table gathers, from Tables~\ref{tab:result_xi1}--\ref{tab:result_xi4}, the results obtained by the forward-looking approach. In the right-hand side of table, ``gap(\%)'' represents the relative gap between the solutions found by both methods, computed as
\begin{equation} \label{adpcplexgap}
100 \left( \frac{F_{\mathrm{flook}}-F_{\mathrm{cplex}}}{F_{\mathrm{cplex}}} \right)\%,
\end{equation}
where $F_{\mathrm{flook}}$ is the best objective function value found by the forward-looking approach and $F_{\mathrm{cplex}}$ is the best objective function value found by CPLEX. The table shows that, within the imposed CPU time limit, for $\xi=1,2,3,4$, CPLEX closed the gap in 7, 5, 4, and 0 instances (out of 10) respectively; while the average gap~\eqref{adpcplexgap} between CPLEX and the forward-looking approach was 5.8\%, 13.4\%, -1.1\%, and -4.7\%. For the instances with $\xi=1$, the forward-looking approach matched the solution found by CPLEX in 5 cases of which 4 are known to be optimal; and none solution was improved. For the instances with $\xi=2$, the forward-looking approach matched 2 solutions (one of them known to be optimal) and improved other 2 solutions. For the instances with $\xi=3$, the forward-looking approach matched 3 solutions (known to be optimal) and improved other 3. For the instances with $\xi=4$, the forward-looking approach improved 5 solutions found by CPLEX.

First of all, we should note that in this experiment we are considering instances with only four periods, which correspond to the smallest instances being considered in this work. Within this set, the cases in which CPLEX wins are concentrated in the instances with $\xi=1,2$, which correspond to the smallest instances and to the instances in which there is little space to exploit leftovers. It is not expected the proposed method to be advantageous when the instance is so small that it can be solved optimally using CPLEX. On the other hand, the numbers show that (a) the proposed method finds solutions close to the optimal solutions when the optimal solutions are known and that, (b) even considering instances with as few as four periods, the larger the~$\xi$, the greater the advantage of using the proposed method.

To corroborate the statements of the previous paragraph, we also experimented running CPLEX in the 20 most difficult instances, with 8 and 12 periods and $\xi \in \{4,P\}$. Table~\ref{tab:ADPvsCPLEX2} shows the results. In 16 out of the 20 instances with $\xi=4$, CPLEX was able to find a feasible solution; while it failed to find a feasible solution in the other 4 instances. Of those 16 instances, the forward-looking approach found better solutions in 15 instances, with an average gap of -33.62\%. Of the total 20 instances with $\xi=P$, CPLEX found a feasible solution in only 2 instances; and in these two cases the forward-looking approach found better solutions, with an average gap of -74.81\%.

\begin{table}[ht!]
\caption{Comparison of the forward-looking approach solutions with the solutions found by CPLEX (two hours of CPU time limit) in the ten instances with four periods and $\xi \in \{1,2,3,4\}$.}
\label{tab:ADPvsCPLEX}
\begin{center}
\resizebox{\textwidth}{!}{
\begin{tabular}{|c|r|rrrrrr|rrrrr|}
\hline
\multirow{3}{*}{\textbf{$\xi$}} & \multirow{3}{*}{Inst.} & \multicolumn{6}{c|}{CPLEX} & \multicolumn{5}{c|}{Forward-looking approach} \\ \cline{3-13}
 &  & \multicolumn{1}{c}{Ceiling of best} & \multicolumn{1}{c}{Best objective} & \multicolumn{1}{c}{Objects} & \multicolumn{1}{c}{Leftovers} & \multicolumn{1}{c}{\multirow{2}{*}{gap (\%)}} & \multicolumn{1}{c|}{CPU} & \multicolumn{1}{c}{Best objective} & \multicolumn{1}{c}{Objects} & \multicolumn{1}{c}{Leftovers} & \multicolumn{1}{c}{CPU} & \multicolumn{1}{c|}{\multirow{2}{*}{gap (\%)}} \\
 &  & \multicolumn{1}{c}{lower bound} & \multicolumn{1}{c}{function value} & \multicolumn{1}{c}{cost} & \multicolumn{1}{c}{value} & \multicolumn{1}{c}{} & \multicolumn{1}{c|}{time} & \multicolumn{1}{c}{function value} & \multicolumn{1}{c}{cost} & \multicolumn{1}{c}{value} & \multicolumn{1}{c}{time} & \multicolumn{1}{c|}{} \\ 
 \hline 
 \hline
\multirow{11}{*}{1} & 1 & 314,108,050 & \textbf{314,108,050} & 9,155 & 0 & 0.0000 & 0.2 & 400,703,843 & 11,679 & 2,647 & 732.3 & 27.5688 \\
 & 2 & 183,065,474 & \textbf{187,422,365} & 6,715 & 0 & 2.3246 & 7,200.0 & \textbf{187,422,365} & 6,715 & 0 & 122.5 & 0.0000 \\
 & 3 & 339,152,904 & \textbf{339,152,904} & 8,916 & 2,820 & 0.0000 & 0.3 & 340,487,089 & 8,951 & 0 & 237.9 & 0.3934 \\
 & 4 & 309,586,584 & \textbf{309,586,584} & 9,677 & 0 & 0.0000 & 0.2 & \textbf{309,586,584} & 9,677 & 0 & 677.9 & 0.0000 \\
 & 5 & 182,258,424 & \textbf{182,258,424} & 6,541 & 0 & 0.0000 & 50.6 & \textbf{182,258,424} & 6,541 & 0 & 1,443.5 & 0.0000 \\
 & 6 & 148,039,222 & \textbf{148,039,222} & 3,914 & 0 & 0.0000 & 0.1 & \textbf{148,039,222} & 3,914 & 0 & 124.9 & 0.0000 \\
 & 7 & 580,789,740 & \textbf{580,790,380} & 12,842 & 1,912 & 0.0001 & 7,200.0 & 607,520,858 & 13,433 & 0 & 916.3 & 4.6024 \\
 & 8 & 80,065,392 & \textbf{186,634,644} & 9,436 & 0 & 57.1005 & 7,200.0 & 191,042,687 & 9,659 & 2,674 & 2,260.8 & 2.3619 \\
 & 9 & 226,123,995 & \textbf{226,123,995} & 4,757 & 0 & 0.0000 & 0.2 & \textbf{226,123,995} & 4,757 & 0 & 221.0 & 0.0000 \\
 & 10 & 288,510,000 & \textbf{288,510,000} & 8,850 & 0 & 0.0000 & 268.4 & 354,815,285 & 10,884 & 3,115 & 470.3 & 22.9820 \\ \cline{2-13}
 & Avg. &  & 276,262,657 & 8,080 &  & 5.9425 & 2,192.0 & 294,800,035 & 8,621 &  & 720.7 & 5.7909 \\ \hline \hline
\multirow{11}{*}{2} & 1 & 277,053,250 & \textbf{277,053,250} & 8,075 & 0 & 0.0000 & 0.3 & \textbf{277,053,250} & 8,075 & 0 & 994.4 & 0.0000 \\
 & 2 & 125,208,746 & \textbf{125,208,746} & 4,486 & 0 & 0.0000 & 1,942.0 & 187,421,443 & 6,715 & 922 & 441.1 & 49.6872 \\
 & 3 & 205,638,834 & \textbf{205,638,834} & 5,406 & 0 & 0.0000 & 0.3 & 339,152,364 & 8,916 & 3,360 & 535.5 & 64.9262 \\
 & 4 & 216,808,300 & \textbf{277,209,196} & 8,665 & 1,484 & 21.7889 & 7,200.0 & \textbf{277,209,196} & 8,665 & 1,484 & 781.3 & 0.0000 \\
 & 5 & 162,301,312 & 235,866,007 & 8,465 & 2,753 & 31.1892 & 7,200.0 & \textbf{182,257,329} & 6,541 & 1,095 & 2,407.4 & -22.7284 \\
 & 6 & 136,049,331 & \textbf{136,049,331} & 3,597 & 0 & 0.0000 & 1.8 & 179,167,551 & 4,737 & 0 & 350.0 & 31.6931 \\
 & 7 & 406,039,028 & \textbf{491,516,168} & 10,868 & 0 & 17.3905 & 7,200.0 & 527,061,892 & 11,654 & 1,912 & 1,656.1 & 7.2319 \\
 & 8 & 80,062,469 & 186,631,619 & 9,436 & 3,025 & 57.1013 & 7,200.0 & \textbf{166,697,412} & 8,428 & 0 & 936.7 & -10.6810 \\
 & 9 & 226,117,517 & \textbf{226,122,466} & 4,757 & 1,529 & 0.0022 & 7,200.0 & 226,122,767 & 4,757 & 1,228 & 429.2 & 0.0001 \\
 & 10 & 249,388,985 & \textbf{249,388,985} & 7,650 & 1,015 & 0.0000 & 551.8 & 284,400,266 & 8,724 & 2,134 & 618.2 & 14.0388 \\ \cline{2-13}
 & Avg. &  & 241,068,460 & 7,141 &  & 12.7472 & 3,849.6 & 264,654,347 & 7,721 &  & 915.0 & 13.4168 \\ \hline \hline
\multirow{11}{*}{3} & 1 & 177,005,290 & \textbf{177,005,290} & 5,159 & 0 & 0.0000 & 4.5 & 277,052,202 & 8,075 & 1,048 & 743.2 & 56.5226 \\
 & 2 & 111,055,089 & \textbf{165,538,722} & 5,931 & 1,419 & 32.9129 & 7,200.0 & 165,540,141 & 5,931 & 0 & 534.1 & 0.0000 \\
 & 3 & 115,486,404 & \textbf{205,637,382} & 5,406 & 1,452 & 43.8398 & 7,200.0 & 205,638,144 & 5,406 & 690 & 400.9 & 0.0000 \\
 & 4 & 127,232,184 & 309,582,248 & 9,677 & 4,336 & 58.9020 & 7,200.0 & \textbf{187,633,080} & 5,865 & 0 & 950.3 & -39.3924 \\
 & 5 & 73,560,960 & 203,212,152 & 7,293 & 0 & 63.8009 & 7,200.0 & \textbf{182,257,427} & 6,541 & 997 & 2,654.2 & -10.3113 \\
 & 6 & 92,931,111 & \textbf{92,931,111} & 2,457 & 0 & 0.0000 & 44.0 & \textbf{92,931,111} & 2,457 & 0 & 215.9 & 0.0000 \\
 & 7 & 352,310,540 & \textbf{352,310,540} & 7,790 & 0 & 0.0000 & 6.3 & \textbf{352,310,540} & 7,790 & 0 & 1,113.0 & 0.0000 \\
 & 8 & 36,551,592 & 203,244,701 & 10,276 & 4,303 & 82.0160 & 7,200.0 & \textbf{166,694,950} & 8,428 & 2,462 & 1,568.4 & -17.9837 \\
 & 9 & 226,118,625 & \textbf{226,122,492} & 4,757 & 1,503 & 0.0017 & 7,200.0 & 226,122,641 & 4,757 & 1,354 & 387.7 & 0.0000 \\
 & 10 & 178,974,000 & \textbf{178,974,000} & 5,490 & 0 & 0.0000 & 9.9 & \textbf{178,974,000} & 5,490 & 0 & 684.5 & 0.0000 \\ \cline{2-13}
 & Avg. &  & 211,455,864 & 6,424 &  & 28.1473 & 4,326.5 & 203,515,424 & 6,074 &  & 925.2 & -1.1165 \\ \hline \hline
\multirow{11}{*}{4} & 1 & 176,987,996 & \textbf{177,003,339} & 5,159 & 1,951 & 0.0087 & 7,200.0 & 277,048,397 & 8,075 & 4,853 & 1,045.8 & 56.5216 \\
 & 2 & 111,048,262 & 169,836,152 & 6,085 & 2,283 & 34.6145 & 7,200.0 & \textbf{165,538,679} & 5,931 & 1,462 & 508.7 & -2.5304 \\
 & 3 & 115,477,259 & \textbf{205,637,085} & 5,406 & 1,749 & 43.8441 & 7,200.0 & 205,637,388 & 5,406 & 1,446 & 366.0 & 0.0001 \\
 & 4 & 127,219,757 & 314,860,300 & 9,842 & 4,964 & 59.5949 & 7,200.0 & \textbf{309,582,225} & 9,677 & 4,359 & 872.1 & -1.6763 \\
 & 5 & 53,604,707 & 276,768,471 & 9,933 & 4,641 & 80.6319 & 7,200.0 & \textbf{182,257,457} & 6,541 & 967 & 1,686.5 & -34.1480 \\
 & 6 & 92,925,615 & \textbf{92,930,733} & 2,457 & 378 & 0.0055 & 7,200.0 & 92,930,797 & 2,457 & 314 & 353.7 & 0.0001 \\
 & 7 & 352,266,598 & 406,035,779 & 8,978 & 3,249 & 13.2425 & 7,200.0 & \textbf{352,308,700} & 7,790 & 1,840 & 1,553.6 & -13.2321 \\
 & 8 & 36,542,003 & 347,683,703 & 17,579 & 11,338 & 89.4899 & 7,200.0 & \textbf{166,694,948} & 8,428 & 2,464 & 1,641.1 & -52.0556 \\
 & 9 & 226,115,028 & \textbf{226,122,389} & 4,757 & 1,606 & 0.0033 & 7,200.0 & 226,122,426 & 4,757 & 1,569 & 470.5 & 0.0000 \\
 & 10 & 178,945,728 & \textbf{178,972,785} & 5,490 & 1,215 & 0.0151 & 7,200.0 & 178,972,975 & 5,490 & 1,025 & 669.4 & 0.0001 \\ \cline{2-13}
 & Avg.  &  & 239,585,074 & 7,569 &  & 32.1450 & 7,200.0 & 215,709,399 & 6,455 &  & 916.7 & -4.7121 \\ \hline
\end{tabular}}
\end{center}
\end{table}

\begin{table}[ht!]
\caption{Comparison of the forward-looking approach solutions with the solutions found by CPLEX (two hours of CPU time limit) in the twenty instances with eight and twelve periods and $\xi \in \{4,P\}$.}
\label{tab:ADPvsCPLEX2}
\begin{center}
\resizebox{\textwidth}{!}{
\begin{tabular}{|c|r|rrrrrr|rrrrr|}
\hline
\multirow{3}{*}{\textbf{$\xi$}} & \multirow{3}{*}{Inst.} & \multicolumn{6}{c|}{CPLEX} & \multicolumn{5}{c|}{Forward-looking approach} \\ \cline{3-13}
 &  & \multicolumn{1}{c}{Ceiling of best} & \multicolumn{1}{c}{Best objective} & \multicolumn{1}{c}{Objects} & \multicolumn{1}{c}{Leftovers} & \multicolumn{1}{c}{\multirow{2}{*}{gap (\%)}} & \multicolumn{1}{c|}{CPU} & \multicolumn{1}{c}{Best objective} & \multicolumn{1}{c}{Objects} & \multicolumn{1}{c}{Leftovers} & \multicolumn{1}{c}{CPU} & \multicolumn{1}{c|}{\multirow{2}{*}{gap (\%)}} \\
 &  & \multicolumn{1}{c}{lower bound} & \multicolumn{1}{c}{function value} & \multicolumn{1}{c}{cost} & \multicolumn{1}{c}{value} & \multicolumn{1}{c}{} & \multicolumn{1}{c|}{time} & \multicolumn{1}{c}{function value} & \multicolumn{1}{c}{cost} & \multicolumn{1}{c}{value} & \multicolumn{1}{c}{time} & \multicolumn{1}{c|}{} \\ 
 \hline 
 \hline
\multirow{20}{*}{4} &  11 & 473,194,584 & 1,693,795,866 & 17,661 & 0 & 72.0631 & 7,200.0 & \textbf{1,007,107,737} & 10,501 & 1,169 & 3,315.6 & -40.5414 \\
 & 12 & 436,284,255 & 2,225,925,343 & 24,622 & 1,945 & 80.3999 & 7,200.0 & \textbf{1,283,103,972} & 14,193 & 0 & 6,345.4 & -42.3564 \\
 & 13 & 372,869,732 & \textbf{612,863,014} & 6,369 & 380 & 39.1594 & 7,200.0 & 741,805,436 & 7,709 & 798 & 1,203.2 & 21.0394 \\
 & 14 & 222,188,546 & 1,658,613,465 & 19,460 & 1,255 & 86.6040 & 7,200.0 & \textbf{1,019,203,304} & 11,958 & 952 & 4,411.3 & -38.5509 \\
 & 15 & 262,867,798 & 2,360,960,250 & 31,425 & 0 & 88.8661 & 7,200.0 & \textbf{1,010,121,118} & 13,445 & 1,732 & 5,956.5 & -57.2157 \\
 & 16 & 383,736,775 & 2,283,682,476 & 27,803 & 338 & 83.1966 & 7,200.0 & \textbf{1,037,648,275} & 12,633 & 1,079 & 2,826.9 & -54.5625 \\
 & 17 & \multicolumn{5}{l}{Solution not found} & 7,200.0 & \textbf{1,031,360,898} & 14,421 & 180 & 5,989.4 & -- \\
 & 18 & 441,237,340 & 1,205,070,336 & 14,832 & 0 & 63.3849 & 7,200.0 & \textbf{1,025,673,270} & 12,624 & 1,482 & 6,018.4 & -14.8869 \\
 & 19 & 618,229,016 & 1,814,594,666 & 21,781 & 2,225 & 65.9302 & 7,200.0 & \textbf{1,026,389,387} & 12,320 & 2,133 & 5,096.1 & -43.4370 \\
 & 20 & 538,241,436 & 1,151,593,956 & 14,044 & 0 & 53.2612 & 7,200.0 & \textbf{1,049,996,019} & 12,805 & 1,176 & 4,192.9 & -8.8224 \\ \cline{2-13}
 & Avg.  &  & 1,667,455,486 & 19,777 &  & 70.3184 & 7,200.0 & 1,023,240,942 & 12,261 &  & 4,535.6 & -31.0371 \\ \cline{2-13}
 & 21 & 997,752,674 & 2,977,411,727 & 23,449 & 1,599 & 66.4893 & 7,200.0 & \textbf{2,243,630,156} & 17,670 & 424 & 3,523.3 & -24.6449 \\
 & 22 & 848,244,208 & 3,940,979,158 & 34,678 & 2,152 & 78.4763 & 7,200.0 & \textbf{1,940,033,795} & 17,071 & 0 & 4,644.8 & -50.7728 \\
 & 23 & 1,113,793,829 & 2,572,336,398 & 22,065 & 1,302 & 56.7011 & 7,200.0 & \textbf{2,073,956,989} & 17,790 & 1,211 & 5,502.0 & -19.3746 \\
 & 24 & 992,737,680 & 3,628,332,840 & 30,584 & 0 & 72.6393 & 7,200.0 & \textbf{2,173,748,840} & 18,323 & 265 & 2,756.2 & -40.0896 \\
 & 25 & 664,024,760 & 5,652,018,841 & 40,516 & 3,675 & 88.2515 & 7,200.0 & \textbf{2,779,137,217} & 19,922 & 1,705 & 4,666.2 & -50.8293 \\
 & 26 & 930,198,690 & 4,379,714,810 & 29,860 & 690 & 78.7612 & 7,200.0 & \textbf{2,658,922,840} & 18,128 & 1,560 & 4,117.1 & -39.2900 \\
 & 27 & \multicolumn{5}{l}{Solution not found} & 7,200.0 & \textbf{2,727,819,151} & 21,042 & 2,603 & 7,119.4 & -- \\
 & 28 & \multicolumn{5}{l}{Solution not found} & 7,200.0 & \textbf{2,421,391,653} & 22,244 & 1,211 & 4,623.5 & -- \\
 & 29 & \multicolumn{5}{l}{Solution not found} & 7,200.0 & \textbf{2,147,856,651} & 15,851 & 1,402 & 5,261.6 & -- \\
 & 30 & \multicolumn{5}{l}{Solution not found} & 7,200.0 & \textbf{2,262,001,370} & 14,611 & 595 & 4,625.8 & -- \\  \cline{2-13}
 & Avg.  &  & 3,858,465,629 &  30,192 &  & 73.5531 & 7,200.0 & 2,342,849,866 & 18,265 &  & 4,684.0 & -37.5002 \\ 
 \hline
 \hline
\multirow{20}{*}{$P$} & 11 & \multicolumn{5}{l}{Solution not found} & 7,200.0 & \textbf{909,955,304} & 9,488 & 824 & 4,170.5 & -- \\
 & 12 & \multicolumn{5}{l}{Solution not found} & 7,200.0 & \textbf{1,254,444,657} & 13,876 & 1,247 & 5,958.0 & -- \\
 & 13 & 191,248,610 & 2,094,633,046 & 21,768 & 14,522 & 90.8700 & 7,200.0 & \textbf{594,579,167} & 6,179 & 1,287 & 1,887.5 & -71.6100 \\
 & 14 & \multicolumn{5}{l}{Solution not found} & 7,200.0 & \textbf{900,474,723} & 10,565 & 1,357 & 4,569.7 & -- \\
 & 15 & \multicolumn{5}{l}{Solution not found} & 7,200.0 & \textbf{1,003,810,128} & 13,361 & 1,802 & 5,501.9 & -- \\
 & 16 & \multicolumn{5}{l}{Solution not found} & 7,200.0 & \textbf{980,726,922} & 11,940 & 798 & 4,128.7 & -- \\
 & 17 & \multicolumn{5}{l}{Solution not found} & 7,200.0 & \textbf{1,025,352,729} & 14,337 & 837 & 6,173.7 & -- \\
 & 18 & 303,165,605 & 4,212,836,468 & 51,852 & 34,828 & 92.8000 & 7,200.0 & \textbf{925,900,439} & 11,396 & 1,769 & 3,356.6 & -78.0200 \\
 & 19 & \multicolumn{5}{l}{Solution not found} & 7,200.0 & \textbf{883,095,903} & 10,600 & 697 & 8,360.2 & -- \\
 & 20 & \multicolumn{5}{l}{Solution not found} & 7,200.0 & \textbf{873,944,862} & 10,658 & 480 & 6,159.4 & -- \\
\cline{2-13}
 & Avg.  &  & 3,153,734,757 & 36,810 &  & 91.8350 & 7,200.0 & 935,228,483 & 11,240 & & 5,026.6 & -74.8150 \\ 
 \cline{2-13}
 & 21 & \multicolumn{5}{l}{Solution not found} & 7,200.0 & \textbf{1,873,119,997} & 14,752 & 451 & 6,038.4 & -- \\
 & 22 & \multicolumn{5}{l}{Solution not found} & 7,200.0 & \textbf{1,727,516,865} & 15,201 & 780 & 9,290.1 & -- \\
 & 23 & \multicolumn{5}{l}{Solution not found} & 7,200.0 & \textbf{1,691,575,639} & 14,510 & 161 & 8,961.1 & -- \\
 & 24 & \multicolumn{5}{l}{Solution not found} & 7,200.0 & \textbf{1,969,220,958} & 16,599 & 1,407 & 3,904.1 & -- \\
 & 25 & \multicolumn{5}{l}{Solution not found} & 7,200.0 & \textbf{2,434,989,295} & 17,455 & 660 & 6,673.6 & -- \\
 & 26 & \multicolumn{5}{l}{Solution not found} & 7,200.0 & \textbf{2,290,036,133} & 15,613 & 642 & 8,996.0 & -- \\
 & 27 & \multicolumn{5}{l}{Solution not found} & 7,200.0 & \textbf{2,391,282,662} & 18,446 & 1,440 & 5,000.3 & -- \\
 & 28 & \multicolumn{5}{l}{Solution not found} & 7,200.0 & \textbf{2,039,308,091} & 18,734 & 213 & 11,452.9 & -- \\
 & 29 & \multicolumn{5}{l}{Solution not found} & 7,200.0 & \textbf{1,970,076,007} & 14,539 & 2,110 & 5,093.2 & -- \\
 & 30 & \multicolumn{5}{l}{Solution not found} & 7,200.0 & \textbf{2,153,321,736} & 13,909 & 99 & 5,106.9 & -- \\ 
\cline{2-13}
 & Avg.  &  & -- & -- &  & -- & 7,200.0 & 2,054,044,738 & 15,976 & & 7,051.7 & -- \\ 
 \hline
 \end{tabular}}
\end{center}
\end{table}

\section{Concluding remarks}

This work contributes to the literature on two-dimensional cutting stock problems with usable leftovers, which is very limited. A forward-looking approach for the multi-period two-dimensional non-guillotine cutting stock problem with usable leftovers, proposed in~\cite{bromro}, was introduced, this being the first method reported in the literature to address this problem. The method solves a sequence of single-period subproblems and differs with a myopic approach in the objective function being minimized. On the one hand, the myopic approach greedily minimizes the cost of the raw material that must be purchased to produce the orders of the period. On the other and, the forward-looking approach takes into consideration the future impact of the decisions of the period. This looking-head feature allows the method to suggest the purchase of some extra raw material whose leftovers are expected to be used in future periods, resulting in a lower overall cost. Numerical experiments shown the efficiency and effectiveness of the method. In summary, the proposed approach greatly improves the solution found with a commercial solver or with a myopic approach in problems with a reasonable number of periods in which usable leftovers can be used over several periods after they have been generated, i.e. a scenario in which leftovers can play a relevant role.

On the one hand, the proposed method can be applied to instances with a large number of periods. On the other hand, solving the single-period subproblems exactly limits the applicability to instances with larger single-period subproblems. Then, devising a heuristic method for the single-period problem would have an immediate impact on methods for solving the multi-period problem. That will be a subject of future work. In another line of research, the problem introduced in~\cite{bromro} and for which a method was developed in the present work, could be modified to take into account situations that sometimes arise in practice. For example, the problem could be modified to allow the anticipated production of items included in future period orders. In this case, storage costs and production capacity limits for each period could be considered.



\section*{Appendix}

Table~\ref{tab:dadosInstancia} describes in detail the thirty instances with four, eight, and twelve periods considered in the present work. Instances were generated with the random instances generator introduced in~\cite{bromro}, where additional twenty five instances with four periods are also described. The number of binary variables, continuous variables, and constraints of each instance, for $\xi \in \{0,1,2,3,4\}$ is given in Table~\ref{tab:BV_CV_CO}. The random instances generator is available at \url{https://github.com/oberlan/bromro2}. 

\setcounter{table}{9}
\begin{table}[ht!]
\caption{Description of the considered thirty instances with four, eight, and twelve periods.}
\label{tab:dadosInstancia}
\begin{center}
\resizebox{\textwidth}{!}{
\begin{tabular}{|c|c|rl|rrll|}
\hline 
\multirow{2}{*}{Inst.} & \multirow{2}{*}{$P$} & \multicolumn{ 2}{c|}{Objects} & \multicolumn{ 4}{c|}{Items} \\ 
\cline{ 3- 8}
& & \multicolumn{1}{l}{$m^s$} & $W_j^s \times H_j^s$ & \multicolumn{1}{l}{$n^s$} & 
\multicolumn{1}{l}{$\tilde{n}^s$} & $d$ & $w_i^s \times h_i^s$ \\ 
\hline	 
\hline	 
\multirow{4}{*}{1} & \multirow{4}{*}{4} & 2 & 77 $\times$ 100, 67 $\times$ 77 & 4 & 2 & \multirow{4}{*}{2} & 2(\underline{6 $\times$ 5}), 2(9 $\times$ 6)\\
& & 2 & 81 $\times$ 36, 95 $\times$ 33 & 6 & 3 &  & 8 $\times$ 11, 2(15 $\times$ 6), 3(18 $\times$ 14)\\
& & 2 & 54 $\times$ 74, 78 $\times$ 100 & 10 & 4 &  & 3(6 $\times$ 8), 3(7 $\times$ 9), 2(17 $\times$ 13), 2(13 $\times$ 8)\\
& & 1 & 53 $\times$ 68 & 7 & 4 &  & 3(10 $\times$ 5), \underline{5 $\times$ 6}, 18 $\times$ 15, 2(16 $\times$ 14)\\ 
\hline
\multirow{4}{*}{2} & \multirow{4}{*}{4} & 3 & 49 $\times$ 82, 34 $\times$ 70, 57 $\times$ 76 & 6 & 3 & \multirow{4}{*}{2} & 2(\underline{7 $\times$ 5}), 19 $\times$ 15, 3(17 $\times$ 15)\\
& & 2 & 39 $\times$ 54, 39 $\times$ 41 & 4 & 3 &  & 17 $\times$ 20, 2(9 $\times$ 20), 20 $\times$ 17\\
& & 2 & 38 $\times$ 72, 85 $\times$ 96 & 7 & 4 &  & 10 $\times$ 10, 3(14 $\times$ 8), 18 $\times$ 20, 2(\underline{6 $\times$ 18})\\
& & 1 & 43 $\times$ 60 & 4 & 2 &  & 14 $\times$ 8, 3(18 $\times$ 7)\\ 
\hline
\multirow{4}{*}{3} & \multirow{4}{*}{4} & 1 & 69 $\times$ 44 & 4 & 3 & \multirow{4}{*}{1} & 15 $\times$ 6, 14 $\times$ 8, 2(8 $\times$ 11)\\
& & 2 & 30 $\times$ 79, 39 $\times$ 92 & 6 & 2 &  & 3(8 $\times$ 17), 3(18 $\times$ 17)\\
& & 2 & 83 $\times$ 89, 65 $\times$ 91 & 8 & 4 &  & 13 $\times$ 11, 3(\underline{8 $\times$ 5}), 2(9 $\times$ 14), 2(18 $\times$ 17)\\
& & 3 & 96 $\times$ 73, 54 $\times$ 65, 95 $\times$ 55 & 4 & 3 &  & 14 $\times$ 14, 2(10 $\times$ 15), 12 $\times$ 13\\ 
\hline
\multirow{4}{*}{4} & \multirow{4}{*}{4} & 2 & 41 $\times$ 97, 85 $\times$ 69 & 4 & 3 & \multirow{4}{*}{3} & 14 $\times$ 12, 2(18 $\times$ 8), 19 $\times$ 15\\
& & 1 & 90 $\times$ 95 & 13 & 5 &  & 3(14 $\times$ 10), 3(\underline{8 $\times$ 10}), 2(19 $\times$ 12), 3(\underline{17 $\times$ 6}), 2(17 $\times$ 9)\\
& & 1 & 75 $\times$ 76 & 6 & 4 &  & 18 $\times$ 12, \underline{5 $\times$ 20}, 2(15 $\times$ 20), 2(9 $\times$ 11)\\
& & 2 & 80 $\times$ 35, 85 $\times$ 60 & 5 & 3 &  & 19 $\times$ 13, 3(16 $\times$ 14), 12 $\times$ 18\\ \hline
\multirow{4}{*}{5} & \multirow{4}{*}{4} & 3 & 91 $\times$ 59, 52 $\times$ 37, 40 $\times$ 66 & 4 & 2 & \multirow{4}{*}{1} & 2(\underline{6 $\times$ 5}), 2(19 $\times$ 14)\\
& & 1 & 88 $\times$ 90 & 13 & 5 &  & 2(20 $\times$ 9), 3(7 $\times$ 7), 2(7 $\times$ 15), 3(19 $\times$ 8), 3(11 $\times$ 16)\\
& & 1 & 83 $\times$ 47 & 10 & 4 &  & 3(20 $\times$ 8), 2(20 $\times$ 9), 3(14 $\times$ 18), 2(17 $\times$ 17)\\
& & 1 & 65 $\times$ 94 & 6 & 2 &  & 3(7 $\times$ 8), 3(17 $\times$ 9)\\ 
\hline
\multirow{4}{*}{6} & \multirow{4}{*}{4} & 1 & 63 $\times$ 39 & 3 & 2 & \multirow{4}{*}{2} & 2(\underline{5 $\times$ 8}), 12 $\times$ 7\\
& & 4 & 81 $\times$ 87, 2(38 $\times$ 30), 81 $\times$ 54 & 5 & 2 &  & 2(14 $\times$ 18), 3(7 $\times$ 19)\\
& & 3 & 83 $\times$ 91, 47 $\times$ 31, 52 $\times$ 71 & 3 & 3 &  & 16 $\times$ 6, 16 $\times$ 9, 7 $\times$ 11\\
& & 3 & 53 $\times$ 56, 44 $\times$ 53, 37 $\times$ 99 & 6 & 4 &  & 3(\underline{11 $\times$ 5}), 14 $\times$ 19, 2(6 $\times$ 12)\\ 
\hline
\multirow{4}{*}{7} & \multirow{4}{*}{4} & 1 & 82 $\times$ 95 & 7 & 5 & \multirow{4}{*}{2} & 12 $\times$ 17, 10 $\times$ 5, 9 $\times$ 17, 3(6 $\times$ 18), 12 $\times$ 20\\
& & 3 & 57 $\times$ 54, 2(33 $\times$ 36) & 8 & 4 &  & 3(20 $\times$ 17), 2(11 $\times$ 8), 2(15 $\times$ 14), 18 $\times$ 5\\
& & 2 & 95 $\times$ 67, 99 $\times$ 57 & 9 & 4 &  & 2(10 $\times$ 17), \underline{5 $\times$ 8}, 3(6 $\times$ 6), 3(14 $\times$ 9)\\
& & 3 & 42 $\times$ 92, 88 $\times$ 100, 85 $\times$ 86 & 11 & 5 &  & 15 $\times$ 15, 2(16 $\times$ 10), 2(\underline{6 $\times$ 5}), 3(16 $\times$ 12), 3(12 $\times$ 17)\\ 
\hline
\multirow{4}{*}{8} & \multirow{4}{*}{4} & 2 & 2(56 $\times$ 33) & 10 & 5 & \multirow{4}{*}{1} & 3(13 $\times$ 17), 2(17 $\times$ 7), 17 $\times$ 10, 7 $\times$ 13, 3(15 $\times$ 10)\\
& & 1 & 70 $\times$ 94 & 8 & 5 &  & 12 $\times$ 8, 2(9 $\times$ 7), 18 $\times$ 5, 3(14 $\times$ 13), 6 $\times$ 9\\
& & 2 & 55 $\times$ 40, 60 $\times$ 59 & 4 & 2 &  & 3(16 $\times$ 9), 11 $\times$ 14\\
& & 1 & 71 $\times$ 53 & 13 & 5 &  & 3(16 $\times$ 19), 2(\underline{5 $\times$ 5}), 2(18 $\times$ 6), 3(11 $\times$ 14), 3(12 $\times$ 18)\\ 
\hline
\multirow{4}{*}{9} & \multirow{4}{*}{4} & 3 & 66 $\times$ 99, 93 $\times$ 54, 30 $\times$ 74 & 4 & 2 & \multirow{4}{*}{2} & 3(\underline{5 $\times$ 16}), 11 $\times$ 16\\
& & 1 & 56 $\times$ 93 & 8 & 4 &  & 3(14 $\times$ 12), 14 $\times$ 10, 3(\underline{10 $\times$ 7}), 19 $\times$ 10\\
& & 3 & 67 $\times$ 68, 43 $\times$ 59, 93 $\times$ 74 & 6 & 3 &  & 2(18 $\times$ 10), 13 $\times$ 17, 3(19 $\times$ 7)\\
& & 3 & 93 $\times$ 92, 86 $\times$ 53, 43 $\times$ 34 & 2 & 2 &  & 14 $\times$ 20, 12 $\times$ 9\\ 
\hline
\multirow{4}{*}{10} & \multirow{4}{*}{4} & 2 & 78 $\times$ 95, 61 $\times$ 90 & 7 & 3 & \multirow{4}{*}{3} & 2(9 $\times$ 19), 2(12 $\times$ 6), 3(\underline{6 $\times$ 12})\\
& & 1 & 62 $\times$ 79 & 7 & 4 &  & 3(20 $\times$ 15), 3(15 $\times$ 7), 16 $\times$ 18\\
& & 2 & 36 $\times$ 60, 35 $\times$ 96 & 6 & 3 &  & 2(16 $\times$ 16), 7 $\times$ 17, 3(\underline{9 $\times$ 8})\\
& & 2 & 84 $\times$ 72, 33 $\times$ 98 & 7 & 4 &  & 2(\underline{11 $\times$ 5}), 3(7 $\times$ 17), 20 $\times$ 16, 19 $\times$ 12\\ 
\hline
\multirow{8}{*}{11} & \multirow{8}{*}{8} & 3 & 61 $\times$ 85, 37 $\times$ 95, 84 $\times$ 46 & 4 & 2 & \multirow{8}{*}{2} & 16 $\times$ 20, 3(\underline{5 $\times$ 6})\\
& & 3 & 72 $\times$ 55, 62 $\times$ 41, 35 $\times$ 33 & 6 & 3 &  & 3(\underline{8 $\times$ 5}), 8 $\times$ 17, 2(14 $\times$ 5)\\
& & 3 & 90 $\times$ 68, 47 $\times$ 44, 52 $\times$ 63 & 3 & 2 &  & 2(14 $\times$ 16), 14 $\times$ 17\\
& & 4 & 2(39 $\times$ 56), 81 $\times$ 81, 61 $\times$ 44 & 10 & 4 &  & 2(19 $\times$ 19), 3(7 $\times$ 15), 2(16 $\times$ 15), 3(18 $\times$ 9)\\
& & 2 & 54 $\times$ 97, 40 $\times$ 86 & 7 & 3 &  & 3(17 $\times$ 7), 13 $\times$ 6, 3(10 $\times$ 6)\\
& & 4 & 2(33 $\times$ 43), 93 $\times$ 77, 84 $\times$ 70 & 9 & 3 &  & 3(16 $\times$ 16), 3(10 $\times$ 11), 3(14 $\times$ 11)\\
& & 3 & 41 $\times$ 74, 86 $\times$ 91, 62 $\times$ 30 & 8 & 3 &  & 3(19 $\times$ 8), 3(8 $\times$ 9), 2(7 $\times$ 6)\\
& & 3 & 100 $\times$ 37, 69 $\times$ 65, 83 $\times$ 62 & 7 & 5 &  & 2(13 $\times$ 18), 7 $\times$ 8, 13 $\times$ 12, 2(12 $\times$ 7), 14 $\times$ 18\\ 
\hline
\multicolumn{8}{|r|}{{Continued on next page}}\\
\hline
\end{tabular}}
\end{center}
\end{table}

\setcounter{table}{9}
\begin{table}[ht!]
\caption{-- continued from previous page}
\label{tab7part2}
\begin{center}
\resizebox{\textwidth}{!}{
\begin{tabular}{|c|c|rl|rrll|}
\hline 
\multirow{2}{*}{Inst.} & \multirow{2}{*}{$P$} & \multicolumn{ 2}{c|}{Objects} & \multicolumn{ 4}{c|}{Items} \\ 
\cline{ 3- 8}
& & \multicolumn{1}{l}{$m^s$} & $W_j^s \times H_j^s$ & \multicolumn{1}{l}{$n^s$} & 
\multicolumn{1}{l}{$\tilde{n}^s$} & $d$ & $w_i^s \times h_i^s$ \\ 
\hline	 
\hline
\multirow{8}{*}{12} & \multirow{8}{*}{8} & 3 & 68 $\times$ 37, 70 $\times$ 43, 97 $\times$ 52 & 7 & 5 & \multirow{8}{*}{3} & 20 $\times$ 14, 14 $\times$ 10, 20 $\times$ 15, 3(17 $\times$ 19), 7 $\times$ 13\\
& & 3 & 88 $\times$ 39, 89 $\times$ 35, 55 $\times$ 79 & 8 & 4 &  & 3(7 $\times$ 17), 3(15 $\times$ 11), 10 $\times$ 12, 20 $\times$ 10\\
& & 2 & 66 $\times$ 77, 58 $\times$ 88 & 11 & 5 &  & 18 $\times$ 9, 3(10 $\times$ 20), 2(18 $\times$ 5), 2(7 $\times$ 12), 3(14 $\times$ 15)\\
& & 2 & 95 $\times$ 69, 85 $\times$ 97 & 8 & 4 &  & 2(20 $\times$ 14), 14 $\times$ 18, 3(8 $\times$ 17), 2(14 $\times$ 15)\\
& & 2 & 30 $\times$ 84, 65 $\times$ 56 & 6 & 3 &  & 3(5 $\times$ 20), 2(12 $\times$ 13), 14 $\times$ 9\\
& & 3 & 75 $\times$ 63, 42 $\times$ 55, 73 $\times$ 89 & 5 & 3 &  & \underline{5 $\times$ 9}, 2(17 $\times$ 15), 2(11 $\times$ 9)\\
& & 3 & 90 $\times$ 57, 67 $\times$ 52, 76 $\times$ 86 & 10 & 4 &  & 3(20 $\times$ 15), 13 $\times$ 19, 3(\underline{16 $\times$ 5}), 3(19 $\times$ 5)\\
& & 2 & 46 $\times$ 91, 88 $\times$ 56 & 10 & 5 &  & 2(10 $\times$ 18), 14 $\times$ 9, 3(11 $\times$ 17), 3(17 $\times$ 9), \underline{9 $\times$ 8}\\ 
\hline	 
\multirow{8}{*}{13} & \multirow{8}{*}{8} & 2 & 58 $\times$ 43, 39 $\times$ 51 & 5 & 3 & \multirow{8}{*}{4} & 10 $\times$ 18, 3(9 $\times$ 9), 12 $\times$ 8\\
& & 3 & 94 $\times$ 47, 97 $\times$ 39, 85 $\times$ 70 & 6 & 3 &  & 8 $\times$ 8, 3(17 $\times$ 6), 2(15 $\times$ 6)\\
& & 2 & 84 $\times$ 72, 85 $\times$ 77 & 6 & 3 &  & 13 $\times$ 18, 3(17 $\times$ 6), 2(5 $\times$ 13)\\
& & 3 & 83 $\times$ 81, 55 $\times$ 67, 81 $\times$ 86 & 7 & 4 &  & 12 $\times$ 12, 3(\underline{13 $\times$ 5}), 15 $\times$ 11, 2(\underline{5 $\times$ 9})\\
& & 3 & 51 $\times$ 61, 97 $\times$ 53, 41 $\times$ 46 & 2 & 2 &  & 18 $\times$ 14, \underline{6 $\times$ 8}\\
& & 2 & 62 $\times$ 45, 60 $\times$ 75 & 3 & 2 &  & 2(6 $\times$ 19), 6 $\times$ 16\\
& & 3 & 44 $\times$ 91, 70 $\times$ 99, 30 $\times$ 51 & 3 & 2 &  & 10 $\times$ 9, 2(\underline{11 $\times$ 7})\\
& & 3 & 96 $\times$ 85, 41 $\times$ 59, 98 $\times$ 73 & 5 & 2 &  & 3(18 $\times$ 9), 2(20 $\times$ 8)\\
\hline
\multirow{8}{*}{14} & \multirow{8}{*}{8} & 3 & 33 $\times$ 32, 57 $\times$ 91, 62 $\times$ 84 & 4 & 2 & \multirow{8}{*}{1} & 3(12 $\times$ 13), 10 $\times$ 10\\
& & 2 & 91 $\times$ 83, 81 $\times$ 68 & 4 & 3 &  & 2(16 $\times$ 18), 16 $\times$ 7, 15 $\times$ 8\\
& & 2 & 70 $\times$ 35, 39 $\times$ 72 & 7 & 4 &  & 8 $\times$ 19, 2(10 $\times$ 10), 3(6 $\times$ 16), 10 $\times$ 6\\
& & 2 & 78 $\times$ 92, 51 $\times$ 93 & 5 & 3 &  & 20 $\times$ 14, 3(15 $\times$ 8), 16 $\times$ 17\\
& & 3 & 50 $\times$ 70, 71 $\times$ 81, 33 $\times$ 47 & 10 & 5 &  & 2(18 $\times$ 5), 13 $\times$ 15, 2(15 $\times$ 5), 3(17 $\times$ 5), 2(15 $\times$ 17)\\
& & 3 & 57 $\times$ 50, 34 $\times$ 86, 94 $\times$ 45 & 8 & 4 &  & 3(19 $\times$ 16), 3(18 $\times$ 12), 14 $\times$ 14, 14 $\times$ 17\\
& & 3 & 68 $\times$ 94, 50 $\times$ 68, 48 $\times$ 53 & 11 & 5 &  & 3(\underline{5 $\times$ 5}), 3(8 $\times$ 16), 3(14 $\times$ 12), 16 $\times$ 20, 11 $\times$ 6\\
& & 2 & 61 $\times$ 64, 73 $\times$ 89 & 6 & 3 &  & 3(7 $\times$ 6), 2(12 $\times$ 15), 16 $\times$ 5\\ \hline
\multirow{8}{*}{15} & \multirow{8}{*}{8} & 2 & 85 $\times$ 40, 55 $\times$ 36 & 5 & 2 & \multirow{8}{*}{4} & 3(17 $\times$ 13), 2(\underline{8 $\times$ 11})\\
& & 3 & 59 $\times$ 53, 92 $\times$ 88, 51 $\times$ 58 & 10 & 5 &  & 2(18 $\times$ 12), 3(\underline{7 $\times$ 18}), 3(11 $\times$ 17), 13 $\times$ 10, \underline{8 $\times$ 11}\\
& & 3 & 98 $\times$ 82, 2(44 $\times$ 49) & 6 & 4 &  & 16 $\times$ 6, 3(18 $\times$ 17), 19 $\times$ 19, 19 $\times$ 16\\
& & 2 & 51 $\times$ 89, 32 $\times$ 70 & 4 & 2 &  & 3(10 $\times$ 17), 13 $\times$ 20\\
& & 4 & 35 $\times$ 51, 38 $\times$ 80, 2(31 $\times$ 49) & 7 & 3 &  & 3(18 $\times$ 14), 2(15 $\times$ 8), 2(\underline{13 $\times$ 5})\\
& & 3 & 67 $\times$ 77, 37 $\times$ 55, 39 $\times$ 78 & 8 & 3 &  & 2(17 $\times$ 6), 3(\underline{10 $\times$ 6}), 3(16 $\times$ 17)\\
& & 2 & 88 $\times$ 70, 54 $\times$ 83 & 11 & 5 &  & 8 $\times$ 20, 2(11 $\times$ 11), 3(11 $\times$ 16), 2(15 $\times$ 10), 3(20 $\times$ 17)\\
& & 2 & 57 $\times$ 83, 45 $\times$ 66 & 6 & 4 &  & 14 $\times$ 14, 3(16 $\times$ 10), 14 $\times$ 20, 10 $\times$ 7\\ 
\hline
\multirow{8}{*}{16} & \multirow{8}{*}{8} & 5 & 31 $\times$ 98, 2(51 $\times$ 39), 2(30 $\times$ 64) & 5 & 3 & \multirow{8}{*}{2} & 2(20 $\times$ 20), 2(20 $\times$ 17), 18 $\times$ 14\\
& & 2 & 86 $\times$ 87, 82 $\times$ 98 & 4 & 2 &  & 3(7 $\times$ 9), 13 $\times$ 5\\
& & 3 & 68 $\times$ 97, 65 $\times$ 65, 78 $\times$ 34 & 10 & 5 &  & 2(17 $\times$ 13), 3(16 $\times$ 12), 12 $\times$ 11, 6 $\times$ 17, 3(\underline{7 $\times$ 5})\\
& & 2 & 54 $\times$ 85, 53 $\times$ 59 & 4 & 2 &  & 12 $\times$ 6, 3(7 $\times$ 11)\\
& & 2 & 43 $\times$ 64, 35 $\times$ 85 & 9 & 3 &  & 3(14 $\times$ 9), 3(16 $\times$ 17), 3(15 $\times$ 18)\\
& & 3 & 82 $\times$ 99, 38 $\times$ 98, 52 $\times$ 53 & 13 & 5 &  & 3(\underline{7 $\times$ 5}), 3(9 $\times$ 10), 3(15 $\times$ 7), 13 $\times$ 10, 3(\underline{6 $\times$ 6})\\
& & 4 & 66 $\times$ 47, 3(35 $\times$ 41) & 6 & 3 &  & 20 $\times$ 7, 2(19 $\times$ 12), 3(20 $\times$ 18)\\
& & 2 & 73 $\times$ 50, 38 $\times$ 84 & 3 & 2 &  & 14 $\times$ 19, 2(17 $\times$ 11)\\ 
\hline
\multirow{8}{*}{17} & \multirow{8}{*}{8} & 2 & 81 $\times$ 37, 33 $\times$ 64 & 6 & 3 & \multirow{8}{*}{2} & 2(9 $\times$ 15), 19 $\times$ 18, 3(11 $\times$ 14)\\
& & 3 & 34 $\times$ 83, 59 $\times$ 86, 72 $\times$ 44 & 5 & 3 &  & 20 $\times$ 15, 14 $\times$ 10, 3(18 $\times$ 14)\\
& & 2 & 55 $\times$ 91, 32 $\times$ 43 & 8 & 4 &  & 17 $\times$ 7, 2(14 $\times$ 20), 2(8 $\times$ 7), 3(8 $\times$ 18)\\
& & 2 & 41 $\times$ 96, 41 $\times$ 86 & 7 & 5 &  & 2(9 $\times$ 9), 18 $\times$ 7, 15 $\times$ 16, 17 $\times$ 18, 2(8 $\times$ 15)\\
& & 2 & 80 $\times$ 86, 74 $\times$ 59 & 11 & 4 &  & 3(14 $\times$ 14), 3(6 $\times$ 20), 3(19 $\times$ 8), 2(11 $\times$ 12)\\
& & 4 & 85 $\times$ 39, 85 $\times$ 63, 2(51 $\times$ 35) & 10 & 4 &  & 2(20 $\times$ 16), 3(14 $\times$ 10), 2(18 $\times$ 20), 3(8 $\times$ 17)\\
& & 2 & 78 $\times$ 53, 62 $\times$ 93 & 9 & 5 &  & 3(20 $\times$ 16), 2(11 $\times$ 5), 2(15 $\times$ 12), 14 $\times$ 14, 9 $\times$ 14\\
& & 2 & 56 $\times$ 66, 52 $\times$ 85 & 15 & 5 &  & 3(\underline{6 $\times$ 8}), 3(\underline{8 $\times$ 5}), 3(11 $\times$ 17), 3(12 $\times$ 16), 3(20 $\times$ 6)\\ 
\hline
\multicolumn{8}{|r|}{{Continued on next page}}\\
\hline
\end{tabular}}
\end{center}
\end{table}

\setcounter{table}{9}
\begin{table}[ht!]
\caption{-- continued from previous page}
\label{tab7part3}
\begin{center}
\resizebox{\textwidth}{!}{
\begin{tabular}{|c|c|rl|rrll|}
\hline 
\multirow{2}{*}{Inst.} & \multirow{2}{*}{$P$} & \multicolumn{ 2}{c|}{Objects} & \multicolumn{ 4}{c|}{Items} \\ 
\cline{ 3- 8}
& & \multicolumn{1}{l}{$m^s$} & $W_j^s \times H_j^s$ & \multicolumn{1}{l}{$n^s$} & 
\multicolumn{1}{l}{$\tilde{n}^s$} & $d$ & $w_i^s \times h_i^s$ \\ 
\hline	 
\hline
\multirow{8}{*}{18} & \multirow{8}{*}{8} & 2 & 45 $\times$ 83, 97 $\times$ 52 & 7 & 3 & \multirow{8}{*}{2} & 15 $\times$ 15, 3(11 $\times$ 13), 3(18 $\times$ 13)\\
& & 2 & 89 $\times$ 87, 88 $\times$ 45 & 8 & 5 &  & 2(18 $\times$ 9), 6 $\times$ 7, 2(12 $\times$ 8), 8 $\times$ 19, 2(\underline{18 $\times$ 6})\\
& & 3 & 2(65 $\times$ 33), 92 $\times$ 72 & 8 & 3 &  & 3(19 $\times$ 20), 2(15 $\times$ 14), 3(9 $\times$ 14)\\
& & 3 & 76 $\times$ 40, 54 $\times$ 71, 43 $\times$ 78 & 9 & 4 &  & 3(7 $\times$ 8), 5 $\times$ 17, 3(6 $\times$ 11), 2(17 $\times$ 15)\\
& & 2 & 72 $\times$ 74, 89 $\times$ 73 & 5 & 2 &  & 3(11 $\times$ 7), 2(20 $\times$ 16)\\
& & 4 & 59 $\times$ 38, 2(44 $\times$ 32), 46 $\times$ 47 & 7 & 4 &  & 6 $\times$ 17, 2(18 $\times$ 16), 2(8 $\times$ 15), 2(18 $\times$ 11)\\
& & 3 & 56 $\times$ 41, 100 $\times$ 45, 40 $\times$ 92 & 2 & 2 &  & 13 $\times$ 20, 18 $\times$ 13\\
& & 2 & 73 $\times$ 77, 83 $\times$ 54 & 6 & 3 &  & 2(\underline{5 $\times$ 7}), 2(16 $\times$ 18), 2(10 $\times$ 9)\\ 
\hline
\multirow{8}{*}{19} & \multirow{8}{*}{8} & 2 & 78 $\times$ 86, 72 $\times$ 67 & 10 & 5 & \multirow{8}{*}{2} & 3(\underline{15 $\times$ 5}), 3(\underline{6 $\times$ 6}), 18 $\times$ 10, 2(8 $\times$ 10), 14 $\times$ 19\\
& & 3 & 53 $\times$ 67, 37 $\times$ 80, 67 $\times$ 56 & 8 & 4 &  & 2(17 $\times$ 5), 2(20 $\times$ 15), 2(15 $\times$ 13), 2(15 $\times$ 9)\\
& & 3 & 57 $\times$ 85, 52 $\times$ 50, 75 $\times$ 37 & 6 & 3 &  & 2(17 $\times$ 9), 2(9 $\times$ 9), 2(12 $\times$ 14)\\
& & 3 & 64 $\times$ 44, 45 $\times$ 96, 75 $\times$ 52 & 10 & 5 &  & 3(18 $\times$ 20), 2(13 $\times$ 9), 8 $\times$ 9, 9 $\times$ 7, 3(14 $\times$ 14)\\
& & 2 & 56 $\times$ 93, 53 $\times$ 49 & 9 & 4 &  & 3(16 $\times$ 10), 3(10 $\times$ 14), 12 $\times$ 17, 2(6 $\times$ 15)\\
& & 2 & 51 $\times$ 89, 65 $\times$ 72 & 5 & 3 &  & 16 $\times$ 14, 18 $\times$ 8, 3(16 $\times$ 5)\\
& & 2 & 92 $\times$ 64, 81 $\times$ 95 & 6 & 3 &  & 3(19 $\times$ 7), 2(6 $\times$ 14), 17 $\times$ 16\\
& & 3 & 62 $\times$ 52, 32 $\times$ 97, 95 $\times$ 35 & 8 & 4 &  & 3(7 $\times$ 16), 2(10 $\times$ 14), 11 $\times$ 12, 2(13 $\times$ 8)\\ 
\hline
\multirow{8}{*}{20} & \multirow{8}{*}{8} & 3 & 75 $\times$ 82, 69 $\times$ 79, 76 $\times$ 64 & 5 & 3 & \multirow{8}{*}{3} & 2(14 $\times$ 10), 2(15 $\times$ 13), 14 $\times$ 12\\
& & 2 & 49 $\times$ 68, 61 $\times$ 79 & 12 & 5 &  & 3(11 $\times$ 18), 2(6 $\times$ 12), 2(\underline{7 $\times$ 7}), 3(\underline{5 $\times$ 12}), 2(13 $\times$ 18)\\
& & 3 & 92 $\times$ 41, 74 $\times$ 51, 78 $\times$ 93 & 7 & 4 &  & \underline{10 $\times$ 5}, 2(13 $\times$ 6), 2(8 $\times$ 10), 2(5 $\times$ 13)\\
& & 3 & 61 $\times$ 85, 45 $\times$ 51, 34 $\times$ 50 & 7 & 3 &  & 3(8 $\times$ 19), 14 $\times$ 10, 3(9 $\times$ 11)\\
& & 2 & 41 $\times$ 50, 63 $\times$ 84 & 6 & 3 &  & 2(13 $\times$ 20), 2(18 $\times$ 12), 2(\underline{10 $\times$ 5})\\
& & 2 & 81 $\times$ 43, 53 $\times$ 45 & 6 & 4 &  & 2(7 $\times$ 14), 13 $\times$ 7, 2(9 $\times$ 11), 19 $\times$ 17\\
& & 3 & 35 $\times$ 82, 2(33 $\times$ 34) & 7 & 5 &  & 6 $\times$ 14, 2(17 $\times$ 19), 19 $\times$ 10, 2(15 $\times$ 9), 11 $\times$ 11\\
& & 3 & 92 $\times$ 52, 83 $\times$ 65, 70 $\times$ 70 & 13 & 5 &  & 3(13 $\times$ 18), 2(16 $\times$ 6), 3(12 $\times$ 8), 3(5 $\times$ 18), 2(19 $\times$ 11)\\ 
\hline	 
\multirow{12}{*}{21} & \multirow{12}{*}{12} & 2 & 65 $\times$ 50, 93 $\times$ 92 & 5 & 2 & \multirow{12}{*}{2} & 2(7 $\times$ 8), 3(12 $\times$ 10)\\
& & 2 & 90 $\times$ 68, 57 $\times$ 69 & 7 & 4 &  & 2(13 $\times$ 6), 3(19 $\times$ 14), 6 $\times$ 11, \underline{6 $\times$ 5}\\
& & 3 & 78 $\times$ 71, 56 $\times$ 70, 62 $\times$ 100 & 6 & 3 &  & 19 $\times$ 15, 2(8 $\times$ 17), 3(15 $\times$ 19)\\
& & 2 & 50 $\times$ 84, 30 $\times$ 49 & 7 & 4 &  & 2(7 $\times$ 7), 14 $\times$ 17, 3(14 $\times$ 13), 8 $\times$ 16\\
& & 2 & 73 $\times$ 99, 44 $\times$ 72 & 4 & 2 &  & 7 $\times$ 13, 3(8 $\times$ 7)\\
& & 3 & 48 $\times$ 50, 70 $\times$ 79, 100 $\times$ 52 & 10 & 5 &  & 17 $\times$ 16, 2(13 $\times$ 17), 2(\underline{5 $\times$ 10}), 2(16 $\times$ 12), 3(6 $\times$ 15)\\
& & 3 & 36 $\times$ 93, 36 $\times$ 77, 92 $\times$ 90 & 4 & 2 &  & 2(13 $\times$ 15), 2(9 $\times$ 18)\\
& & 3 & 74 $\times$ 65, 47 $\times$ 70, 100 $\times$ 34 & 4 & 2 &  & 3(15 $\times$ 18), 11 $\times$ 9\\
& & 2 & 50 $\times$ 81, 70 $\times$ 87 & 5 & 3 &  & 16 $\times$ 10, 2(16 $\times$ 17), 2(10 $\times$ 13)\\
& & 2 & 52 $\times$ 86, 46 $\times$ 48 & 9 & 5 &  & 11 $\times$ 14, 2(19 $\times$ 8), 7 $\times$ 14, 2(15 $\times$ 6), 3(15 $\times$ 19)\\
& & 2 & 93 $\times$ 47, 31 $\times$ 89 & 5 & 3 &  & 13 $\times$ 16, 15 $\times$ 18, 3(18 $\times$ 7)\\
& & 2 & 81 $\times$ 92, 37 $\times$ 80 & 11 & 5 &  & 3(9 $\times$ 14), 2(16 $\times$ 8), 2(5 $\times$ 19), 15 $\times$ 7, 3(14 $\times$ 17)\\ 
\hline
\multirow{12}{*}{22} & \multirow{12}{*}{12} & 2 & 73 $\times$ 35, 72 $\times$ 91 & 5 & 4 & \multirow{12}{*}{2} & \underline{6 $\times$ 5}, 7 $\times$ 8, 16 $\times$ 12, 2(11 $\times$ 8)\\
& & 2 & 39 $\times$ 63, 54 $\times$ 63 & 8 & 3 &  & 2(13 $\times$ 13), 3(7 $\times$ 19), 3(11 $\times$ 7)\\
& & 3 & 96 $\times$ 44, 63 $\times$ 56, 54 $\times$ 53 & 5 & 4 &  & 8 $\times$ 20, 15 $\times$ 11, 18 $\times$ 8, 2(14 $\times$ 9)\\
& & 2 & 45 $\times$ 82, 69 $\times$ 37 & 12 & 5 &  & 3(17 $\times$ 17), 3(19 $\times$ 11), 13 $\times$ 11, 3(9 $\times$ 11), 2(7 $\times$ 14)\\
& & 3 & 72 $\times$ 62, 63 $\times$ 36, 37 $\times$ 97 & 5 & 4 &  & 18 $\times$ 13, 19 $\times$ 15, 2(18 $\times$ 19), 15 $\times$ 14\\
& & 2 & 39 $\times$ 37, 84 $\times$ 42 & 6 & 2 &  & 3(17 $\times$ 6), 3(10 $\times$ 5)\\
& & 3 & 2(31 $\times$ 38), 98 $\times$ 38 & 13 & 5 &  & 2(8 $\times$ 18), 3(8 $\times$ 16), 3(6 $\times$ 13), 2(16 $\times$ 7), 3(8 $\times$ 7)\\
& & 3 & 99 $\times$ 67, 94 $\times$ 93, 65 $\times$ 87 & 12 & 5 &  & 3(14 $\times$ 6), 20 $\times$ 19, 2(20 $\times$ 14), 3(17 $\times$ 17), 3(12 $\times$ 14)\\
& & 2 & 78 $\times$ 66, 42 $\times$ 95 & 6 & 3 &  & 3(9 $\times$ 5), 18 $\times$ 13, 2(\underline{6 $\times$ 5})\\
& & 2 & 78 $\times$ 50, 84 $\times$ 44 & 6 & 3 &  & 3(13 $\times$ 13), 12 $\times$ 9, 2(15 $\times$ 16)\\
& & 3 & 76 $\times$ 51, 70 $\times$ 88, 76 $\times$ 57 & 6 & 2 &  & 3(15 $\times$ 12), 3(7 $\times$ 12)\\
& & 3 & 71 $\times$ 40, 44 $\times$ 52, 55 $\times$ 58 & 6 & 3 &  & \underline{5 $\times$ 18}, 3(12 $\times$ 6), 2(6 $\times$ 17)\\
\hline
\multicolumn{8}{|r|}{{Continued on next page}}\\
\hline
\end{tabular}}
\end{center}
\end{table}

\setcounter{table}{9}
\begin{table}[ht!]
\caption{-- continued from previous page}
\label{tab7part4}
\begin{center}
\resizebox{\textwidth}{!}{
\begin{tabular}{|c|c|rl|rrll|}
\hline 
\multirow{2}{*}{Inst.} & \multirow{2}{*}{$P$} & \multicolumn{ 2}{c|}{Objects} & \multicolumn{ 4}{c|}{Items} \\ 
\cline{ 3- 8}
& & \multicolumn{1}{l}{$m^s$} & $W_j^s \times H_j^s$ & \multicolumn{1}{l}{$n^s$} & 
\multicolumn{1}{l}{$\tilde{n}^s$} & $d$ & $w_i^s \times h_i^s$ \\ 
\hline	 
\hline
\multirow{12}{*}{23} & \multirow{12}{*}{12} & 3 & 100 $\times$ 62, 68 $\times$ 83, 86 $\times$ 66 & 4 & 2 & \multirow{12}{*}{2} & 3(5 $\times$ 11), 20 $\times$ 15\\
& & 2 & 82 $\times$ 51, 65 $\times$ 68 & 8 & 5 &  & 2(8 $\times$ 19), 2(20 $\times$ 18), 19 $\times$ 11, 14 $\times$ 7, 2(19 $\times$ 5)\\
& & 2 & 66 $\times$ 60, 60 $\times$ 63 & 4 & 2 &  & \underline{12 $\times$ 5}, 3(17 $\times$ 14)\\
& & 3 & 81 $\times$ 52, 32 $\times$ 97, 97 $\times$ 46 & 7 & 3 &  & 2(20 $\times$ 10), 3(11 $\times$ 10), 2(13 $\times$ 18)\\
& & 2 & 34 $\times$ 57, 39 $\times$ 95 & 6 & 4 &  & 2(13 $\times$ 18), 2(13 $\times$ 15), 6 $\times$ 12, 20 $\times$ 17\\
& & 2 & 38 $\times$ 92, 33 $\times$ 95 & 6 & 4 &  & 2(19 $\times$ 9), 11 $\times$ 17, 2(17 $\times$ 9), 17 $\times$ 17\\
& & 3 & 77 $\times$ 44, 37 $\times$ 100, 50 $\times$ 37 & 9 & 3 &  & 3(9 $\times$ 16), 3(5 $\times$ 20), 3(19 $\times$ 9)\\
& & 3 & 86 $\times$ 62, 92 $\times$ 99, 72 $\times$ 43 & 5 & 2 &  & 2(19 $\times$ 5), 3(15 $\times$ 17)\\
& & 2 & 58 $\times$ 34, 57 $\times$ 88 & 7 & 4 &  & 10 $\times$ 17, 6 $\times$ 15, 2(5 $\times$ 12), 3(10 $\times$ 10)\\
& & 3 & 2(51 $\times$ 45), 50 $\times$ 53 & 9 & 5 &  & 3(19 $\times$ 6), 9 $\times$ 16, \underline{5 $\times$ 8}, 3(20 $\times$ 20), 15 $\times$ 10\\
& & 3 & 98 $\times$ 92, 84 $\times$ 46, 35 $\times$ 45 & 6 & 3 &  & 11 $\times$ 20, 2(12 $\times$ 15), 3(15 $\times$ 6)\\
& & 2 & 37 $\times$ 35, 41 $\times$ 54 & 2 & 2 &  & 14 $\times$ 6, 14 $\times$ 9\\ 
\hline
\multirow{12}{*}{24} & \multirow{12}{*}{12} & 3 & 69 $\times$ 73, 63 $\times$ 95, 62 $\times$ 94 & 8 & 4 & \multirow{12}{*}{2} & 19 $\times$ 20, 3(13 $\times$ 12), 14 $\times$ 7, 3(14 $\times$ 19)\\
& & 2 & 69 $\times$ 32, 39 $\times$ 59 & 8 & 4 &  & 8 $\times$ 9, 2(10 $\times$ 8), 3(18 $\times$ 14), 2(10 $\times$ 19)\\
& & 3 & 97 $\times$ 33, 78 $\times$ 42, 56 $\times$ 30 & 7 & 3 &  & 17 $\times$ 14, 3(15 $\times$ 10), 3(20 $\times$ 12)\\
& & 3 & 87 $\times$ 55, 36 $\times$ 76, 33 $\times$ 56 & 4 & 2 &  & 3(10 $\times$ 6), 15 $\times$ 20\\
& & 3 & 100 $\times$ 84, 2(36 $\times$ 41) & 10 & 5 &  & 15 $\times$ 18, 3(8 $\times$ 8), 2(13 $\times$ 16), 20 $\times$ 15, 3(15 $\times$ 17)\\
& & 3 & 85 $\times$ 67, 92 $\times$ 35, 46 $\times$ 98 & 5 & 3 &  & 8 $\times$ 19, 19 $\times$ 6, 3(19 $\times$ 19)\\
& & 2 & 52 $\times$ 75, 56 $\times$ 60 & 10 & 4 &  & 3(14 $\times$ 18), 3(\underline{8 $\times$ 6}), \underline{5 $\times$ 15}, 3(9 $\times$ 17)\\
& & 3 & 35 $\times$ 53, 67 $\times$ 54, 62 $\times$ 93 & 4 & 2 &  & 11 $\times$ 7, 3(9 $\times$ 7)\\
& & 2 & 97 $\times$ 66, 69 $\times$ 39 & 4 & 3 &  & 7 $\times$ 18, 8 $\times$ 8, 2(19 $\times$ 17)\\
& & 2 & 83 $\times$ 38, 54 $\times$ 66 & 7 & 3 &  & 2(18 $\times$ 7), 3(20 $\times$ 13), 2(19 $\times$ 17)\\
& & 2 & 87 $\times$ 51, 33 $\times$ 55 & 4 & 2 &  & 2(9 $\times$ 20), 2(15 $\times$ 7)\\
& & 3 & 68 $\times$ 68, 39 $\times$ 87, 82 $\times$ 78 & 6 & 3 &  & 19 $\times$ 14, 2(5 $\times$ 18), 3(13 $\times$ 8)\\ 
\hline
\multirow{12}{*}{25} & \multirow{12}{*}{12} & 3 & 86 $\times$ 45, 57 $\times$ 40, 64 $\times$ 87 & 9 & 4 & \multirow{12}{*}{1} & 15 $\times$ 11, 3(14 $\times$ 20), 3(9 $\times$ 16), 2(15 $\times$ 7)\\
& & 2 & 70 $\times$ 31, 95 $\times$ 99 & 8 & 5 &  & 7 $\times$ 6, 2(12 $\times$ 20), 19 $\times$ 8, 3(15 $\times$ 8), 7 $\times$ 18\\
& & 3 & 49 $\times$ 36, 83 $\times$ 98, 35 $\times$ 51 & 4 & 2 &  & 2(10 $\times$ 16), 2(20 $\times$ 12)\\
& & 4 & 61 $\times$ 63, 97 $\times$ 89, 2(34 $\times$ 40) & 12 & 5 &  & 20 $\times$ 15, 3(14 $\times$ 18), 3(16 $\times$ 15), 3(9 $\times$ 6), 2(8 $\times$ 16)\\
& & 3 & 33 $\times$ 65, 68 $\times$ 56, 90 $\times$ 82 & 10 & 4 &  & 3(12 $\times$ 11), 3(20 $\times$ 13), 12 $\times$ 20, 3(6 $\times$ 13)\\
& & 2 & 83 $\times$ 83, 79 $\times$ 81 & 5 & 3 &  & 3(15 $\times$ 19), 11 $\times$ 14, 11 $\times$ 15\\
& & 2 & 51 $\times$ 77, 33 $\times$ 95 & 6 & 3 &  & 2(\underline{5 $\times$ 5}), 2(7 $\times$ 12), 2(8 $\times$ 14)\\
& & 2 & 32 $\times$ 35, 99 $\times$ 81 & 6 & 3 &  & 2(17 $\times$ 17), 3(14 $\times$ 7), 7 $\times$ 13\\
& & 3 & 47 $\times$ 58, 72 $\times$ 81, 83 $\times$ 51 & 2 & 2 &  & 14 $\times$ 6, 5 $\times$ 17\\
& & 3 & 42 $\times$ 99, 75 $\times$ 47, 57 $\times$ 87 & 10 & 5 &  & 2(6 $\times$ 20), 2(15 $\times$ 6), 3(17 $\times$ 14), 19 $\times$ 14, 2(19 $\times$ 12)\\
& & 2 & 66 $\times$ 59, 54 $\times$ 86 & 4 & 2 &  & 5 $\times$ 18, 3(5 $\times$ 20)\\
& & 3 & 55 $\times$ 58, 99 $\times$ 45, 67 $\times$ 73 & 6 & 3 &  & 2(11 $\times$ 15), 3(20 $\times$ 13), 13 $\times$ 19\\
\hline	 
\multirow{12}{*}{26} & \multirow{12}{*}{12} & 2 & 51 $\times$ 42, 79 $\times$ 85 & 5 & 4 & \multirow{12}{*}{5} & 6 $\times$ 13, 8 $\times$ 15, 2(16 $\times$ 7), 15 $\times$ 15\\
& & 3 & 95 $\times$ 82, 100 $\times$ 90, 54 $\times$ 75 & 3 & 2 &  & 2(\underline{18 $\times$ 5}), 7 $\times$ 17\\
& & 2 & 85 $\times$ 35, 69 $\times$ 83 & 4 & 2 &  & 7 $\times$ 19, 3(17 $\times$ 13)\\
& & 2 & 90 $\times$ 100, 81 $\times$ 96 & 11 & 5 &  & 2(13 $\times$ 12), 2(12 $\times$ 19), 2(20 $\times$ 17), 2(16 $\times$ 19), 3(14 $\times$ 6)\\
& & 3 & 79 $\times$ 91, 51 $\times$ 40, 85 $\times$ 79 & 8 & 5 &  & 13 $\times$ 15, 19 $\times$ 7, 2(14 $\times$ 15), 2(6 $\times$ 19), 2(20 $\times$ 7)\\
& & 3 & 78 $\times$ 59, 85 $\times$ 31, 85 $\times$ 56 & 10 & 5 &  & 2(17 $\times$ 11), 3(\underline{10 $\times$ 9}), \underline{5 $\times$ 19}, 3(15 $\times$ 11), 18 $\times$ 12\\
& & 2 & 81 $\times$ 76, 66 $\times$ 70 & 5 & 3 &  & 2(\underline{12 $\times$ 6}), 2(19 $\times$ 16), 11 $\times$ 20\\
& & 2 & 80 $\times$ 52, 74 $\times$ 68 & 3 & 3 &  & 14 $\times$ 6, 14 $\times$ 17, 13 $\times$ 14\\
& & 3 & 83 $\times$ 95, 45 $\times$ 48, 95 $\times$ 63 & 5 & 3 &  & 7 $\times$ 10, 3(19 $\times$ 8), 18 $\times$ 16\\
& & 2 & 79 $\times$ 82, 79 $\times$ 36 & 7 & 3 &  & 2(17 $\times$ 19), 2(13 $\times$ 11), 3(\underline{6 $\times$ 10})\\
& & 3 & 32 $\times$ 85, 45 $\times$ 97, 78 $\times$ 86 & 8 & 4 &  & 2(14 $\times$ 18), 3(17 $\times$ 19), 2(12 $\times$ 15), 7 $\times$ 13\\
& & 2 & 45 $\times$ 42, 36 $\times$ 71 & 7 & 3 &  & 9 $\times$ 15, 3(14 $\times$ 8), 3(19 $\times$ 10)\\ 
\hline
\multicolumn{8}{|r|}{{Continued on next page}}\\
\hline
\end{tabular}}
\end{center}
\end{table}

\setcounter{table}{9}
\begin{table}[ht!]
\caption{-- continued from previous page}
\label{tab7part5}
\begin{center}
\resizebox{\textwidth}{!}{
\begin{tabular}{|c|c|rl|rrll|}
\hline 
\multirow{2}{*}{Inst.} & \multirow{2}{*}{$P$} & \multicolumn{ 2}{c|}{Objects} & \multicolumn{ 4}{c|}{Items} \\ 
\cline{ 3- 8}
& & \multicolumn{1}{l}{$m^s$} & $W_j^s \times H_j^s$ & \multicolumn{1}{l}{$n^s$} & 
\multicolumn{1}{l}{$\tilde{n}^s$} & $d$ & $w_i^s \times h_i^s$ \\ 
\hline	 
\hline
\multirow{12}{*}{27} & \multirow{12}{*}{12} & 5 & 47 $\times$ 71, 71 $\times$ 96, 3(32 $\times$ 51) & 10 & 4 & \multirow{12}{*}{1} & 16 $\times$ 9, 3(19 $\times$ 13), 3(17 $\times$ 12), 3(18 $\times$ 17)\\
& & 2 & 62 $\times$ 65, 38 $\times$ 91 & 3 & 2 &  & 2(20 $\times$ 18), 11 $\times$ 5\\
& & 2 & 100 $\times$ 62, 69 $\times$ 62 & 7 & 3 &  & 18 $\times$ 5, 3(13 $\times$ 19), 3(17 $\times$ 15)\\
& & 2 & 61 $\times$ 47, 84 $\times$ 91 & 11 & 5 &  & 6 $\times$ 6, 3(20 $\times$ 5), 15 $\times$ 12, 3(17 $\times$ 18), 3(7 $\times$ 15)\\
& & 3 & 90 $\times$ 82, 42 $\times$ 52, 91 $\times$ 35 & 12 & 5 &  & 3(13 $\times$ 13), 5 $\times$ 18, 3(8 $\times$ 8), 2(9 $\times$ 15), 3(10 $\times$ 18)\\
& & 2 & 93 $\times$ 96, 95 $\times$ 54 & 11 & 5 &  & 2(8 $\times$ 15), 2(16 $\times$ 15), 15 $\times$ 13, 3(11 $\times$ 5), 3(10 $\times$ 5)\\
& & 2 & 67 $\times$ 97, 72 $\times$ 65 & 5 & 2 &  & 3(9 $\times$ 18), 2(14 $\times$ 14)\\
& & 2 & 43 $\times$ 81, 58 $\times$ 100 & 5 & 4 &  & 2(11 $\times$ 6), 18 $\times$ 17, 9 $\times$ 7, 8 $\times$ 13\\
& & 3 & 37 $\times$ 58, 48 $\times$ 40, 54 $\times$ 93 & 4 & 2 &  & 16 $\times$ 20, 3(10 $\times$ 13)\\
& & 3 & 63 $\times$ 69, 71 $\times$ 52, 50 $\times$ 36 & 4 & 2 &  & 2(15 $\times$ 17), 2(19 $\times$ 19)\\
& & 2 & 89 $\times$ 50, 94 $\times$ 56 & 8 & 4 &  & 3(\underline{5 $\times$ 5}), 14 $\times$ 11, 13 $\times$ 11, 3(5 $\times$ 20)\\
& & 2 & 91 $\times$ 67, 57 $\times$ 72 & 7 & 4 &  & 15 $\times$ 17, 18 $\times$ 16, 2(7 $\times$ 18), 3(13 $\times$ 19)\\ 
\hline
\multirow{12}{*}{28} & \multirow{12}{*}{12} & 2 & 93 $\times$ 73, 38 $\times$ 66 & 9 & 4 & \multirow{12}{*}{3} & 3(15 $\times$ 13), 13 $\times$ 11, 3(15 $\times$ 5), 2(8 $\times$ 15)\\
& & 3 & 94 $\times$ 36, 53 $\times$ 41, 100 $\times$ 64 & 5 & 2 &  & 2(20 $\times$ 16), 3(6 $\times$ 12)\\
& & 2 & 69 $\times$ 98, 92 $\times$ 99 & 8 & 3 &  & 2(17 $\times$ 19), 3(8 $\times$ 10), 3(8 $\times$ 17)\\
& & 3 & 75 $\times$ 42, 36 $\times$ 41, 66 $\times$ 47 & 3 & 2 &  & 2(19 $\times$ 12), 14 $\times$ 17\\
& & 3 & 2(35 $\times$ 40), 59 $\times$ 64 & 9 & 5 &  & 19 $\times$ 11, 17 $\times$ 11, 6 $\times$ 20, 3(18 $\times$ 17), 3(11 $\times$ 6)\\
& & 2 & 71 $\times$ 51, 53 $\times$ 31 & 6 & 2 &  & 3(19 $\times$ 14), 3(15 $\times$ 15)\\
& & 2 & 73 $\times$ 55, 71 $\times$ 61 & 6 & 3 &  & 2(14 $\times$ 18), 2(\underline{5 $\times$ 19}), 2(15 $\times$ 16)\\
& & 2 & 93 $\times$ 34, 35 $\times$ 74 & 5 & 3 &  & 2(12 $\times$ 17), 9 $\times$ 15, 2(19 $\times$ 9)\\
& & 3 & 99 $\times$ 49, 2(37 $\times$ 69) & 14 & 5 &  & 3(14 $\times$ 5), 2(\underline{7 $\times$ 5}), 3(15 $\times$ 15), 3(19 $\times$ 18), 3(9 $\times$ 19)\\
& & 2 & 65 $\times$ 81, 31 $\times$ 61 & 12 & 4 &  & 3(11 $\times$ 13), 3(7 $\times$ 8), 3(6 $\times$ 15), 3(\underline{6 $\times$ 9})\\
& & 2 & 79 $\times$ 48, 75 $\times$ 73 & 4 & 2 &  & 20 $\times$ 19, 3(12 $\times$ 7)\\
& & 2 & 89 $\times$ 72, 58 $\times$ 91 & 12 & 5 &  & 2(15 $\times$ 14), 2(10 $\times$ 17), 2(7 $\times$ 18), 3(11 $\times$ 20), 3(15 $\times$ 18)\\ 
\hline
\multirow{12}{*}{29} & \multirow{12}{*}{12} & 3 & 70 $\times$ 66, 90 $\times$ 86, 36 $\times$ 44 & 7 & 5 & \multirow{12}{*}{1} & 12 $\times$ 20, 2(8 $\times$ 20), 15 $\times$ 16, 2(9 $\times$ 6), 12 $\times$ 9\\
& & 3 & 75 $\times$ 85, 47 $\times$ 59, 32 $\times$ 38 & 6 & 4 &  & 14 $\times$ 19, 8 $\times$ 11, 7 $\times$ 10, 3(6 $\times$ 5)\\
& & 3 & 99 $\times$ 44, 45 $\times$ 83, 65 $\times$ 95 & 5 & 3 &  & 10 $\times$ 6, 15 $\times$ 20, 3(16 $\times$ 10)\\
& & 3 & 86 $\times$ 72, 48 $\times$ 81, 72 $\times$ 42 & 4 & 4 &  & 9 $\times$ 12, 10 $\times$ 12, 11 $\times$ 14, 7 $\times$ 14\\
& & 2 & 99 $\times$ 35, 48 $\times$ 43 & 6 & 3 &  & \underline{5 $\times$ 5}, 2(10 $\times$ 11), 3(6 $\times$ 10)\\
& & 3 & 39 $\times$ 43, 72 $\times$ 55, 52 $\times$ 60 & 6 & 4 &  & 2(18 $\times$ 12), 2(11 $\times$ 6), 5 $\times$ 15, 9 $\times$ 13\\
& & 2 & 30 $\times$ 34, 81 $\times$ 84 & 4 & 2 &  & 17 $\times$ 10, 3(6 $\times$ 7)\\
& & 3 & 81 $\times$ 48, 46 $\times$ 32, 38 $\times$ 36 & 9 & 5 &  & 9 $\times$ 15, 11 $\times$ 9, 3(5 $\times$ 18), 2(13 $\times$ 12), 2(13 $\times$ 6)\\
& & 3 & 89 $\times$ 65, 99 $\times$ 66, 46 $\times$ 66 & 6 & 5 &  & 5 $\times$ 9, 2(8 $\times$ 16), 11 $\times$ 5, 6 $\times$ 16, 10 $\times$ 11\\
& & 3 & 40 $\times$ 92, 46 $\times$ 49, 70 $\times$ 67 & 8 & 4 &  & 19 $\times$ 15, 20 $\times$ 15, 3(8 $\times$ 17), 3(12 $\times$ 10)\\
& & 3 & 76 $\times$ 42, 66 $\times$ 90, 85 $\times$ 60 & 10 & 4 &  & 2(9 $\times$ 9), 3(11 $\times$ 14), 3(20 $\times$ 9), 2(14 $\times$ 14)\\
& & 5 & 91 $\times$ 86, 2(46 $\times$ 39), 2(41 $\times$ 41) & 11 & 4 &  & 3(16 $\times$ 20), 2(19 $\times$ 16), 3(6 $\times$ 7), 3(20 $\times$ 15)\\ 
\hline
\multirow{12}{*}{30} & \multirow{12}{*}{12} & 3 & 34 $\times$ 50, 34 $\times$ 38, 98 $\times$ 33 & 3 & 2 & \multirow{12}{*}{2} & 2(6 $\times$ 7), 16 $\times$ 8\\
& & 3 & 49 $\times$ 78, 53 $\times$ 70, 84 $\times$ 100 & 2 & 2 &  & 8 $\times$ 19, 9 $\times$ 14\\
& & 3 & 79 $\times$ 96, 69 $\times$ 43, 76 $\times$ 73 & 8 & 5 &  & 20 $\times$ 5, 3(\underline{5 $\times$ 7}), 17 $\times$ 10, 2(12 $\times$ 12), 5 $\times$ 13\\
& & 2 & 50 $\times$ 98, 60 $\times$ 59 & 9 & 4 &  & 2(5 $\times$ 8), 3(20 $\times$ 13), 2(18 $\times$ 16), 2(13 $\times$ 15)\\
& & 3 & 36 $\times$ 100, 90 $\times$ 41, 73 $\times$ 97 & 5 & 4 &  & 8 $\times$ 15, 16 $\times$ 19, 2(17 $\times$ 11), 7 $\times$ 7\\
& & 3 & 82 $\times$ 96, 51 $\times$ 40, 55 $\times$ 47 & 6 & 3 &  & 3(9 $\times$ 8), 20 $\times$ 18, 2(10 $\times$ 9)\\
& & 3 & 50 $\times$ 78, 77 $\times$ 35, 66 $\times$ 79 & 4 & 2 &  & 3(9 $\times$ 7), 11 $\times$ 10\\
& & 2 & 44 $\times$ 45, 76 $\times$ 54 & 11 & 5 &  & 8 $\times$ 17, 3(11 $\times$ 7), 3(8 $\times$ 20), 12 $\times$ 14, 3(14 $\times$ 11)\\
& & 3 & 62 $\times$ 71, 93 $\times$ 67, 90 $\times$ 93 & 4 & 2 &  & 15 $\times$ 13, 3(15 $\times$ 15)\\
& & 3 & 89 $\times$ 62, 75 $\times$ 86, 63 $\times$ 40 & 3 & 2 &  & 17 $\times$ 9, 2(8 $\times$ 18)\\
& & 3 & 38 $\times$ 59, 59 $\times$ 71, 100 $\times$ 51 & 4 & 2 &  & 15 $\times$ 13, 3(\underline{10 $\times$ 5})\\
& & 5 & 35 $\times$ 99, 2(46 $\times$ 94), 2(61 $\times$ 51) & 10 & 4 &  & 3(19 $\times$ 16), 4(15 $\times$ 20), 3(18 $\times$ 17)\\ 
\hline
\end{tabular}}
\end{center}
\end{table}

\end{document}